\documentclass[12pt]{report}
\usepackage{german}
\usepackage{makeidx}
\usepackage{amssymb}
\usepackage{stmaryrd}
\usepackage{wasysym}
\usepackage{fontenc}
\makeindex
\author{Andreas Johann Raab\\
Luisenstrasse 60, 80798 M"unchen,\\
Federal Republic of Germany\\
E-mail:andreas@andreasjohannraab.de}

\title{Ein gesuchter, dennoch bislang unbekannter elementarer Satz}
\begin{document}

\maketitle
\newpage
\begin{center}
Ein Beitrag zur Theorie kontinuierlicher dynamischer 
Systeme,\newline
insbesondere zur Ergodenfrage derselben.
 
\end{center}\newpage

\newpage
{\bf The Subject}
\newline\newline
In this tractatus we present the elementary theorem 1.1, which we call the 
elementary quasiergodic theorem. 
We call it elementary, because it is limited to finit-dimensional real phase-spaces and
we call it the elementary quasiergodic theorem, because it answers the following question:
Given a continuous dynamical system, which has piecewise differentiable trajectories $\tau\subset\zeta$
in a finit-dimensional real phase-space $\zeta\subset\mathbb{R}^{n}$, which is compact, let its 
trajectories define a field of normed tangents, which is continuous in almost each point of the phase-space 
$\zeta$:
Which are the sets $\chi\subset\zeta$ of the phase-space, which are quasiergodic?\newline
The set of all trajectories or fixed points of this continuous dynamical system forms a partition $\gamma$ of 
its phase-space $\zeta$.
A set $\chi\subset\zeta$ is
quasiergodic, if and only if the following implication is true:
$$x,y\in \chi\ \land\ x\in\tau\in\gamma\ \Rightarrow y\in \mathbf{cl}(\tau)\ .$$
Thus we regard both, closed trajectories and fixed points, as quasiergodic,
though in the generic case a quasiergodic set turns out to be a multidimensional manifold. 
Let $\tau_{1}\subset\zeta$ and $\tau_{2}\subset\zeta$ be two of those trajectories. We will show, that
if in this case 
$\tau_{1}$ and the topological hull $\mathbf{cl}(\tau_{2})$ are not disjoint,
then $$\mathbf{cl}(\tau_{2})=\mathbf{cl}(\tau_{1})$$ 
coincides.
Let us itemize the premisses of the  elementary quasiergodic theorem:
\begin{enumerate}
\item The phase-space $\zeta\subset\mathbb{R}^{n}$ is bounded. 
\item Each trajektory $\tau$ representing a single development
is piecewise differentiable.
\item The field of normed tangents of the trajektories is continuous in almost each point of the phase-space 
$\zeta$.
If and only if the field of normed tangents of the trajektories
is continuous in almost each point of the phase-space, we say, that
the field of normed tangents is locally parallel
in almost each point of the phase-space. 
\end{enumerate}
Whereat this form of continuity in almost each point of the phase-space stands for the circumstance, that
the Lebesgue-measure of the set of
points 
of discontinuity of the field of normed tangents
is zero. If this three itemized premisses are true and if $\chi$ is one of the quasiergodic sets, then the equivalences
\begin{displaymath}
\tau_{1}\cap\chi\not=\emptyset \wedge \tau_{2}\cap\chi\not=\emptyset\ \Leftrightarrow 
\mathbf{cl}(\tau_{1})=\mathbf{cl}(\tau_{2})
\end{displaymath}
and
\begin{displaymath}
\chi=\mathbf{cl}(\tau)\ \Leftrightarrow\ \tau\cap\chi\not=\emptyset
\end{displaymath} 
are true. 
We will show, that if this three itemized premisses are true, then there exists a partition 
$$\Gamma=\{\mathbf{cl}(\tau):\tau\in\gamma\}$$ of the phase-space $\zeta\subset\mathbb{R}^{n}$.
This is the statement of the theorem of the existence of zimmers 2.1.2 and it is the
main statement of the elementary quasiergodic theorem 1.1.
Each set $\chi\in \Gamma$ 
we call a zimmer.\footnote{You can not find this term in the literature,
because it is introduced in this tractatus. 'Zimmer' is the german word for room or apartment,
which F.H"olderlin uses in the title of one of his short Scardanelli-poems. You find this poem in 
the german version of this abstract.} 
Any zimmer is the topological hull of a trajectory or a fixed point $\tau\in\gamma$.
Only closed trajectories or fixed points $\tau\in\gamma$ are identical to their 
topological hulls $\mathbf{cl}(\tau)$. We shall call any zimmer trivial, if and only if it is
a trajectory or a fixed point
$\mathbf{cl}(\tau)=\tau$. Otherwise we shall call any zimmer non-trivial.\newline
The set of the non-trivial zimmers of the phase-space $\zeta\subset\mathbb{R}^{n}$ 
is the subset $\Gamma$ of those 
topological hulls, which are not one- or zero-dimensional.
Anyway each zimmer is an attractor. 
\newline
Thus the theorem of the existence of zimmer 2.1.2. already answers the 
question, which
are the sets $\chi\subset\zeta$ of the phase-space, which are quasiergodic.
However we  
fish for another determination of the quasiergodic sets of a given dynamical system
fulfilling the three itemized criteria.
We search for
a determination of its quasiergodic sets formed by its invariants.
We find, that there is the linear homogenous partial differential equation of first order     
(\ref{miymmo}), 
which allows us the construction of the quasiergodic sets, because  
the quasiergodic sets are identical to minimal invariant manifolds. This is the other part
of the statement of the elementary quasiergodic theorem:
The quasiergodic sets are identical to the minimal invariant manifolds. 
The conception of minimal invariant manifolds is presented in the first chapter of this tractatus.
\newline
We told, that zimmers are attractors. Furthermore we will show in this tractatus, that any 
non-trivial zimmer is a sensitve attractor.
This is the statement of the theorem of sensitiveness of non-trivial zimmers
2.2.2.

\newpage
{\Large {\bf Unser Gegenstand}}
\begin{flushright}
{\em An Zimmern\\
\quad\\
Die Linien des Lebens sind verschieden,\\
Wie Wege sind, und wie der Berge Grenzen.\\
Was hier wir sind, kann dort ein Gott erg"anzen\\
Mit Harmonien und ewigem Lohn und Frieden.}\\
\quad \\
Friedrich H"olderlin\\
\end{flushright}
Wir beweisen hier einen elementaren Satz, den Satz 1.1, den wir den elementaren 
Quasiergodensatz\index{elementarer Quasiergodensatz} nennen.
Dessen Aussage k"onnen wir insofern als die L"osung 
eines altbekannten Problemes ansehen:\newline
Wenn ein dynamisches System vorliegt,
dessen Zustandsentwicklungen
st"uckweise stetig differenzierbare Trajektorien 
eines endlichdimensionalen reellen Zustandsraumes
darstellen,
deren Einheitstangenten in fast jedem Zustand existieren und die dabei ein fast "uberall
stetiges Einheitstangentenfeld festlegen,
dann
beantwortet dieser Satz n"amlich die folgende Frage:
\newline
Welche Teilmengen $\chi$ dieses beschr"ankten und
endlichdimensional-reellen Zustandsraumes $\xi\subset \mathbb{R}^{n}$
sind so beschaffen,
dass innerhalb dieser Teilmengen $\chi\subset\xi$ Quasiergodik gegeben ist?\newline
Wir betrachten ein in der Form kontinuierliches dynamisches System,
dass
alle jeweiligen Zustandsentwicklungen, 
die durch lauter jeweilige st"uckweise glatte und zusammenh"angende Trajektorien $\tau_{1}$ 
wiedergegeben sind,
fast "uberall ein Feld der Tangenten an dieselben festlegen,
das auf dem Zustandsraum fast "uberall stetig ist.    
Jede Trajektorie $\tau_{1}$, die $\chi$ schneidet,
verl"auft dann so,
dass jede andere Trajektorie $\tau_{2}$ der Trajektorie $\tau_{1}$ genau dann beliebig nahe kommt,
wenn diese andere Trajektorie $\tau_{2}$ ebenfalls $\chi$ schneidet.
Genau dann liegt $\tau_{2}$
in der abgeschlossenen H"ulle $ \mathbf{cl}(\tau_{1})$ von $\tau_{1}$, die $\chi$ ist.
Listen wir die Voraussetzungen auf! 
\begin{enumerate}
\item Der Zustandsraum $\xi$  ist beschr"ankt. 
\item Jede Trajektorie $\tau$, die eine Zustandsentwicklung 
repr"asentiert, ist
st"uckweise glatt. 
\item Die Trajektorien verlaufen fast "uberall im Zustandsraum lokal parallel.
Genau dann, wenn die Zustandsentwicklungen so beschaffen sind,
dass die Einheitstangenten an die mit den Zustandsentwicklungen identifizierbaren Trajektorien 
im jeweiligen Zustandsraum $\xi$
fast "uberall gegen einander konvergieren, sagen wir, dass diese
Trajektorien fast "uberall lokal parallel verlaufen.\index{lokale Parallelit\"at}
\end{enumerate}
Die Modifizierung "uberall gegebener lokaler Parallelit"at zur
fast "uberall vorliegenden lokalen Parallelit"at bezieht sich hierbei auf 
das Lebesgue-Mass des jeweiligen Zustandsraumes. Die Menge der Zust"ande,
in denen die Trajektorien kein stetiges Einheitstangentenfeld festlegen, ist
eine Lebesgue-Nullmenge.
Wenn 1., 2. und 3. gilt,
gibt es immer eine Partition $\Gamma$ des Zustandsraumes
$\xi$ von der Art,
dass f"ur alle jeweils Zustandsentwicklungen darstellenden Trajektorien $\tau_{1},\tau_{2}\subset\xi$ und f"ur alle Teilmengen
dieser Partition
$\chi\in \Gamma$ die "Aquivalenz
\begin{equation}\label{iplim}
\tau_{1}\cap\chi\not=\emptyset \wedge \tau_{2}\cap\chi\not=\emptyset\ \Leftrightarrow 
\mathbf{cl}(\tau_{1})=\mathbf{cl}(\tau_{2})
\end{equation}
gilt.
Daher gilt auch die "Aquivalenz
\begin{equation}
\chi=\mathbf{cl}(\tau)\ \Leftrightarrow\ \tau\cap\chi\not=\emptyset
\end{equation} 
f"ur alle Zustandsentwicklungen $\tau$ und f"ur alle Teilmengen $\chi\in \Gamma$ des Zustandsraumes
$\xi$, die 
Elemente dieser Partition $\Gamma$ sind.
\newline
Die allgemeine Beantwortung der Frage, welche Teilmengen $\chi$ dieses beschr"ankten und
endlichdimensional-reellen Zustandsraumes $\xi\subset \mathbb{R}^{n}$ unter 
der Voraussetzung, dass 1., 2. und 3. gelten,
so beschaffen sind,
dass innerhalb dieser Teilmengen $\chi\subset\xi$ Quasiergodik vorliegt, ist
uns damit gegeben: Diese Teilmengen $\chi$ sind die abgeschlossenen H"ullen $\chi=\mathbf{cl}(\tau)$ von Zustandsentwicklungen
darstellenden Trajektorien $\tau$. Dieselben nennen wir Zimmer.\index{Zimmer} Zimmer sind also kompakt.
Insofern haben wir zu zeigen, dass unter den skizzierten Voraussetzungen der 
Beschr"anktheit des Zustandsraumes,
der st"uckweisen Gl"atte der Trajektorien und
der fast "uberall gegebenen lokalen Parallelit"at der Trajektorien die "Aquivalenz
(\ref{iplim}) gilt.\newline 
Es bleibt dann aber noch die Frage,
welche jeweiligen Teilmengen derartige Zimmer 
im Phasenraum eines jeweils konkret vorgelegten dynamischen Systemes 
sind.
Die jeweilige Beantwortung dieser konkreten Frage 
scheint darin zu bestehen,
immer eine jeweilige
lineare, homogene, parzielle Differenzialgleichung 1.Ordnung formulieren zu k"onnen,
n"amlich die parzielle Differenzialgleichung (\ref{miymmo}), 
\index{Invariantengleichung einer Flussfunktion} deren L"osungen die Mannigfaltigkeiten aus der Partition $\Gamma$ festlegen, die die mit den Zimmern 
koinzidierenden minimalen invarianten Mannigfaltigkeiten sind.
Gerade die Einsichten dieses Traktates lassen uns ahnen, dass
mittels dieser parziellen Differenzialgleichung (\ref{miymmo}) die Zimmer im Einzelfall zu bestimmen, 
allerdings mit Schwierigkeiten behaftet ist. Diese Schwierigkeiten fassen wir als eine eigene Problematik und als ein
eigenes Thema auf,
das diese Abhandlung nur streift.\newline 
Unter den skizzierten Voraussetzungen 1.--3.
sind jene Zimmer Mannigfaltigkeiten und daher genau dann eindimensionale Mannigfaltigkeiten,
wenn sie geschlossene Trajektorien sind. Andernfalls zeigen sie sensitives 
Verhalten, 
wie dies aus
dem 
Satz 2.2.2 hervorgeht.
Exakt alle diese sensitven, nicht eindimensionalen Zimmer nennen wir nicht-triviale Zimmer. Die Menge
$\mathbf{cl}(\tau)\setminus\tau$ ist f"ur alle in diesen mehrdimensionalen Zimmern $\mathbf{cl}(\tau)$
liegenden Trajektorien $\tau$
nicht leer.
\newline
Zimmer sind Attraktoren\footnote{Sowohl die "ubliche, hier geltende Definition von 
Attraktoren als auch die Definition von sensitiven Attraktoren findet sich hier im
Unterabschnitt 2.2.}
und
nicht-triviale Zimmer sind gem"ass dem Satz 2.2.2 sensitive Attraktoren.\index{Attraktor}\index{sensitiver Attraktor}
Diese Attraktoren, die Zimmer, sind insofern elementar,\index{elementarer Attraktor}
als es im Phasenraum des jeweiligen Systemes keine kleineren Attraktoren gibt.
Die Nondekomponierbarkeit schon von Attraktoren bedingt diese Elementarit"at, 
welche durch deren Sensitivit"at im Fall mehrdimensionaler Attraktoren dramatisiert ist.\index{Nondekomponierbarkeit sensitiver Attraktoren}
Denn
entweder nicht-triviale Zimmer oder aber geschlossene Trajektorien, die jeweils eine periodische Zustandsentwicklung darstellen,
sind insofern "uberall im Zustandsraum:
\newline
\newline{\em
Der Zustandsraum eines Systemes mit den Eigenschaften $1.-3.$ ist nach $(1)$ und $(2)$
partioniert in elementare Attraktoren, die entweder geschlossene Trajektorien oder aber
im generischen Fall sensitive Attraktoren sind: Im Sinn dieser Aussage sind sensitive Attraktoren
der universelle Normalfall und Periodizit"at ist die Ausnahme von diesem Normalfall.}
\newline\newline
Wenn ein System mit den Eigenschaften 1.-3. gegeben ist, k"onnen
wir uns demnach in die ber"uhmteste Dramatik werfen: Ob dann ein Zustand der eines 
\begin{center}
{\bf sensitiven Attraktors oder eines Zyklus oder ein Fixpunkt ist, das ist nun die Frage\ ,}
\end{center}
die sich nach dem elementaren Quasiergodensatz f"ur 
jeden Zustand eines
jeweils vorgelegten 
dynamischen Systemes mit den Eigenschaften 1.-3. stellt.
\newpage
\tableofcontents
\chapter{Was ist das Quasiergodenproblem?}
\begin{flushright}
Allen Geduldigen: Bestandsaufnahme.\\
\end{flushright}
Sind wir nun reich?
Nein, denn es ist kein Preisgeld auf die L"osung des 
Quasiergodenproblemes gesetzt. 
Und deshalb ist auch nicht mit wenigstens derjeniger 
Genauigkeit bestimmt, die rechtliche Anspr"uche geltend machbar liesse, was als eine L"osung des 
Quasiergodenproblemes anerkannt werden m"usste; wenngleich auch eine juristisch tragf"ahige
Formulierung der Bedingungen daf"ur, das Preisgeld zu beanspruchen, immer eine solche
w"are, die naturgem"ass noch nicht auf der Sachverhaltskenntis
beruhen k"onnte, die erst nach der Probleml"osung vorliegen kann.
\newline
Was ist das Quasiergodenproblem
"uberhaupt?
Wir k"onnen und 
m"ussen wohl darauf antworten, dass 
die Negation oder die Best"atigung der 
sogenannten 
Quasiergodenhypothese als die
L"osung des
n"amlichen Problemes gelten darf.
Zu Beginn des Jahrzehntes des ersten Weltkrieges formulierten
P. und T. Ehrenfest\index{Ehrenfest, P. und T.} schliesslich
die 
Quasiergodenhypothese
als die minimalinvasive Abwandlung 
der unhaltbaren Ergodenhypothese\index{Ergodenhypothese} Ludwig Boltzmanns\index{Boltzmann, Ludwig} 
aus dem Jahre 1885. \newline
Ludwig Boltzmann hatte seinerzeit jene 
Behauptung aufgestellt,
dass jede Trajektorie der Energiehyperfl"ache eines 
Vielteichensystemes dieselbe durchlaufe, die
als Ergodenhypothese bekannt ist.
Diese f"ur den Physiker augenscheinlich sehr greifbare Fassung 
der Ergodenhypothese ist
dadurch bestimmt, 
dass ihr Gelten eine hinreichende 
Begr"undung der 
Gleichheit der sogenannten
Zeitmittel und Scharmittel
erm"oglicht.
Boltzmann deckt also mit der Ergodenhypothese seinen 
Bedarf nach
einer Aussage, die die Identit"at der
Zeitmittel mit den Scharmitteln\index{Zeitmittel-Scharmittel-Identit\"at}
tr"agt. \newline
Der Siegeszug der Mengenlehre\index{Mengenlehre} beginnt erst
nach Boltzmanns Fassung der Ergodenhypothese
und dieser Siegeszug
etabliert ein moderneres Bild,
das wir kulturhistorisch mit der mathematischen
Modernit"at schlechthin gleichsetzen k"onnen. 
In diesem Bild ist die Energiehyperfl"ache
formal identifiziert. Sie ist eine Mannigfaltigkeit eines 
geradzahligdimensionalen reellen Raumes $\mathbb{R}^{2n}$ f"ur $n\in \mathbb{N}$, deren
Kontinuit"at ein "Aquivalent der durch die Hamiltonfunktion eines jeweiligen
Vielteichensystemes vorgegebenen Kontinuit"at ist. Und daher ist die Energiehyperfl"ache
eine ungeradzahligdimensionale, n"amlich eine $2n-1$-dimensionale 
Mannigfaltigkeit.
\newline
Dabei gilt die eing"angige Aussage,
dass die Dimension
einer
Mannigfaltigkeit 
eine topologische Invariante ist,
dass es also keine Hom"oomorphismen zwischen
Mannigfaltigkeiten verschiedener Dimension gibt. Jene eing"angige Aussage wurde
aber erst sp"ater, im Zuge der Entwicklung 
der mengentheoretischen Topologie, bewiesen.
Es gibt
also auch keinen Hom"oomorphismus zwischen einer mehr als
eindimensionalen 
jeweiligen Energiehyperfl"ache
und einer eindimensionalen Trajektorie. 
\newline
Was allerdings der Fall sein kann, ist,
dass die Trajektorien allen Punkten der 
Energiehyperfl"ache beliebig nahe kommen.
Und gerade dies ist die Idee von P. und T. Ehrenfest,
dass die Trajektorien, die die Systementwicklung 
eines Vielteichensystemes darstellen,
allen Punkten derjenigen 
Hyperfl"ache, 
die durch alle Invariante des Vielteichensystemes
bestimmt ist,
beliebig nahe kommt. Die Behauptung dieses Gedankens 
als die
Quasiergodenhypothese zu bezeichnen, die als die minimalinvasive Abschw"achung 
der unhaltbaren Ergodenhypothese Ludwig Boltzmanns
gelten darf, lag nahe und diese Benennung verfestigte sich bis heute.
Wobei gerade dies die Erkenntnis von P. und T. Ehrenfests ist,
dass dann, wenn nur diese Abschw"achung der Ergodenhypothese, die Quasiergodenhypothese, 
gilt, dass dann derjenige Zweck auch erf"ullt werden kann,
den die Ergodenhypothese erf"ullen soll:
Wenn die Quasiergodenhypothese gilt, folgt die 
Identit"at der
Zeitmittel und Scharmittel auch.\newline
Es sind die Belange des Physikers,
die hier einen Bedarf 
nach der Beantwortung einer Frage
anmelden, 
die offenbar ihrem Wesen nach eine 
mathematische Frage ist. Diese Herkunft
aus der Relevanz f"ur die Physik
pr"agt dabei zun"achst das Erscheinungsbild
einer im Kern mathematischen Frage.
\section{Grundbegriffe f"ur die formale Fassung des Quasiergodenproblemes}
\subsection{Flussfunktionen und trajektorielle Partitionen}
Wir ahnen, dass jene besagte Herkunft der Quasiergodenfrage
aus dem Bedarf der Physik
einerseits
dieser wesentlich mathematischen Frage eine
Auspr"agung gibt. Und wir ahnen, dass andererseits die entsprechende,
abstrahiert emanzipierte, rein mathematische Frage,
die sich aus der Quasiergodenfrage in ihrer 
bedarfsgepr"agten Form ergibt, 
eine
allgemeinere ist, als deren bedarfsgepr"agte Form. 
Wir bleiben in dieser Abhandlung
den Wurzeln der Quasiergodenfrage in der Physik
nahe.\newline
Eine Teilmenge 
$\zeta\subset \mathbb{R}^{n}$
eines 
nicht notwendig geradzahligdimensionalen reellen Raumes $\mathbb{R}^{n}$ f"ur $n\in \mathbb{N}$
soll unseren
Zustandsraum darstellen; und dass dabei eine Abbildung
\begin{equation}\label{kurza}
\begin{array}{c}
\Psi:\zeta\times \mathbb{R}\to\zeta\ ,\\
(z,t)\mapsto\Psi(z,t)
\end{array}
\end{equation}
die Eigenschaft hat,
dass die Mengen
\begin{equation}\label{kurzb}
\Psi(z,\mathbb{R})=\{\Psi(z,t):t\in\mathbb{R}\}
\end{equation}
den Zustandsraum $\zeta$ partionieren, dass also die Implikation 
\begin{equation}\label{kurzc}
\begin{array}{c}
\Psi(z_{1},\mathbb{R}),\Psi(z_{2},\mathbb{R})\in
\{\Psi(z,\mathbb{R}):z\in\zeta\}:=[\Psi]\Rightarrow\\
\Psi(z_{1},\mathbb{R})\cap\Psi(z_{2},\mathbb{R})\in\{\emptyset,\Psi(z_{1},\mathbb{R})\}
\end{array}
\end{equation}
gilt, dies modelliert den
Determinismus.\index{Determinismus} Jede solche Abbildung $\Psi$, die gem"ass (\ref{kurza})-(\ref{kurzc})
beschaffen ist, nennen wir eine Flussfunktion im weiteren Sinn 
binnen $\zeta\subset \mathbb{R}^{n}$.\index{Flussfunktion im weiteren Sinn binnen einer Menge}
Sie sei ferner insofern normiert,
als f"ur sie
\begin{equation}\label{kure}
\Psi(\mbox{id},0)=\mbox{id}
\end{equation}
gelte. Ferner habe jede Flussfunktion $\Psi$ ihre Viabilit"at,\index{Viabilit\"at einer Flussfunktion} die Eigenschaft, dass f"ur alle $z\in\zeta$
die Alternative
\begin{equation}\label{kure}
\begin{array}{c}
\Psi(z,\mathbb{R})=\{z\}\ \lor\\
\exists\ \varepsilon(z)\in\mathbb{R}^{+}\ \forall\ \vartheta\in\ ]0, \varepsilon(z)[\ \Psi(z,\vartheta)\not=z
\end{array}
\end{equation}
wahr ist.
Dabei modelliert 
die Kontinuit"at jeder solchen Abbildung $\Psi$, jeder 
Flussfunktion 
binnen $\zeta$ also,
die Kontinuit"at des Determinismus.
Wenn die Flussfunktion $\Psi$ binnen $\zeta$ eine 
als Zeitableitung interpretierbare, stetige 
parzielle Ableitung $\partial_{2}\Psi$ hat,
so modelliert sie geeignetenfalles
den
Determinismus 
eines Hamiltonischen Vielteichensystemes. 
Denn 
dessen Hamiltonfunktion ist in dem Fall, dass $n=2m$ geradzahlig ist, zu einer 
auf dem Zustandsraum $\zeta\subset \mathbb{R}^{n}$
einmal differenzierbaren   
$\mathcal{C}^{1}$-Funktion $H$ "aquivalent, 
deren Gradient $\nabla H$ die 
Ableitung 
\begin{equation}\label{kurzd}
\partial_{2}\Psi(\mbox{id},0)=(\sigma_{3}\olessthan \mathbf{1}_{m})\nabla H
\end{equation}
ist, wobei
\begin{displaymath}
\sigma_{3} = \left( \begin{array}{cc}
0 & -1\\
1 & 0 \end{array} \right)
\end{displaymath}
die symplektische Matrix oder die 3.te der Pauli-Matrizen $\sigma_{3}$ und wobei $\mathbf{1}_{m}$ die $m\times m$-Einheitsmatrix ist und
$\olessthan$ das Kroneckerprodukt\footnote{Das Zeichen \glqq$\olessthan$\grqq wird kaum verwendet, um das 
Kroneckerprodukt zu bezeichnen, das zwar nicht selten durch das Zeichen \glqq$\otimes$\grqq
notiert wird, das allerdings gel"aufigerweise
das 
direkte Produkt bezeichnet. Das direkte Produkt verwenden wir in diesem Traktat zwar gar nicht.
Nichtsdestotrotz k"onnte die Bezeichnung des Kroneckerproduktes mit \glqq$\otimes$\grqq den
Leser in die Irre f"uhren. Da es sich bei dem Kroneckerprodukt $A \olessthan B$,
das im Symbolverzeichnis\index{Notationskonvention} definiert ist, 
um eine Art des Pumpens der Matrix $B$ in die Matrix $A$ handelt, 
halten wir diese Bezeichnungswahl f"ur ideogrammatisch vertretbar.} bezeichne. Dann ist die Gleichung (\ref{kurzd}) 
das "Aquivalent 
der Hamiltongleichungen, wobei die ersten $m$ Komponenten der Flussfunktion die nicht 
notwendigerweise kartesischen Ortskoordinaten und die zweiten $m$ Komponenten der Flussfunktion die
ihnen jeweils entsprechenden kanonischen Impulskomponenten sind.\newline
Exakt jede eindimensionale Mannigfaltigkeit eines $\mathbb{R}^{n}$
nat"urlicher Dimension $n\in\mathbb{N}$ nennen wir hier eine Trajektorie.\index{Trajektorie}
Wir nehmen in diese Abhandlung
einfachheitshalber diese Redeweise, obwohl diese Bezeichnungsweise
damit die Trajektorie aus ihrem wesentlichen partitiv-kollektiven Kontext herausreisst, als eine Menge des Mengensystemes
$\{\Psi(z,\mathbb{R}):z\in\zeta\}$ den Determinismus zu modellieren. Da wir hier aber
ohnehin den konstanten Kontext eines Trajektorienkollektives betrachten, das 
Heine-Descartesch
ist, ist diese einfache Sprechweise
hier m"oglich. 
Wir 
formulieren den Begriff der Quasizyklizit"at und den der Zyklizit"at einer Trajektorie,
einer eindimensionalen Mannigfaltigkeit des $\mathbb{R}^{n}$
nat"urlicher Dimension $n\in\mathbb{N}$: Wir nennen jede Trajektorie $\tau\subset\mathbb{R}^{n}$
genau dann zyklisch oder geschlossen, wenn sie ihre bez"uglich der nat"urlichen Topologie $\mathbf{T}(n)$ abgeschlossene H"ulle
$$\mathbf{cl}(\tau)=\tau$$
ist. Der Begriff der Quasizyklizit"at einer Trajektorie
ergibt sich,
wenn wir die Redeweise erl"autern,
dass eine Trajektorie $\tau\subset\mathbb{R}^{n}$ einem Punkt $y\in\mathbb{R}^{n}$ beliebig
nahe komme:
Genau dann, wenn es eine positive reelle Zahl $\varepsilon^{+}$ gibt, f"ur die f"ur alle positiv-reellen Zahlen 
$\varepsilon\in]0,\varepsilon^{+}[$
\begin{equation}
\mathbf{card}(\{x\in\mathbb{R}^{n}:||x-y||=\varepsilon\}\cap\tau)=\mathbf{card}(\mathbb{N})
\end{equation}
gilt, sagen wir, dass $\tau$ dem Punkt $y$ beliebig nahe komme.\index{Beliebiges Nahe-Kommen}
Quasizyklisch nun nennen wir eine Trajektorie genau dann,
wenn es einen Punkt gibt,
dem sie beliebig nahe kommt. Existieren Trajektorien, die nur einem in ihr liegenden
Punkt beliebig nahe kommen, die aber zyklisch sind, weil sie dabei zugleich mit ihrer abgeschlossenen H"ulle 
identisch sind? Das k"onnen wir uns zun"achst vielleicht vorstellen. Die Objektivierbarkeit dieser Imagination
k"onnen wir aber auf diese Weise negieren:\newline
Jede quasizyklische Trajektorie ist nicht zyklisch. W"are eine Trajektorie $\tau\subset\mathbb{R}^{n}$
n"amlich zugleich quasizyklisch und zyklisch, so existierte ein Punkt $y$, dem sie 
beliebig nahe kommt und dieser l"age in ihrer abgeschlossenen H"ulle $\mathbf{cl}(\tau)$, die mit der Trajektorie $\tau$
wegen deren gleichzeitiger Zyklizit"at identisch w"are, sodass $y\in\tau$ w"are.
Nehmen wir an, die Trajektorie $\tau$ sei glatt und habe in allen ihren Punkten eine Tangente, mithin auch in
$y$: Dann f"uhrte die Annahme, dass es keine 
positive reelle Zahl $\varepsilon^{+}$ g"abe, f"ur die f"ur alle 
$\varepsilon\in]0,\varepsilon^{+}[$
$$\mathbf{card}(\{x\in\mathbb{R}^{n}:||x-y||=\varepsilon\}\cap\tau)=2$$
ist, auf einen Widerspruch zum Zwischenwertsatz.\index{Zwischenwertsatz} F"ur jeden Hom"oomrphismus $f$
des $\mathbb{R}^{n}$ auf denselben gilt dabei offenbar die "Aquivalenz
$$\mathbf{card}(\{x\in\mathbb{R}^{n}:||x-y||=\varepsilon\}\cap\tau)=\mathbf{card}(\mathbb{N})\ \Leftrightarrow$$
$$\mathbf{card}(f(\{x\in\mathbb{R}^{n}:||x-y||=\varepsilon\})\cap f(\tau))=\mathbf{card}(\mathbb{N})\ .$$
Wir k"onnen die glatte Trajektorie $\tau$ also zu ihrem hom"oomorphen Bild $f(\tau)$ deformieren und stellen fest,
dass 
Zyklizit"at und Quasizyklizit"at einander ausschliessen.\index{Zyklizit\"at}\index{Quasizyklizit\"at}
\newline
Periodizit"at und Quasiperiodizit"at sind zun"achst relative Eigenschaften der Zust"ande der Wertemenge $\zeta$
einer Flussfunktion $\Psi$
im Verh"altnis zu derselben: Periodisch ist $z\in\zeta$ bez"uglich $\Psi$ genau dann, wenn
es eine nullverschiedene reelle Zahl $T\in\mathbb{R}\setminus\{0\}$ gibt, f"ur die 
$$\Psi(z,T)=z$$
ist. Quasiperiodisch nennen wir $z\in\zeta$ bez"uglich $\Psi$ genau dann, wenn es eine streng monoton wachsende Folge
$\{t_{j}\}_{j\in \mathbb{N}}\in \mathbb{R}^{\mathbb{N}}$ gibt, f"ur die
$$\lim_{j\to\infty} \Psi(z,t_{j})=z$$ ist
und f"ur die es dabei keine nullverschiedene reelle Zahl $T\in\mathbb{R}\setminus\{0\}$ gibt, f"ur die
$$\{t_{j}\}_{j\in \mathbb{N}}\in T\mathbb{N}^{\mathbb{N}}$$
gilt. Offenbar gilt f"ur eine stetige, d.h., f"ur eine sowohl in ihrer ersten als auch in ihrer zweiten Ver"anderlichen stetige Flussfunktion $\Psi$,
dass $z$ bez"uglich $\Psi$ 
genau dann periodisch ist, wenn $\Psi(z,\mathbb{R})=\mathbf{cl}(\Psi(z,\mathbb{R}))$ 
zyklisch ist. Wenn $z$ mit $\Psi(z,\mathbb{R})=\{z\}$ ein Fixpunkt ist, ist $\Psi(z,\mathbb{R})$ ja nicht eindimensional
und keine Trajektorie.\newline
Zyklizit"at und Quasizyklizit"at schliessen einander aus, wobei
jede in einem Kompaktum beinhaltete Trajektorie wegen des Satzes von Heine-Borel\index{Satz von Heine-Borel} 
zyklisch oder quasizyklisch ist. 
Im kompakten Zustandsraum ist es der Satz von Heine-Borel,
welcher die Poincar$\acute{e}$sche Wiederkehr\index{Poincar$\acute{e}$sche Wiederkehr} impliziert, die oft aus dem Liouville-Theorem\index{Liouville-Theorem} abgeleitet wird,
das in dem spezielleren Fall gilt, wenn eine Hamiltonische 
Vielteilchen-Interpretation\index{Vielteilchen-Interpretation}\index{Hamiltonische Vielteilchen-Interpretation} 
gem"ass (\ref{kurzd})
m"oglich ist. Daher, und weil ausserdem noch vor Newton der Gedanke des Determinismus im Raum nat"urlicher 
Topologisierung -- selbstverst"andlich, ohne damals von nat"urlichen
Topologien die Rede sein konnte -- nennen wir jede Menge von Trajektorien $[\Psi]$ genau dann Heine-Descartesch, wenn sie den 
Bedingungen (1)-(3) der Einleitung gen"ugt.\index{Heine-Descartessche Kollektivierung}\index{Heine-Descartessche trajektorielle Partition}
Damit kann dann aber auch f"ur jede hom"oomorphe Deformation eines Hamiltonischen Vielteilchen-Systemes
dessen Poincar$\acute{e}$sche Wiederkehr gezeigt werden, weil dieses Wiederkehrverhalten
eine topologische Invariante ist. 
Das Liouville-Theorem, welches dissipative Variabilit"at jeweiliger Lebesgue-positiver 
Zustandsraumteilmengen ausschliesst, veranlasst zu der Idee, die Dissipationsfreiheit\index{Dissipationsfreiheit} 
in Form der Masserhaltung in verallgemeinerten Zustandsr"aumen zu generalisieren und 
die Masserhaltung
zur Voraussetzung von Ergodens"atzen zu machen, wie beispielsweise dem Birkhoffschen 
Ergodentheorem.\index{Birkhoffsches Ergodentheorem}\newline
Dass die
Periodizit"at und die Quasiperiodizit"at eines Zustandes $z\in\zeta$ f"ur jede endlichdimensional-reelle Flussfunktion $\Psi$, deren 
Trajektorien $\Psi(z,\mathbb{R})$ glatt sind,
einander ebenfalls 
ausschliessende Alternativen sind, ist offensichtlich.
Daher gilt f"ur stetige Flussfunktionen $\Psi$
mit kompaktem Zustandsraum $\zeta$:
Die Quasiperiodizit"at bzw. die Periodizit"at eines Zustandes $z\in\zeta$ bez"uglich $\Psi$
und die Quasizyklizit"at bzw. die Zyklizit"at der Trajektorie $\Psi(z,\mathbb{R})$ sind zueinander "aquivalent.\newline
Wir k"onnen in der einander ausschliessenden Alternativit"at der Quasizyklizit"at oder Zyklizit"at 
der Trajektorien kompakter Zustandsr"aume eine Vorform des Sachverhaltes sehen, den wir hier zeigen
werden; des Sachverhaltes n"amlich, dass die Determinismen, die durch stetige Flussfunktionen $\Psi$
mit kompaktem Zustandsraum $\zeta$ beschrieben werden, so beschaffen sind,
dass f"ur jeden Zustand $z\in\zeta$ die Menge $\mathbf{cl}(\Psi(z,\mathbb{R}))$
entweder ein sensitiver Attraktor ist; oder aber die Menge $\mathbf{cl}(\Psi(z,\mathbb{R}))=\Psi(z,\mathbb{R})$ ist ein Zyklus
oder sie beschreibt einen Fixpunkt.
Daher behandeln wir hier vorzugsweise stetige Flussfunktionen
mit kompaktem Zustandsraum.
Und so bleibt es nicht aus, dass wir manchmal kurzerhand auch zyklische Trajektorien $\Psi(z,\mathbb{R})$ als periodische 
und quasizyklische Trajektorien $\Psi(z,\mathbb{R})$ als quasiperiodische Trajektorien
bezeichnen.\index{Quasiperiodizit\"at}\index{Periodizit\"at}
\index{Periodizit\"at bez\"uglich einer Flussfunktion}
\newline
Wenn wir f"ur mit Hilfe von Flussfunktionen modellierte Formen des
Determinismus eine geeignete Quasiergodenaussage
formulieren und wir diese Aussage schliesslich zeigen k"onnen, dann
wird auch das Quasiergodenproblem f"ur Vielteichensysteme
beantwortet.
Wir bemerken aber, dass hier eine wesentlich mathematische Frage
erst einmal formuliert werden muss. Es ist wie bei einem Portrait: 
Diese Formulierung muss erst einmal getroffen
werden.
\newline
Als eine mathematische Fassung des 
Quasiergodenproblemes k"onnen wir dabei nur eine Fragestellung
anerkennen, die so beschaffen ist, dass
ihre Beantwortung diejenigen Fragen mitbeantwortet, welche die
Ergodik 
eines Vielteichensystemes
betreffen.
Insofern sind wir zuallererst auf der Suche
nach dem rein
mathematischen Problem, dessen
L"osung den Bedarf
des Physikers im Hinblick auf die    
Ergodik 
eines Vielteichensystemes
deckt.
Wir suchen eine 
mathematisch formulierte Frage.\newline
Deren 
im Hinblick auf jene   
Ergodik bedarfsdeckende
Beantwortung muss dabei aber nicht alle 
Fragen der statistischen 
Mechanik eines Vielteichensystemes
kl"aren.
Es k"onnte ja schliesslich auch so sein, dass die 
Beantwortung jener mathematisch formulierten Frage
so ausf"allt, dass die historische und unsch"arfer gefasste
Quasiergodenhypothese als falsch anzusehen ist.
Dann erg"aben sich f"ur den 
Physiker neue Aufgaben, die den Aufbau der statistischen 
Mechanik betreffen.\newline
Wir verraten schon, was anders kaum zu erwarten ist.
Die Quasiergodenhypothese ist im weiteren Sinn nicht als falsch anzusehen. 
Indess, im chaosfreien Szenario ist sie K.O. gem"ass dem Trivialit"atssatz 2.2.3\index{Trivialit\"atssatz}: Diese Arbeit 
wird demonstrieren, dass 
Flussfunktionen $\Psi$, die
Heine-Descartessche Kollektivierungen $\{\Psi(z,\mathbb{R}):z\in\zeta\}$
festlegen,
durchaus chaotisch sein k"onnen. Sind sie aber nicht chaotisch, so bleibt insofern, als
dann f"ur jede Trajektorie $\tau\in \{\Psi(z,\mathbb{R}):z\in\zeta\}$ deren Zyklizit"at im Sinn negierter
Reichhaltigkeit
$$\mathbf{cl}(\tau)\setminus\tau=\emptyset$$
gilt,
von der Quasiergodik nichts "ubrig, aber auch gar nichts. 
\newline 
Dies dramatisiert die skizzierten historischen Einw"ande gegen Boltzmanns Hypothese, die 
darauf basieren, dass es keine Hom"oomorphismen zwischen Mannigfaltigkeiten verschiedener Dimensionen gibt, 
vor denen gerade die Ehrenfestsche minimalinvasive Abwandlung der Boltzmannschen Ergodenhypothese zur 
Quasiergodenhypothese retten soll.   
Was dazu verf"uhrt, menschlicher Merkf"ahigkeit entgegenzukommen, und 
die Parole von dem \glqq Ohne Chaos nichts los!\grqq zu pr"agen.\index{Ohne Chaos nichts los!}
Exakt jeder Hom"oomorphismus $\gamma$ eines endlichdimensionalen reellen 
Zustandsraumes $\zeta\subset\mathbb{R}^{n}$ auf denselben
ist h"ullenerhaltend in dem Sinn, dass f"ur alle Zustandsraumteilmengen $Z\subset\zeta$
\begin{equation}\label{thisssi}
\gamma^{-1}\circ\mathbf{cl}\circ\gamma(Z)=Z
\end{equation}  
ist und die mengenwertige Abbildung
$$\gamma^{-1}\circ\mathbf{cl}\circ\gamma=\mbox{id}$$
insofern
die Identit"at auf der Potenzmenge $2^{\zeta}$ des Zustandsraumes ist; die Hom"oomorphie der Bijektion $\gamma$ ist
auf des letzteren nat"urliche Topologie $\mathbf{T}(n)\cap\zeta$
bezogen. Daraus ergibt sich unmittelbar die Partitivit"at des Mengensystemes
$$\{\mathbf{cl}(\Psi(z,\mathbb{R})):z\in\zeta\}\ ,$$
das eine Partition des Zustandsraumes $\zeta$ ist, wenn die
Bijektionen $\Psi^{t}:=\Psi(\mbox{id},t)$ f"ur alle $t\in\mathbb{R}$
Hom"oomorphismen sind. Von dieser
Vorraussetzung gehen wir hier aber offenbar nicht aus.
Denn, dass Heine-Descartessche Kollektivierungen vorliegen, ist durch die Stetigkeit aller Fl"usse $\Psi^{t}$ f"ur $t\in\mathbb{R}$ nicht
impliziert, deren Definitionsmenge nicht beschr"ankt sein muss; und umgekehrt impliziert die Stetigkeit der
Fl"usse $\Psi^{t}$ f"ur alle $t\in\mathbb{R}$
auch nicht, dass $\{\Psi(z,\mathbb{R}):z\in\zeta\}$ Heine-Descartesch ist.\index{Heine-Descartessche Kollektivierung}   
Denken wir an diejenige Partition $Q_{2}$ der abgeschlossenen, im Ursprung zentrierten 
Einheitskreisscheibe $\bigcup Q_{2}$, welche das Mengensystem der
ebenfalls im Ursprung zentrierte Kreise 
des Radius $r\in[0,1]$ ist. Der Kreis mit Radius null ist
der Ursprung $\{0\}\not=0$. Diese Partition $Q_{2}$ ist auf verschiedene Weise
durch Flussfunktionen $\xi_{0}$ und $\xi_{1}$ in dem Sinn parametrisierbar, dass f"ur beide Indizes
$j\in\{0,1\}$ 
$$Q_{2}=\{\xi_{j}(z,\mathbb{R}):z\in\bigcup Q_{2}\}$$
ist. Es gibt aber offenbar solche Flussfunktionen $\xi_{0}$, f"ur welche die Bijektionen
$\xi_{0}^{t}=\xi_{0}(\mbox{id},t)$ f"ur alle $t\in\mathbb{R}$ stetig sind und solche
Flussfunktionen $\xi_{1}$, f"ur welche die Bijektionen
$\xi_{1}^{t}=\xi_{1}(\mbox{id},t)$ f"ur alle $t\in\mathbb{R}$ stetig sind. Wir erkennen hierbei, dass
eine Flussfunktion $\Psi$, f"ur welche die Implikation
\begin{equation}
t\in\mathbb{R}\Rightarrow (\Psi^{t})^{-1}\circ\mathbf{cl}\circ\Psi^{t}=\mbox{id}
\end{equation} 
gilt, keinerlei Sensitivit"at\index{Sensitivit\"at} auftreten kann. Sensitivit"at
kann nur vorliegen, wenn die Existenzaussage
\begin{equation} 
\exists\ t_{\star}\in\mathbb{R}:\  (\Psi^{t_{\star}})^{-1}\mathbf{cl}\Psi^{t_{\star}}\not=\mbox{id}
\end{equation}
zutrifft, die
eine Differenz behauptet, welche den Bruch einer dynamischen Symmetrie\index{Symmetriebruch}  
formuliert. Gehen wir von dem gem"ass (\ref{kurza})-(\ref{kurzc})
festgelegten Modell des Determinismus aus und nennen dabei
die eindimensionalen glatten Mannigfaltigkeiten 
$\Psi(z,\mathbb{R})$
der Menge
\begin{equation}\label{kurze}
[\Psi]:=\{\Psi(z,\mathbb{R}):z\in\zeta\}
\end{equation}
Trajektorien
und die Partition $[\Psi]$ eine trajektorielle Partition oder Kollektivierung der jeweiligen Flussfunktion!\index{trajektorielle Partition}
Wir nennen hier ferner exakt jede endlichdimensional-reelle Flussfunktion $\Psi$, deren
Kollektivierung $[\Psi]$ so beschaffen ist, dass 
es eine Flussfunktion $\tilde{\Psi}$ in der nicht leeren Faser
\begin{equation}\label{knurrze}
[\mbox{id}]^{-1}(\{[\Psi]\})=\Bigl\{\hat{\Psi}\in(\bigcup [\Psi])^{(\bigcup [\Psi])\times\mathbb{R}}:[\hat{\Psi}]=[\Psi] \Bigr\}\not=\emptyset
\end{equation} 
gibt, eine insensitive Kollektivierung\index{insensitive Kollektivierung} und
jede Flussfunktion genau dann insensitiv eichbar,\index{insensitiv eichbare Flussfunktion} wenn deren 
Kollektivierung insensitiv ist; und konsequenterweise bezeichnen wir exakt jedes Element einer Faser der Form (\ref{knurrze}) 
als eine 
insensitiv geeichte Flussfunktion.\index{insensitiv geeichte Flussfunktion} 
\newline 
Vielleicht gehen wir aber innerhalb dieses Rahmens zu weit,
wenn wir schon hier mehr verraten, n"amlich, dass sich durch den elementaren Quasiergodensatz
das Gesamtbild der statistischen 
Mechanik ver"andert, obwohl wir uns aus Gr"unden thematischer Bindung zun"achst strikt untersagen m"ussen, auf diese
Thematik einzugehen.
Die G"ultigkeit 
der Quasiergodenhypothese
zeigt sich n"amlich zuletzt in einer differenzierten 
Form, die 
neue Aufgaben
an die Theorie der Vielteichensysteme stellt:
Im generischen Fall n"amlich negiert der
elementare Quasiergodensatz selbst bei 
auftretendem Chaos
das quasiergodische Entwicklngsverhalten innerhalb der
gesamten Energiehyperfl"ache!\newline 
Das im Allgemeinen
dem elementaren Quasiergodensatz nicht gem"asse
quasiergodische Entwicklngsverhalten, dass ein
klassisches Vielteichensystem innerhalb der
gesamten Energiehyperfl"ache quasiergodisch sein soll,
kann im klassischen Modell nur durch systemfremde, durch externe statistische 
Fluktuationen ohne resultierende Energieeinspeisung in das jeweilige System begr"undet werden, {\em durch 
apriorische Annahme systemischer Unabgeschlossenheit;} oder aber in einem 
quantentheoretischen oder in einem semiklassischen
Modell dadurch, der Unsch"arferelation gem"asse Fluktuationen geltend zu machen: Wenn die genauere Analyse
erg"abe, dass klassische systemfremde, externe, statistische 
Fluktuationen anzusetzen, keine 
Erkl"arung f"ur
die im klassischen Rahmen empirisch gl"anzend best"atigte  
Boltzmannsche teilchenstatistische 
Ph"anomenlogie liefert, die ja auf der 
Annahme 
der Quasiergodik innerhalb der gesamten Energiehyperfl"ache
basiert, dann ist die Boltzmannsche Statistik wesentlich
quantentheoretisch. 
In diesem Fall w"are die Boltzmannsche Statistik wesentlich
quantentheoretisch, obwohl in ihr bekanntlich das Wirkungsquantum nicht explizit vorkommt.
\footnote{Ludwig Boltzmann entdeckt die Quantentheorie also gegebenenfalles 
-- je nachdem, wie die besagte Analyse ausf"allt --
zweimal beinahe explizit und 
zweimal implizit:\newline
1. Mit der Bestimmung der Konstante des Stefan-Boltzmannschen Strahlungsgesetzes, die
mittels des Planckschen Strahlungsgesetzes auf das Wirkungsquantum zur"uckf"uhrbar ist, entdeckt 
Boltzmann implizit das Wirkungsquantum als Experimentator.\newline 2. Dabei,
die Quasiergodik innerhalb der gesamten Energiehyperfl"ache anzunehmen,
die
f"ur
seine empirisch best"atigte  
teilchenstatistische 
Ph"anomenlogie notwendig ist, st"osst Boltzmann auf die Notwendigkeit von Vakuumfluktationen.\newline
Die besagte genauere Analyse ist aber eine Untersuchung jenseits unseres Rahmens.}
\newline
In unserem Traktat spielen bestimmte Partitionen eine grosse Rolle, eine Hauptrolle, n"amlich als eben diese trajektorielle Partitionen.
Wir benutzen als Notationskonvention\index{Notationskonvention},\footnote{Die wenigen speziellen, nicht allgemein verfestigten oder die
weniger gel"aufigen Schreibweisen, die wir in dieser Abhandlung verwenden, 
sind an deren Schluss in einem Symbolverzeichnis aufgelistet.  
Ferner sind 
diese weniger gel"aufigen Schreibweisen
im fortlaufenden
Text eingef"uhrt und die Stellen, an denen dies geschieht, sind im Index unter dem Stichwort \glqq Notationskonvention\grqq  vermerkt. } 
dass 
\begin{equation}\label{ausfp}
\begin{array}{c}
\mathbf{part} (A):=\Bigl\{Q\subset 2^{A}:\bigcup Q=A\ \land\\
q_{1},q_{2}\in Q\Rightarrow q_{1}\cap q_{2}\in\{\emptyset,q_{1}\}\Bigr\}
\end{array}
\end{equation}
die Menge aller Partitionen einer Menge $A$ notiere. Da eine Partition einer Menge ein Mengensystem, d.h., eine Menge von Mengen ist, ist
$\mathbf{part} (A)$ ein Mengensystem, dessen Elemente allesamt Mengensysteme sind.\index{Partition}
Es gilt also
$$\alpha\in\mathbf{part} (A)\ \Rightarrow\ \alpha\setminus\{\emptyset\}\in\mathbf{part} (A)\ \land\ 
\alpha\cup\{\emptyset\}\in\mathbf{part} (A)\ ,$$
sofern $A\not=\emptyset$ ist. Denn es ist 
$$\mathbf{part} (\emptyset)=\Bigl\{\emptyset,\{\emptyset\},\{\emptyset,\{\emptyset\}\}\Bigr\}\ .$$ 
Dieses Traktat ist an der Grenze, nicht der des Autors. Es behandelt den Determinismus n"amlich mit den Methoden der Analysis, die in die der mengentheoretischen
Topologie "ubergehen.
Bereits in diesem Traktat spielen daher zwar Topologien eine durchaus grosse Rolle und zwar die 
nat"urlichen, diejenigen, deren Basis euklidische Kugeln endlichdimensionaler und reeller R"aume sind. Die Mengensysteme, die Topologien sind,
teilen sich hier noch nicht die beiden Hauptrollen mit den Partitionen,
mit denen erstere sich aber 
die beiden Hauptrollen
in der abstrakten Lehre des Determinismus 
teilen, einer Theorie, die der Dialog 
der beiden Typen von Mengensystemen ist. 
Daher erw"ahnen dies wir schon hier:
Zu jeder Partition $\alpha\in\mathbf{part} (A)$ gibt es das 
Mengensystem $\alpha^{\cup}$, wenn wir f"ur jedes Mengensystem $\beta$
\begin{equation}\label{ausfq}
\beta^{\cup}:=\Bigl\{\bigcup a:a\subset\beta\Bigr\}
\end{equation}
setzten,
das, gleich, ob $\emptyset\in\beta$ oder ob $\emptyset\not\in\beta$,
allemal die leere Menge als Element hat, da $\emptyset\subset\beta$ und $\bigcup\emptyset=\emptyset$ ist.
Im Allgemeinen ist $\beta^{\cup}$ zwar wie eine Topologie abgeschlossen gegen jedwede Vereinigung $\bigcup b$ einer
Teilmenge $b\subset\beta$ und es ist ausserdem $\emptyset\in \beta$.
Im Allgemeinen ist dennoch
$\beta^{\cup}$ 
offenbar keine Topologie, wohingegen
$\alpha^{\cup}$ f"ur jede Partition $\alpha\in \mathbf{part} (A)$ eine Topologie ist, falls
$A\not=\emptyset$ ist.
Denn das Mengensystem $\alpha^{\cup}$ ist abgeschlossen gegen Schnitte endlich vieler seiner Elemente,
wobei $\bigcup\alpha^{\cup}\not=\emptyset$ ist. Deswegen nennen wir exakt 
$\alpha^{\cup}$ f"ur jede Partition $\alpha\in\mathbf{part} (A)$ f"ur $A\not=\emptyset$ die
von der Partition $\alpha$ generierte Topologie und bezeichnen jede Topologie $\mathbf{T}$
genau dann als partitionsgeneriert, wenn es eine Partition $\vartheta\in\mathbf{part}(\bigcup\mathbf{T})$ gibt, f"ur die
$$\vartheta^{\cup}=\mathbf{T}$$
ist.
Das Verh"altnis einer trajektoriellen Partition zu einer Flussfunktion $\Psi$, die diese Partition gem"ass (\ref{kurze}) festlegt, dr"ucken wir dadurch aus, dass wir
sagen, dass $[\Psi]$ die Kollektivierung oder das Trajektorienkollektiv der Flussfunktion $\Psi$ ist.  
Und wir nennen dabei exakt jede stetige bzw. st"uckweise differenzierbare bzw. differenzierbare Funktion $\Psi$, f"ur die (\ref{kurza})-(\ref{kurzc}) gilt, eine 
$\mathcal{C}^{0}$- bzw. $\mathcal{C}^{1/2}$- bzw.
$\mathcal{C}^{1}$-Flussfunktion.\index{Flussfunktion}\index{differenzierbare Flussfunktion}\index{stetige Flussfunktion}\index{st\"uckweise differenzierbare Flussfunktion} Und entprechend nennen wir genau dann das Mengensystem $[\Psi]$ eine trajektorielle $\mathcal{C}^{0}$- bzw. $\mathcal{C}^{1/2}$- bzw. $\mathcal{C}^{1}$-Partition.\index{trajektorielle $\mathcal{C}^{1}$-Partition}\index{trajektorielle $\mathcal{C}^{0}$-Partition}\index{trajektorielle $\mathcal{C}^{1/2}$-Partition}\index{trajektorielle $\mathcal{C}^{1}$-Partition}
\subsection{Minimale invariante Mannigfaltigkeiten und die Relativierung der Hom"oomorphie}
F"ur jede Flussfunktion, gleich welche Kontinuit"at sie hat, sind Invariante
als deren Invariante konzipierbar.\index{Invariante einer Flussfunktion}
Wir verstehen unter einer Invariante\index{Invariante einer trajektoriellen Partition} 
jeder trajektoriellen Partition $[\Psi]$ gem"ass (\ref{kurze})
eine besondere Abbildung, die auf dem Zustandsraum $\zeta \subset \mathbb{R}^{ n}$ 
definiert ist und deren Werte
reelle Zahlen sind: Die Besonderheit einer Invariante einer trajektoriellen Partition ist es, dass
sie auf jeweils jeder 
Trajektorie $\tau$ der trajektoriellen Partition  $[\Psi]$
konstant ist. Die Invarianten der jeweiligen trajektoriellen Partition  $[\Psi]$ sind also 
die
Elemente der Menge von Funktionen
\begin{equation}\label{invoxox}
\mbox{Inv}^{\star} ([\Psi]):= \Bigl\{J\in\ \mathbb{R}^{ \zeta}: x\in\tau\in [\Psi]\Rightarrow J (\tau)=\{J (x)\}\Bigr\}\ .
\end{equation} 
Diese Menge umfasst die Menge $\bigcup\{\{c\}^{ \zeta}:c\in \mathbb{R}\}$, die Menge
aller Funktionen, die undifferenziert auf dem ganzen Zustandsraum $\zeta$ konstant sind und exakt welche
wir als triviale Invariante\index{triviale Invariante} bezeichnen. Mit der Menge der Funktionen
\begin{equation}
\mbox{Inv} ([\Psi]):=\mbox{Inv}^{\star} ([\Psi]) \setminus \bigcup \{\{c\}^{ \zeta}:c\in \mathbb{R}\}
\end{equation} 
sind daher alle nichttriviale Invariante notiert. F"ur alle halb-oder ganzzahligen Indizes 
$j\in \{0,1/2,1\}$ und f"ur alle trajektoriellen Partitionen $[\Psi]$
sei mit 
\begin{equation}
\mbox{Inv}^{ j} ([\Psi]):=\mbox{Inv} ([\Psi]) \cap\ \{f|\zeta : f\in \mathcal{C}^{ j} ( \mathbb{R}^{ n}, \mathbb{R})\}
\end{equation} 
die Menge aller
nichttrivialer Invarianter 
der trajektoriellen Partition $[\Psi]$ notiert, die die
Kontinuit"at\index{Notationskonvention} $j\in \{0,1/2,1\}$ haben: Denn dabei sei
$\mathcal{C}^{0} ( \mathbb{R}^{ n}, \mathbb{R})$ die Menge aller stetiger, 
auf dem $\mathbb{R}^{ n}$ definierter
Funktionen, deren Werte reelle Zahlen sind und
$\mathcal{C}^{1} ( \mathbb{R}^{ n}, \mathbb{R})$ sei die in $\mathcal{C}^{0} ( \mathbb{R}^{ n}, \mathbb{R})$   
enthaltene Menge der differenzierbaren reellen Funkrionen, die auch in der Menge
$\mathcal{C}^{1/2}( \mathbb{R}^{ n}, \mathbb{R})$ 
aller Lipschitz-stetiger reeller Funktionen liegt. 
Wie zumeist in der Literatur bezeichne dabei $g|A$ auch in dieser Abhandlung f"ur jede Funktion $g$ und jede Menge $A$  
die Restriktion der
Funktion $g$ auf die Menge $A$. $g|A$ ist demnach die leere Menge, wenn die Menge $A$ und die 
Definitionsmenge der Funktion $g$ disjunkt sind. Halbzahlige Indizes 
$q\in\{-1,-1/2,0,1/2,1,3/2\dots\}$ nennen wir
Kontinuit"atsindizes,\index{Kontinuit\"atsindex} wenn sie die Kontinuit"at von Funktionen kennzeichnen,
die wir dabei als von $\mathcal{C}^{q}$-Kontinuit"at bezeichnen.
Der Kontinuit"atsindex $-1/2$ kennzeichnet st"uckweise Stetigkeit; und selbst die 
st"uckweise Stetigkeit muss f"ur eine Funktion der 
Kontinuit"at $\mathcal{C}^{-1}$ nicht gegeben sein.
\newline
Wie ist es m"oglich, dass wir hier frei von der $\mathcal{C}^{q}$-Kontinuit"at von Abbildungen $\phi$ f"ur einen jeweiligen Kontinuit"atsindex
$q$ reden, frei, d.h., ohne dabei die Referenz auszusprechen, auf welche jeweils 
zugrundegelegte Topologisierung $\mathbf{T}(\mathbf{P}_{1}\phi)$ bzw. $\mathbf{T}(\mathbf{P}_{2}\phi)$
der jeweiligen Definitionsmenge $\mathbf{P}_{1}\phi$ bzw.
Wertemenge $\mathbf{P}_{2}\phi$ wir uns beziehen? 
Dies ist deshalb m"oglich, weil wir vereinbaren, dass dann, wenn wir diese Referenz
nicht ausdr"ucklich angeben, 
die Topologisierungen
\begin{equation}\label{ausfy}
\begin{array}{c}
\mathbf{T}(n(\mathbf{P}_{1}\phi))\cap\mathbf{P}_{1}\phi =\mathbf{T}(\mathbf{P}_{1}\phi)\ ,\\
\mathbf{T}(n(\mathbf{P}_{2}\phi))\cap\mathbf{P}_{1}\phi =\mathbf{T}(\mathbf{P}_{2}\phi)
\end{array}
\end{equation}
stillschweigend angenommen seien,\footnote{Wir sind darauf bedacht, zwischen 
den Werten $f^{-1} (a)\in X$ der Umkehrfunktion $f^{-1}$ einer  
Funktion $f\in Y^{X}$ 
und den
Urbildmengen 
$$f^{-1} (A):=\{x\in \mathbf{P}_{1}f:f(x)\in A\}\ \subset X\ ,$$
die f"ur beliebige
Funktionen $f\in Y^{X}$ und f"ur beliebige Mengen $A$ definiert seien,
zu unterscheiden. Die 
Werte $f^{-1} (a)\in X$ der Umkehrfunktion $f^{-1}$ sind
nur 
f"ur Wertemengenelemente $a\in \mathbf{P}_{2}f\subset Y$ definiert.
Wenn $b\in Y\setminus \mathbf{P}_{2}f$ 
zwar ein Element der Menge $Y$ ist, in die die Funktion $f\in Y^{X}$ abbildet,
$b$ jedoch kein Element der  
Wertemenge $\mathbf{P}_{2}f$ ist, 
so ist $f^{-1} (b)$ nicht definiert. Die Urbildmenge $f^{-1} (\{b\})=\emptyset$ ist
allerdings sehr wohl definiert.
Insbesondere die Urbildmengen einelementiger Mengen $\{a\}\subset \mathbf{P}_{2}f$ der jeweiligen Wertemenge
$$f^{-1} (\{a\})=\{f^{-1} (a)\}\not=f^{-1} (a)$$
und die Werte der Umkehrfunktion $f^{-1} (a)$
wollen wir auseinanderhalten. }
falls $\phi$ eine Abbildung ist, deren 
Definitionsmenge $\mathbf{P}_{1}\phi$ im $\mathbb{R}^{n(\mathbf{P}_{1}\phi)}$ und
deren Wertemenge $\mathbf{P}_{2}\phi$ im $\mathbb{R}^{n(\mathbf{P}_{2}\phi)}$ liegt,
wobei $n(\mathbf{P}_{1}\phi)$ und $n(\mathbf{P}_{2}\phi)$ nat"urliche Zahlen seien und wobei 
$\mathbf{T}(n(\mathbf{P}_{1}\phi))$ bzw. $\mathbf{T}(n(\mathbf{P}_{2}\phi))$
die nat"urlichen Topologien des $\mathbb{R}^{n(\mathbf{P}_{1}\phi)}$ bzw. des 
$\mathbb{R}^{n(\mathbf{P}_{2}\phi)}$ seien. F"ur jede Funktion $\phi$ wollen wir hier f"urderhin mit
$\mathbf{P}_{1}\phi$ deren Definitionsmenge 
und mit $\mathbf{P}_{2}\phi$ deren Wertemenge bezeichnen.\index{Notationskonvention} 
Denn f"ur jede nat"urliche Zahl $n\in\mathbb{N}$ und f"ur jedes $n$-Tupel $a=(a_{1},a_{2},\dots a_{n})$ 
-- gleich, welcher Art dessen jeweilige Komponenten sein m"ogen -- sei f"ur alle $j\in\{1,2,\dots n\}$
\begin{equation}\label{ausfxx}
\mathbf{P}_{j}a=a_{j}\ .
\end{equation}
Da jede Funktion $\phi$ eine Menge von Paaren ist, ist $\mathbf{P}_{1}$ in dem Ausdruck $\mathbf{P}_{1}\phi$ 
als mengenweise angewandter Operator verstanden, sodass $\mathbf{P}_{1}\phi$ die Menge aller
ersten Komponenten der Menge von Paaren $\phi$, mithin die Definitionsmenge von $\phi$ ist.
Die Angabe 
der Referenztopologisierung kann daher zwar im Allgemeinen nicht
unterbleiben, hier in diesem Traktat in der Regel allerdings schon,
weil wir in demselben in der Regel von 
Abbildungen reden, die Teilmengen endlichdimensionaler reeller R"aume
auf ebensolche abbilden.
\newline
Im Allgemeinen jedoch bezeichne
f"ur je zwei Topologien $\mathbf{T}_{1}$ und $\mathbf{T}_{2}$
und 
zwei Mengen ${\rm A}_{1}\subset \bigcup\mathbf{T}_{1}$ und ${\rm A}_{2}\subset\bigcup\mathbf{T}_{2}$
der zugeh"origen topologischen R"aume
$(\bigcup\mathbf{T}_{1},\mathbf{T}_{1})$ und $(\bigcup\mathbf{T}_{2},\mathbf{T}_{2})$ 
\begin{equation}\label{ausfx}
\begin{array}{c}
\mathcal{C}_{\mathbf{T}_{1},\mathbf{T}_{2}}({\rm A}_{1},{\rm A}_{2}):=\Bigl\{f\in{\rm A}_{1}^{{\rm A}_{2}}:\\
{\rm U}\in\mathbf{T}_{2}\cap{\rm A}_{2}\Rightarrow f^{-1}({\rm U})\in\mathbf{T}_{1}\cap{\rm A}_{1}\Bigr\}\ , 
\end{array}
\end{equation}
die Menge der 
bez"uglich der Topologien $\mathbf{T}_{1}$ und $\mathbf{T}_{2}$
stetigen Abbildungen der Menge ${\rm A}_{1}$ auf die Menge ${\rm A}_{2}$.
Wobei f"ur jede Topologie $\mathbf{T}$ und jede Menge ${\rm A}$
$$\mathbf{T}\cap {\rm A}:=\{{\rm U}\cap {\rm A}:{\rm U}\in \mathbf{T}\}$$
das Mengensystem aller Schnitte einer offenen Menge der Topologie $\mathbf{T}$
mit der Menge ${\rm A}$ bezeichnet. In der 
Allgemeinheit beliebiger Topologien brauchen wir in dieser Abhandlung keine
h"oheren Kontinuit"aten als die der Stetigkeit.
\newline
Enth"alt $\zeta$ keine offene Menge,
\footnote{K"onnen wir sagen, dass der jeweilige Zustandsraum $\zeta\subset\mathbb{R}^{n}$ f"ur $n\in\mathbb{N}$ die Kontinuit"aten Invarianter einer trajektoriellen Partition $[\Psi]$
{\em im Allgemeinen} unterscheide? Dazu m"ussten wir doch eigentlich auf der Potenzmenge des $\mathbb{R}^{n}$ 
eine Sigma-Algebra von Mengensystemen $\sigma(n)$ von Teilmengen des $\mathbb{R}^{n}$ einrichten;
$\sigma(n)$ sollte 
das Mengensystem $X_{n}$ aller Teilmengen des $\mathbb{R}^{n}$, die eine offene Menge ${\rm U}\in\mathbf{T}(n)$ enthalten,
als Element haben,
wobei $\mathbf{T}(n)$ die nat"urliche Topologie des $\mathbb{R}^{n}$ ist;
wir m"ussten diese Sigma-Algebra $\sigma(n)$ mit einem 
Wahrscheinlichkeitsmass $p_{n}$ bewerten. Wenn dann $p_{n}(X_{n})=1$ ist, k"onnen wir sagen,
dass der jeweilige Zustandsraum $\zeta$ die Kontinuit"aten Invarianter 
trajektorieller Partitionen $[\Psi]$, die von $\mathcal{C}^{\infty}$- Kontinuit"at sind,
$p_{n}$-fast immer unterscheide.}
so sind die Mengen
$\mbox{Inv} ([\Psi])$ und $\mbox{Inv}^{ j}\ ([\Psi])$ f"ur alle 
Kontinuit"atsindizes $j\in \{-1,-1/2,0,\dots\}$ gleich, andernfalls differieren die 
Mengen
$\mbox{Inv}^{ j(1)}\ ([\Psi])$ und $\mbox{Inv}^{ j(2)}\ ([\Psi])$, falls die Kontinuit"atsindizes
$j(1), j(2)\in \{-1,-1/2,0,\dots\}$ verschieden sind und -- falls die Kontinuit"at der trajektoriellen Partition $[\Psi]$
hinreichend hoch ist. Was aber gibt uns ein Kriterium daf"ur, dass die Kontinuit"at der 
jeweiligen
trajektoriellen Partition $[\Psi]$ \glqq hinreichend hoch\grqq  ist?
\newline  
Wir k"onnen 
die beschriebene Auffassung der Invarianten einer trajektoriellen Partition $[\Psi]$
topologisch phrasieren:
Das Mengensystem
$$\{\phi^{-1}(\{a\}):a\in \mathbf{P}_{2}\phi\}\ \in\ \mathbf{part}(\mathbf{P}_{1}\phi)$$
ist das Niveauliniensystem der Funktion $\phi$ und offenbar eine Partition
deren Definitionsmenge, sodass
$$\{\phi^{-1}(\{a\}):a\in \mathbf{P}_{2}\phi\}^{\cup}$$
gem"ass der Ausf"uhrungen bei
(\ref{ausfp}) und (\ref{ausfq}) eine partitionsgenerierte Topologie\index{partitionsgenerierte Topologie} ist, exakt welche wir die 
Niveaulinientopologie der Funktion $\phi$ nennen.\index{Niveaulinientopologie einer Funktion}
Betrachten wir speziell eine Funktion $f$, deren Definitionsmenge 
$\mathbf{P}_{2}f=\zeta$ der Zustandsraum $\zeta$ einer trajektoriellen Partition $[\Psi]$ ist und deren Wertemenge $\mathbf{P}_{2}f\subset\mathbb{R}$
auf dem Zahlenstrahl $\mathbb{R}$ liegt:
Genau dann, wenn die Niveaulinientopologie der Funktion $f$ 
so beschaffen ist, dass die Niveaulinientopologie 
$\{f^{-1}(\{a\}):a\in \mathbf{P}_{2}f\}^{\cup}$
der Funktion $f$
gr"ober ist, als die von der trajektoriellen Partition $[\Psi]$ generierte Topologie $[\Psi]^{\cup}$,
so ist $f$ eine Invariante der trajektoriellen Partition $[\Psi]$.
Die Niveaulinientopologie 
$\{f^{-1}(\{a\}):a\in \mathbf{P}_{2}f\}^{\cup}$
der Funktion $f$ gilt dabei genau dann als gr"ober als die von der trajektoriellen Partition $[\Psi]$ generierte Topologie $[\Psi]^{\cup}$,
wenn die Implikation
\begin{equation}
a\in \{f^{-1}(\{a\}):a\in \mathbf{P}_{2}f\}^{\cup}\Rightarrow a\in [\Psi]^{\cup}
\end{equation} 
wahr ist. Der Satz "uber impizite Funktionen sagt uns dabei, dass die 
Kontinuit"at der Funktion $f$ die ihrer Niveaulinien, d.h., der Elemente
ihres Niveauliniensystemes ist.
Daher k"onnen wir die letztgestellte Frage nach einem Kriterium daf"ur, dass die Kontinuit"at der 
jeweiligen
trajektoriellen Partition $[\Psi]$ daf"ur hinreichend hoch ist, dass f"ur $j(1)\not=j(2)$ die Verschiedenheit
\begin{equation}\label{ausfr}
\mbox{Inv}^{ j(1)}\ ([\Psi])\not=\mbox{Inv}^{ j(2)}\ ([\Psi])
\end{equation} 
gilt, so beantworten:
Falls die trajektorielle Partition $[\Psi]$ nicht nirgends mindestens von der Kontinuit"at 
$\max\{j(1),j(2)\}$ ist, d.h., falls es eine offene Menge ${\rm U}$ des Zustandsraumes $\zeta$ gibt,
f"ur die gilt, dass alle Elemente des Spurmengensystemes 
$$[\Psi]\cap{\rm U}=\{\tau\cap{\rm U}:\tau\in [\Psi]\}$$
Trajektorienabschnitte der Kontinuit"at
$\mathcal{C}^{\max\{j(1),j(2)\}}$ sind, so gilt die Verschiedenheit (\ref{ausfr}).
Wir wollen nicht dar"uber streiten, ob Phrasierungen Perspektiven "offnen k"onnen.
\newline 
Allemal ist $\mbox{Inv}^{ 1/2} ([\Psi])$ die Menge
st"uckweise differenzierbarer und dabei stetiger Invarianter der trajektoriellen Partition $[\Psi]$, die keine triviale Invariante sind 
und $\mbox{Inv}^{ 0} ([\Psi])$ ist die Menge stetiger Invarianter der trajektoriellen Partition $[\Psi]$
und $\mbox{Inv}^{ 1} ([\Psi])$ ist die Menge stetig differenzierbarer Invarianter derselben, 
die keine triviale Invariante sind.
\newline
Nun wollen wir uns einen
einfachen, jedoch sehr wichtigen Sachverhalt klar machen, der die Quasiergodik einer trajektoriellen Partition 
betrifft, dessen Invariante wir nun treffend formuliert haben:
Welche 
im Hinblick auf die Quasiergodik relevanten Gegebenheiten
finden wir in dem Fall vor, dass eine nicht leere Menge nichttrivialer stetiger Invarianter einer trajektoriellen Partition existiert? 
\newline
Es sei $J\in \mbox{Inv}^{ 0}([\Psi])\not=\emptyset$ eine stetige Invariante der trajektoriellen Partition  $[\Psi]$.
F"ur alle Werte $\omega_{1}, \omega_{2} \in J\ (\mathbb{R})$ dieser stetigen Invariante $J$ ist
die Mannigfaltigkeit
\begin{equation}
J^{-1}\ (\{\omega_{1}\})=\Bigr\{x\in \zeta:J(x)=\omega_{1}\Bigr\}
\end{equation} 
eine abgeschlossene Menge.
Offensichtlich ist die "Aquivalenz
\begin{equation} 
J^{-1}\ (\{\omega_{1}\})\ \cap\ J^{-1}\ (\{\omega_{2}\}) = \emptyset\  
\Leftrightarrow \omega_{1}\not=\omega_{2}
\end{equation} 
wahr.
"Uberdies gilt f"ur die stetige nichttriviale Invariante $J$, dass f"ur alle ihre Werte 
$$\omega_{1}\in J\ (\mathbb{R}),\ \omega_{2}\in J\ (\mathbb{R}) \setminus \{ \omega_{1}\}\not=\emptyset$$ eine positive Zahl  $\delta\in \mathbb{R}^{ +}$ existiert, 
f"ur die die Disjunktionsbehauptung
\begin{equation}
(J^{-1}\ (\{\omega_{2}\})\  +\ \mathbb{B}_{ \delta} (0))\ \cap\ (J^{-1}\ (\{\omega_{2}\})\ +\ \mathbb{B}_{ \delta} (0))\ =\ \emptyset
\end{equation}
wahr ist. 
Es ist also offensichtlich ausgeschlossen, dass eine in 
der Punktmenge
$J^{-1}\ (\{\omega_{1}\})$ des Zustandsraumes 
$\zeta$
enthaltene Trajektorie keinem Punkt
$x\in J^{-1}\ (\{\omega_{2}\})$ beliebig nahe liegt, wenn $\omega_{1}\not=\omega_{2}$ ist.
Dass ein Zustand $x\in J^{-1}\ (\{\omega_{2}\})$ einer Trajektorie $\theta$
aus dem Spurmengensystem 
$$[\Psi]\cap J^{-1}\ (\{\omega_{2}\}):=\Bigl\{\tau\cap J^{-1}\ (\{\omega_{2}\}):\tau\in[\Psi]\Bigr\}$$
$$=\Bigl\{\tau\in[\Psi]:\tau\subset J^{-1}\ (\{\omega_{2}\})\Bigr\}\not=\emptyset$$
beliebig nahe liegt,
heisst dabei gerade, dass $x$ nicht in der abgeschlossenen H"ulle 
$\mathbf{cl} (\theta)$
der Trajektorie $\theta$ liegt und dass der Zustand $x\in J^{-1}\ (\{\omega_{2}\})$ von der abgeschlossenen H"ulle $\mathbf{cl} (\theta)$ bez"uglich der 
nat"urlichen Topologie des Zustandsraumes separiert ist. Daher gilt auch
die Implikation
\begin{equation}
\begin{array}{c}
\omega_{1}\in J\ (\mathbb{R})\ \land\ \omega_{2}\in J\ (\mathbb{R}) \setminus \{ \omega_{1}\}\ \land\ 
x\in J^{-1}\ (\{\omega_{1}\})\\
\Rightarrow\\
x\not\in \mathbf{cl} (J^{-1}\ (\{\omega_{1}\}))\ .
\end{array}
\end{equation}
{\em Es ist also
von vorne herein klar, dass Quasiergodik nur innerhalb der Urbilder stetiger Invarianter einer trajektoriellen Partition  vorliegen kann,
weil solche Urbilder abgeschlossen sind und weil
diese Urbilder
im Zustandsraum separiert sind.}\newline\newline
Urbilder stetiger Invarianter einer trajektoriellen Partition spielen also eine entscheidende Rolle bei der formalen Fassung der Quasiergodenhypothese,\index{Quasiergodenhypothese}
weswegen wir uns dieselben anschauen wollen:
\newline
Ein nicht leeres Urbild $J^{-1}\ (\{\omega\})$ f"ur $J \in \mbox{Inv}^{ p} ([\Psi])$
nennen wir genau dann eine invariante Mannigfaltigkeit\index{invariante Mannigfaltigkeit} der Invariante $J \in \mbox{Inv}^{ p} ([\Psi])$,
wenn $p \in\{ 1/2,1\}$ ist.
Sei nun f"ur $j\in \mathbb{N}$ und f"ur einen halb-oder ganzzahligen Index $q\in\{ 1/2,1\}$
\begin{equation}\label{nimon}
\beta:= (\beta_{ k})_{\ 1\leq k\leq j}\ \in \ (\mbox{Inv}^{ q} ([\Psi]))^{j} 
\end{equation}
ein $j$-Tupel Invarianter der trajektoriellen Partition  $[\Psi]$ mit $\mathcal{C}^{ q}$- Kontinuit"at.\newline
Dabei seien erstens dieses Tupels Komponenten so beschaffen, dass
$\{\beta_{\ 1},...\beta_{\ j}\}$
unabh"angige Invariante sind, dass es also f"ur alle $k \in \{1,..j\}$ keine $\mathcal{C}^{ q}$-Funktion  
$g_{ k}\in \ \mathcal{C}^{ q}\ (\mathbb{R}^{ j-1},\ \zeta)$
von der Art gibt, dass 
\begin{equation}\label{nimom}
\beta_{ k} = g_{ k} \ (\beta_{\ 1},..\beta_{\ k-1},\beta_{\ k+1},..\beta_{\ j})
\end{equation} 
ist;
dennoch gebe es ein Tupel $\omega=(\omega_{ k})_{\ 1\leq k\leq j}\in \mathbb{R}^{j}$, f"ur das
der Schnitt
\begin{equation}\label{mimom}
\beta^{ -1} (\{\omega\})\ :=\ \bigcap \Bigl\{\beta_{ k}^{ -1}\ (\omega_{ k})\ :k\in\{1,..j\}\Bigr\}\ \not=\emptyset
\end{equation}
nicht leer ist. Die Bedingung (\ref{mimom}) ist eine Koh"arenzbedingung
an die unabh"angige 
Menge $\mathcal{C}^{ q}$-Invarianter $\{\beta_{\ 1},...\beta_{\ j}\}$, die
deren Unabh"angigkeit entgegnet.\newline
Es sei 
f"ur jede trajektorielle Partition  
und f"ur jeden Kontinuit"atsindex $q\in\{ 1/2,1\}$
mit 
$\mbox{Ui}^{ q}([\Psi])$ die folgende Hilfskonstruktion bezeichnet:
$\mbox{Ui}^{ q}([\Psi])$ sei die Menge
aller $j$-Tupel $\beta$ f"ur ein $j\in\mathbb{N}$, die gem"ass (\ref{nimon}) beschaffen sind und die
dabei die 
Menge aller jeweiligen
Komponenten 
$\{\beta_{\ 1},...\beta_{\ j}\}$ haben, die
eine
Menge
unabh"angiger $\mathcal{C}^{ q}$-Invarianter der 
trajektoriellen Partition $[\Psi]$ ist, wobei
es ein Tupel $\omega=(\omega_{ k})_{\ 1\leq k\leq j}\in \mathbb{R}^{j}$ gebe, 
f"ur das der Schnitt $\beta^{ -1} (\{\omega\})$ gem"ass (\ref{mimom})
nicht leer ist.  
Genau dann, wenn f"ur $\beta\in \mbox{Ui}^{ q}([\Psi])$ die Implikation
\begin{equation}\label{mioom}
\begin{array}{c}
\beta_{\ j+1}\in\mbox{Inv}^{q} ([\Psi])\Rightarrow\\
\beta\oplus\beta_{j+1}=(\beta_{\ 1},...\beta_{\ j},\beta_{\ j+1})\not\in \mbox{Ui}^{ q}([\Psi])
\end{array}
\end{equation}
wahr ist,
nennen wir das Tupel Invarianter $\beta\in \mbox{Ui}^{ q}([\Psi])$ maximal. Mit $\oplus$ notieren wir dabei hier 
die Kokatenation\index{Notationskonvention} und
$\overline{\mbox{Ui}}^{ q}([\Psi])$ bezeichne nun 
f"ur jede trajektorielle Partition $[\Psi]$ 
und f"ur jeden der beiden halb-oder ganzzahligen Kontinuit"atsindizes $q\in\{ 1/2,1\}$
die Menge aller maximalen Tupel der Menge $\mbox{Ui}^{ q}([\Psi])$.
Genau dann,
wenn $\beta\in\overline{\mbox{Ui}}^{ q}([\Psi])$ ein maximales Tupel ist und
es ein Element
$\omega=(\omega_{ k})_{\ 1\leq k\leq j}$ in der Wertemenge $\beta(\zeta)$ gibt,
f"ur das der (\ref{mimom}) gem"asse Schnitt 
$\beta^{ -1} (\{\omega\})$
nicht leer ist,
nennen wir denselben, d.h., den Schnitt
\begin{equation}\label{mimmo}
\beta^{ -1} (\{\omega\})\ :=\ \bigcap \Bigl\{\beta_{ k}^{ -1}\ (\omega_{ k})\ :k\in\{1,..j\}\Bigr\}\ ,
\end{equation}
eine
minimale invariante 
$\mathcal{C}^{ q}$- Mannigfaltigkeit der trajektoriellen Partition \index{minimale invariante Mannigfaltigkeit
einer trajektoriellen Partition } $[\Psi]$. Wir klassifizieren
minimale invariante Mannigfaltigkeiten analog 
zur Klassifikation invarianter Mannigfaltigkeiten
nach ihrer Kontinuit"at:
Ist $q=1/2$, so ist $\beta^{ -1} (\{\omega\})$ eine st"uckweise differenzierbare minimale invariante Mannigfaltigkeit und f"ur $q=1$ 
ist sie eine differenzierbare Mannigfaltigkeit. F"ur die beiden Kontinuit"atsindizes $q\in\{ 1/2,1\}$
sei die Menge all solcher 
minimaler invarianter 
$\mathcal{C}^{ q}$- Mannigfaltigkeiten\index{Notationskonvention} einer jeweiligen trajektoriellen Partition 
f"urderhin notiert mit $\mathcal{M}^{  q}\ ([\Psi])$.\newline
Wenn wir die minimalen invarianten 
$\mathcal{C}^{ 1}$- Mannigfaltigkeiten
bestimmen wollen, so brauchen wir nur alle L"osungen $\beta_{j}\in\mathcal{C}^{ 1}(\zeta,\mathbb{R})$ zu bestimmen,
die f"ur alle $x\in \zeta$ der 
Gleichung
\begin{equation}\label{miymmo}
\partial_{2}\Psi(x,0)^{\top}\nabla\beta_{j}(x)=0
\end{equation}
gen"ugen, die wir die Invariantengleichung 
der jeweiligen Flussfunktion $\Psi$ nennen.\index{Invariantengleichung einer Flussfunktion} Auf den ersten Blick meinen wir, 
dass, wenn $\partial_{2}\Psi(x,0)\not=0$ ist und $n$ die Dimension desjenigen $\mathbb{R}^{n}$,
der den jeweiligen Zustandsraum $\zeta$ enth"alt, $n-1$
unabh"angige Richtungen existieren, die alle die Richtungen der Gradienten $\nabla\beta_{j}(x)$ sein k"onnen,
sodass es $n-1$ unabh"angige L"osungen $\beta_{j}\in\mathcal{C}^{ 1}(\zeta,\mathbb{R})$ gibt.
\newline
Was aber, wenn die Trajektorie $\Psi(x,\mathbb{R})$ der trajektoriellen Partition $[\Psi]$
so verl"auft, dass sie nicht nur durch $x$, sondern beliebig nahe an $x$ vorbei l"auft? --
Was die Redeweise daf"ur ist, dass 
die Trajektorie $\Psi(x,\mathbb{R})$
durch jede hinreichend kleine Umgebung ${\rm U}$ von $x$ nicht nur durch $x$ 
verl"auft, sondern zusammenh"angende Trajektorienausschnitte $\theta\subset\Psi(x,\mathbb{R})\cap {\rm U}$ 
mit $\theta\not\ni x$ existieren. 
Und zwar unendlich viele, die in $x$ dicht liegen:\newline
Wir k"onnen uns dann gut vorstellen, dass diese Trajektorienausschnitte $\theta\subset\Psi(x,\mathbb{R})\cap {\rm U}$
in der Umgebung ${\rm U}$
eine $\mathcal{C}^{ 1}$-Fl"ache oder eine $\mathcal{C}^{ 1}$-Hyperfl"ache $F_{\Psi}(x)$ 
der Dimension $\nu\geq 2$
quasi aufspannen, in der $x$ liegt und deren Einheitstangenten in $x_{1},x_{2},\dots\in{\rm U}$
$$\partial_{2}\Psi(x_{1},\mathbb{R})/||\partial_{2}\Psi(x_{1},\mathbb{R})||,\ 
\partial_{2}\Psi(x_{2},\mathbb{R})/||\partial_{2}\Psi(x_{2},\mathbb{R})||,\dots$$
sind. Dann m"ussen die Gradienten $\nabla\beta_{j}(x)$ der besagten $\mathcal{C}^{ 1}$-Invarianten $\beta_{j}\in\mathcal{C}^{ 1}(\zeta,\mathbb{R})$
alle im Normalraum der $\mathcal{C}^{ 1}$-Hyperfl"ache $F_{\Psi}(x)$ an der Stelle $x$ liegen.
Da letztere die Dimension $\nu$ hat, kann es h"ochstens $n-\nu$
unabh"angige L"osungen $\beta_{j}\in\mathcal{C}^{ 1}(\zeta,\mathbb{R})$ geben.
Diese suggestive "Uberlegung f"uhrt uns vor Augen, dass die minimalen invarianten Mannigfaltigkeiten nicht die
Trajektorien selbst sein m"ussen. Wir verweisen auf den 
Satz {\bf A.1} des Anhanges.
\newline
Und wir weisen auf das Problem der Bestimmung konkreter 
minimaler invarianter Mannigfaltigkeiten hin, das leider nicht lediglich die Aufgabe ist, eine parzielle Differenzialgleichung zu l"osen.
Schon eine Differenzialgleichung zu l"osen, kann uns mit Unannehmlichkeiten konfrontieren. 
Insofern als die Handhabung dieser Unannehmlichkeiten ein klassischer Problemkreis der mathematischen Physik 
ist, der sich zu einem Schwerpunkt der Numerik vergegenw"artigte, ragt die Frage der Behandlung
der Invariantengleichung (\ref{miymmo}) von anderen Arbeitsgebieten in dieses Traktat.
Wir nehmen an, dass
uns ein Orakel an die Hand gegeben ist,
das uns zu einer vorgelegten $\mathcal{C}^{ 1}$-Flussfunktion $\Psi$
die Menge aller st"uckweise differenzierbarer
L"osungen $\beta_{j}\in\mathcal{C}^{1/2}(\zeta,\mathbb{R})$ der Invariantengleichung (\ref{miymmo})
dieser jeweiligen Flussfunktion $\Psi$ sagt, das uns also
die Menge
\begin{equation}\label{amimmo}
\begin{array}{c}
{\rm Inv}^{q}([\Psi])=\Bigl\{\lambda\in\mathcal{C}^{q}(\mathbf{P}_{2}\Psi,\mathbb{R}):\\
x\in\mathbf{P}_{2}\Psi\Rightarrow
\partial_{2}\Psi(x,0)^{\top}\nabla\lambda(x)=0\Bigr\}
\end{array}
\end{equation}
f"ur beide Kontinuit"atsindizes $q\in\{1/2,1\}$
angibt. Dann sollte es uns doch m"oglich sein,
minimale invariante $\mathcal{C}^{ q}$-Mannigfaltigkeiten 
zu bestimmen. 
Denn f"ur jedes Mengensystem $X$ sind die Mengensysteme
\begin{equation}\label{bmimmo}
\begin{array}{c}
X[1]:=\{\Lambda\subset X:\bigcap\Lambda\not=\emptyset\}\ ,\\
X[2]:=\{\bigcap\Lambda:\Lambda\in X[1]\}\ ,\\
X[3]:=\{\Gamma\in X[2]:\Theta\in X[2]\land\Theta\supset\Gamma\Rightarrow\Gamma=\Theta\}
\end{array}
\end{equation}
wohldefiniert. Deren Konstruktion kommt von dem Mengensystem $X$ 
"uber den ersten Auswahlschritt aller Teilmengensysteme mit nicht leerem Schnitt
und "uber den anschliessenden Auswahlschritt der kleinsten aller dieser jeweiligen nicht leeren Schnitte
zu dem Mengensystem $X[3]$.
Mit ${\rm Inv}^{q}([\Psi])$ ist uns f"ur beide Kontinuit"atsindizes $q\in\{1/2,1\}$
auch das  
Mengensystem
\begin{equation}\label{cmimmo}
\begin{array}{c}
\Bigl\{\lambda^{ -1}(\omega): (\lambda,\omega)\in{\rm Inv}^{q}([\Psi])\times\mathbb{R}\Bigr\}
\end{array}
\end{equation}
zur Hand, die die leere Menge als Element hat, wenn es eine Invariante
$\lambda\in {\rm Inv}^{q}([\Psi])$ gibt, f"ur die 
$\mathbf{P}_{2}\lambda\not=\mathbb{R}$ ist, was aber f"ur die Konstruktion des 
Mengensystemes 
\begin{equation}\label{dmimmo}
\begin{array}{c}
\Bigl\{\lambda^{ -1}(\omega): (\lambda,\omega)\in{\rm Inv}^{q}([\Psi])\times\mathbb{R}\Bigr\}[3]\ =
\mathcal{M}^{q}([\Psi])
\end{array}
\end{equation} gem"ass (\ref{bmimmo})
ohne Relevanz ist. Die Konstruktion des Mengensystemes $X[3]$
zu dem Mengensystem $X$ 
ist offenbar nur in dem Fall, dass $X$ eine endliche Menge endlicher Mengen ist,
programmierbar.
Dann ist $X$
ein "Aquivalent einer endlichen Menge ${\rm X}$
endlicher Mengen der Potenzmenge $2^{\mathbb{N}}$, die 
eine Matrix aus der Menge $$\{0,1\}^{\mathbf{card}({\rm X})\times\max\bigcup{\rm X}}$$
mit $\mathbf{card}({\rm X})\max\bigcup{\rm X}$ bits zu kodieren vermag. 
Wir sehen:
Die Konstruktion des Mengensystemes minimaler nicht leerer Schnitte ${\rm X}[3]$, das eine Matrix aus der Menge
$$\{0,1\}^{\mathbf{card}({\rm X}[3])\times\max\bigcup{\rm X}}$$
kodiert, ist m"oglich. Dabei ist schon in diesem finiten und programmierbaren Fall offenbar ein
$\mathbf{NP}$-Problem gegeben, 
weil $2^{\mathbf{card}({\rm X})}$ Schnitte von Mengen aus ${\rm X}$ gebildet und "uberpr"uft werden m"ussen.
Aber die Repr"asentierbarkeit des kontinuumsm"achtigen Mengensystemes
kontinuumsm"achtiger Mengen (\ref{cmimmo}) durch
eine endliche Menge endlicher Mengen wird durch das m"ogliche beschriebene Ph"anomen,
dass eine durch einen Zustand $x$ verlaufende Trajektorie
beliebig nahe an demselben vorbei verl"auft,
kritisch. Wir sehen: Die Bestimmung minimaler 
invarianter Mannigfaltigkeiten ist mit der L"osung der Invariantengleichung (\ref{miymmo}) 
noch 
nicht im Wesentlichen bereits getan.
\newline
Wir bemerken, dass 
die Menge aller 
eine minimale invariante Mannigfaltigkeit
$m\in\mathcal{M}^{  q}\ ([\Psi])$ schneidenden Trajektorien der trajektoriellen Partition $[\Psi]$
die st"uckweise stetige minimale invariante Mannigfaltigkeit $m$ partioniert: Es ist 
\begin{equation}[\Psi]_{ m}:= 
\Bigl\{\tau\in\ [\Psi]:\tau \cap m\ \not=\emptyset\ \Bigr\}\ \in\ \mathbf{part} (m)\ ,
\end{equation}
wobei 
f"ur jedes Mengensystem\index{Mengensystem} $W$ und $w\subset\bigcup W$
\begin{equation}\label{not1}
W_{ w}:= \{v\in W:\ v\cap w \not= \emptyset\}
\end{equation}
die $W$-Auswahl durch $w$ sei.\index{Notationskonvention} 
Es gilt dar"uber hinaus auch die Implikation
\begin{equation} 
m\ \in\ \mathcal{M}^{0}\ ([\Psi])\ \Rightarrow\ \Bigl\{\mathbf{cl} (\tau):\tau\in [\Psi] \Bigr\}_{ m} \in\  
\mathbf{part}\ (m)\ .
\end{equation}
Offensichtlich existieren trajektorielle Partitionen,
deren
minimale invariante Mannigfaltigkeiten unbeschr"ankt sind.
F"ur diese unbeschr"ankten trajektoriellen Partitionen kann der Fall vorliegen, dass
nur triviale Quasiergodik gegeben ist: Beispielsweise ist jedes diffeomorphe Bild der
trivialen trajektoriellen Partition
$(\mathbb{R}^{ 2},\{(c,\mathbb{R}):c\in \mathbb{R}\})$ eine trajektorielle Partition, f"ur die insofern
nur triviale Quasiergodik gegeben ist, als $[\Psi]$
lauter abgeschlossene Trajektorien als Elemente hat. Ferner existieren offensichtlich auch
beschr"ankte trajektorielle Partitionen, deren Trajektorien lauter abgeschlossene Trajektorien sind.
Letztere sind dann aber geschlossene Trajektorien; solche Trajektorien, die zum Kreis $\mathbb{S}^{1}$ relativ hom"oomorph sind.\index{triviale Quasiergodik}
"Uber das, was wir mit 
relativer Hom"oomorphie\index{relative Hom\"oomorphie} meinen, lassen wir uns hierbei nicht im Unklaren und wir
gestatten uns die folgende exkursive Erl"auterung:\newline
Sind $\mathbf{T}_{1}$ und $\mathbf{T}_{2}$ zwei Topologien und liegen in
deren zugeh"origen topologischen R"aumen
$(\bigcup\mathbf{T}_{1},\mathbf{T}_{1})$ und $(\bigcup\mathbf{T}_{2},\mathbf{T}_{2})$ die Mengen
${\rm A}_{1}\subset \bigcup\mathbf{T}_{1}$ und ${\rm A}_{2}\subset\bigcup\mathbf{T}_{2}$,
so nennen wir diese beiden Mengen ${\rm A}_{1}$ und ${\rm A}_{2}$
genau dann zueinander relativ hom"oomorph bez"uglich der Topologien 
$\mathbf{T}_{1}$ und $\mathbf{T}_{2}$, wenn dieselben bez"uglich
der Relativtopologien $\mathbf{T}_{1}\cap {\rm A}_{1}$ und $\mathbf{T}_{1}\cap {\rm A}_{2}$
zueinander hom"oomorph
sind. 
Exakt dann notieren wir\index{Notationskonvention} 
\begin{equation}\label{eydita}
{\rm A}_{1}\sim_{\mathbf{T}_{1}}\sim_{\mathbf{T}_{2}}{\rm A}_{2} 
\end{equation}
und exakt andernfalls setzen wir ${\rm A}_{1}\not\sim_{\mathbf{T}_{1}}\sim_{\mathbf{T}_{2}}{\rm A}_{2}$,
wohingegen
\begin{equation}\label{eyditb}
{\rm A}_{1}\sim_{\mathbf{T}_{1},\mathbf{T}_{2}}{\rm A}_{2}
\end{equation}
die Hom"oomorphie der Mengen ${\rm A}_{1}$ und ${\rm A}_{2}$
bez"uglich der Topologien 
$\mathbf{T}_{1}$ und $\mathbf{T}_{2}$ notiere und dementsprechend
${\rm A}_{1}\not\sim_{\mathbf{T}_{1},\mathbf{T}_{2}}{\rm A}_{2}$ gerade die 
Negation der Hom"oomorphie (\ref{eyditb}) schriftlich darstelle.
Da das Mengensystem
$$\mathbf{T}\cap {\rm A} =\{{\rm U}\cap {\rm A}:{\rm U}\in \mathbf{T}\}$$
$$=\Bigl\{{\rm U}\cap ({\rm A}\cap\bigcup\mathbf{T}):{\rm U}\in \mathbf{T}\Bigr\}=\mathbf{T}\cap ({\rm A}\cap\bigcup\mathbf{T})$$
ist, k"onnen wir ohne Einschr"ankung der Allgemeinheit von der relevanten 
Teilmenge ${\rm A}\cap\bigcup\mathbf{T}\subset {\rm A}$ ausgehen, die
in dem topologischen Raum $(\bigcup\mathbf{T},\mathbf{T})$ liegt.\footnote{Allemal: Denn, selbst wenn
${\rm A}$ und $\bigcup\mathbf{T}$ disjunkt sind und in diesem Disjunktionsfall ${\rm A}\cap\bigcup\mathbf{T}=\emptyset$
ist, gilt die Inklusion
${\rm A}\cap\bigcup\mathbf{T}=\emptyset\subset \bigcup\mathbf{T}$. Der Disjunktionsfall liegt genau dann vor, wenn
die 
Identit"at  
$$\mathbf{T}\cap {\rm A}=\{\emptyset\}$$
zutrifft; sodass 
das Mengensystem $\mathbf{T}\cap {\rm A}$
exakt im Disjunktionsfall keine Topologie ist, falls wir den Begriff der Topologie so verfassen, wie er 
am weitesten verbreitet ist, n"amlich so, dass f"ur jede Topologie $\mathbf{T}$
die Reichhaltigkeitsaussage $\bigcup\mathbf{T}\not=\emptyset$ gelte. 
Verlangten wir neben der finiten Schnittabgeschlossenheit einer Topologie
und deren Abgeschlossenheit gegen"uber Vereinigungen, dass 
lediglich die abgeschw"achte Reichhaltigkeitsaussage
$\mathbf{T}\not=\emptyset$ gelten solle, so w"are das Mengensystem
$\mathbf{T}\cap {\rm A}=\{\emptyset\}$ auch im besagten Disjunktionsfall eine Topologie.
Halten wir uns hier aber an die Majorit"at und befinden, dass 
exakt im besagten Disjunktionsfall, exakt in dem ${\rm A}\cap\bigcup\mathbf{T}=\emptyset$ ist,
das Mengensystem $\mathbf{T}\cap {\rm A}$ keine Topologie ist!
Im Fall, dass ${\rm A}$ und $\bigcup\mathbf{T}$ einander schneiden, ist das
Mengensystem
$\mathbf{T}\cap {\rm A}$ eine Topologie, exakt welche bekanntlich als die Spurtopologie oder die 
Relativtopologie\index{Relativtopologie}\index{Spurtopologie}
der Topologie $\mathbf{T}$
bez"uglich ${\rm A}$
bezeichnet wird. }
Die relative Hom"oomorphie ist insofern schw"acher als die 
Hom"oomorphie, als die Implikation
\begin{equation}\label{eyditc}
{\rm A}_{1}\sim_{\mathbf{T}_{1},\mathbf{T}_{2}}{\rm A}_{2}\ \Rightarrow\ 
{\rm A}_{1}\sim_{\mathbf{T}_{1}}\sim_{\mathbf{T}_{2}}{\rm A}_{2}
\end{equation}
wahr ist, die aber im Allgemeinen nicht umkehrbar ist.
Die relative Hom"oomorphie ist dabei in dem Sinn eine universelle "Aquivalenzrelation, als die Aussagen
\begin{equation}\label{eyditd}
\begin{array}{c}
{\rm A}_{1}\sim_{\mathbf{T}_{1}}\sim_{\mathbf{T}_{1}}{\rm A}_{1}\ ,\\
{\rm A}_{1}\sim_{\mathbf{T}_{1}}\sim_{\mathbf{T}_{2}}{\rm A}_{2}\Leftrightarrow
{\rm A}_{2}\sim_{\mathbf{T}_{2}}\sim_{\mathbf{T}_{1}}{\rm A}_{1}\ ,\\
{\rm A}_{1}\sim_{\mathbf{T}_{1}}\sim_{\mathbf{T}_{2}}{\rm A}_{2}\ \land\ {\rm A}_{2}\sim_{\mathbf{T}_{2}}\sim_{\mathbf{T}_{3}}{\rm A}_{3}\\
\Rightarrow
{\rm A}_{1}\sim_{\mathbf{T}_{1}}\sim_{\mathbf{T}_{3}}{\rm A}_{3}
\end{array}
\end{equation} f"ur alle 
Topologien $\mathbf{T}_{1}$ und $\mathbf{T}_{2}$ und alle Mengen
${\rm A}_{1}\subset\bigcup\mathbf{T}_{1}$ und ${\rm A}_{2}\subset\bigcup\mathbf{T}_{2}$
wahr sind.\newline
Nun redeten wir aber soeben von zum Kreis $\mathbb{S}^{1}$ relativ hom"oomorphen Trajektorien,
ohne die Referenz dieser besprochenen relativen Hom"oomorphie ausdr"ucklich festzulegen. Wir sagten ja nicht,
bez"uglich welcher Topologie des Raumes, in dem die jeweilige Trajektorie liegt 
und bez"uglich welcher Topologie des Raumes, in dem der Kreis $\mathbb{S}^{1}$ liegt, die jeweilige Trajektorie zu
dem Kreis $\mathbb{S}^{1}$
relativ hom"oomorph sei. Wenn wir den Kreis $\mathbb{S}^{1}$ kurzerhand in der Ebene situieren, die
der $\mathbb{R}^{2}$ objektiviert, 
dann greifen wir mit diesem allzu naheliegenden Griff auch zu kurz. Inwiefern wir da zu kurz greifen und f"ur welchen Zweck
wir da zu kurz greifen?\newline
Da bekanntlich kein Hom"oomorphismus existiert,
der den $\mathbb{R}^{2}$ auf den $\mathbb{R}^{k}$ abbildet, wenn $k\in\mathbb{N}\setminus\{2\}$ ist, ist
das freie Reden "uber die Hom"oomorphie von Punktmengen, die im $\mathbb{R}^{n}$ f"ur eine beliebige 
Dimension $n\in\mathbb{N}\setminus\{2\}$ liegen,
nicht auf die Weise m"oglich, die uns die folgende Festlegung erlaubt:
Es sei 
$$\mathbb{S}^{1}(\mathbb{R}^{2}):=\{x\in\mathbb{R}^{2}:||x||_{\mathbb{R}^{2}}=1\}$$ und f"ur alle 
$n\in\mathbb{N}\setminus\{1,2\}$ sei
$$\mathbb{S}^{1}(\mathbb{R}^{n}):= \mathbb{S}^{1}(\mathbb{R}^{2})\oplus (\delta_{1k})_{3\leq j\leq n}$$ 
der in den $\mathbb{R}^{n}$ eingebettete Kreis,
wobei $||\mbox{id}||_{\mathbb{R}^{2}}$ die euklidische Norm des $\mathbb{R}^{2}$ sei. Ferner sei ${\rm A}$ eine Punktmenge, die auf die Weise homogen sei, dass es 
eine nat"urliche Zahl
$n({\rm A})\in \mathbb{N}$
gibt, f"ur die
${\rm A}\subset\mathbb{R}^{n({\rm A})}$ gilt.
Wenn wir nun "uber eine solche, in dieser Form homogene Punktmenge ${\rm A}$
sagen, dass sie
zu dem Kreis $\mathbb{S}^{1}$ hom"oomorph bzw. diffeomorph bzw. relativ hom"oomorph sei,
dann meinen wir, dass ${\rm A}$ 
zu der passenden Einbettung des ebenen Kreises $\mathbb{S}^{1}(\mathbb{R}^{2})$ in den $\mathbb{R}^{n({\rm A})}$, also
zu der objektaktualisierten Menge $$\mathbb{S}^{1}(\mathbb{R}^{n({\rm A})})\subset\mathbb{R}^{n({\rm A})}$$ 
hom"oomorph bzw. diffeomorph bzw. relativ hom"oomorph sei.
Die euklidische Norm $||\mbox{id}||_{\mathbb{R}^{2}}$ legt dabei die euklidische Topologie des $\mathbb{R}^{2}$ fest, die
wir auch als die nat"urliche Topologie des $\mathbb{R}^{2}$ bezeichnen. Unsere Darlegungen beschr"anken sich in dieser
Abhandlung auf Zustandsr"aume ${\rm Z}$, die in reellen R"aumen $\mathbb{R}^{n}$
endlicher Dimension $n\in\mathbb{N}$ liegen, wobei ${\rm Z}$ durch deren jeweilige nat"urliche Topologie $\mathbf{T}(n)$ topologisiert sei.\newline 
Da die 
Trajektorien $\tau$, von denen hier die Rede ist, immer solche sind, die in einem reellen Raum $\mathbb{R}^{n}$
endlicher Dimension $n\in\mathbb{N}$ liegen,
und da
die nat"urliche Topologie $\mathbf{T}(n(\tau))$
desjenigen Raumes $\mathbb{R}^{n(\tau)}$,
in dem die jeweilige Trajektorie $\tau\subset: \mathbb{R}^{n(\tau)}$ liegt, 
mit demselben durch die jeweilige Trajektorie $\tau$ festgelegt ist,
k"onnen wir uns hier ersparen, die 
Referenztopologien
ausdr"ucklich festzuhalten,
wenn wir "uber die relative Hom"oomorphie einer Trajektorie zu einer anderen
Menge eines topologischen 
Raumes sprechen: Es ist $\mathbf{T}(n(\tau))$ die jeweilige Topologie, auf die
sich hier gegebenfalls eine jede Aussage "uber die relative Hom"oomorphie einer Trajektorie $\tau$ zu einer anderen
Menge eines topologischen 
Raumes bezieht.
\newline
Wir legen vor diesem Hintergrund der thematischen Beschr"ankung 
dieser Abhandlung f"ur alle Punktmengen ${\rm A}\subset\mathbb{R}^{n({\rm A})}$ und 
${\rm B}\subset\mathbb{R}^{n({\rm B})}$ 
reeller R"aume endlicher Dimension $n({\rm A}),n({\rm B})\in\mathbb{N}$ 
die folgende Schreibweise fest:\index{Notationskonvention} Es gelte
\begin{equation}\label{eydite}
\begin{array}{c}
{\rm A}\sim_{\mathbf{T}(n({\rm A}))}\sim_{\mathbf{T}(n({\rm B}))}{\rm B}\ \Leftrightarrow:
{\rm A}\sim\sim{\rm B}\ ,\\
{\rm A}\not\sim_{\mathbf{T}(n({\rm A}))}\sim_{\mathbf{T}(n({\rm B}))}{\rm B}\ \Leftrightarrow:
{\rm A}\not\sim\sim{\rm B}
\end{array}
\end{equation}
und entsprechend
\begin{equation}\label{eyditf}
\begin{array}{c}
{\rm A}\sim_{\mathbf{T}(n({\rm A})),\mathbf{T}(n({\rm B}))}{\rm B}\ \Leftrightarrow:
{\rm A}\sim{\rm B}\ , \\
{\rm A}\not\sim_{\mathbf{T}(n({\rm A})),\mathbf{T}(n({\rm B}))}{\rm B}\ \Leftrightarrow:
{\rm A}\not\sim{\rm B}\ 
\end{array}
\end{equation}
und schliesslich auch
\begin{equation}\label{eyditz}
\begin{array}{c}
{\rm A}\sim\mathbb{S}^{1}(\mathbb{R}^{n({\rm A})})\ \Leftrightarrow:\ {\rm A}\sim\mathbb{S}^{1}\ ,\\
{\rm A}\not\sim\mathbb{S}^{1}(\mathbb{R}^{n({\rm A})})\ \Leftrightarrow:\ {\rm A}\not\sim\mathbb{S}^{1}\ ,\\
{\rm A}\sim\sim\mathbb{S}^{1}(\mathbb{R}^{n({\rm A})})\ \Leftrightarrow:\ {\rm A}\sim\sim\mathbb{S}^{1}\ ,\\
{\rm A}\not\sim\sim\mathbb{S}^{1}(\mathbb{R}^{n({\rm A})})\ \Leftrightarrow:\ {\rm A}\not\sim\sim\mathbb{S}^{1}\ .\\
\end{array}
\end{equation}
Da f"ur stetige Flussfunktionen $\Psi$
mit kompaktem Zustandsraum $\zeta$
die Quasiperiodizit"at bzw. die Periodizit"at eines Zustandes $z\in\zeta$ bez"uglich $\Psi$
und die Quasizyklizit"at bzw. die Zyklizit"at der Trajektorie $\Psi(z,\mathbb{R})$ zueinander "aquivalent sind,
ist auch die Geschlossenheit einer Trajektorie $\tau\in [\Psi]$ zu deren
relativer Hom"oomorphie zum Kreis "aquivalent: Denn es gelten f"ur jede stetige Flussfunktion $\Psi$
mit kompaktem Zustandsraum $\zeta\subset\mathbb{R}^{n}$ und f"ur jeden Zustand $z\in\zeta$ die "Aquivalenzen
\begin{equation}\label{ayditz}
\begin{array}{c}
\Psi(z,\mathbb{R})=\mathbf{cl}(\Psi(z,\mathbb{R}))\Leftrightarrow\\
\exists\ T\in\mathbb{R}\setminus\{0\}:\  \Psi(z,T)=z\ ,\\
\quad\\
\exists\ T\in\mathbb{R}\setminus\{0\}:\  \Psi(z,T)=z\Leftrightarrow\\ 
\Psi(z,]0,T])=\Psi(z,\mathbb{R})\ \land\ \Psi(z,0)=\Psi(z,T)\\ 
\quad\\
\Psi(z,]0,T])=\Psi(z,\mathbb{R})\ \land\ \Psi(z,0)=\Psi(z,T)\Leftrightarrow\\
\Psi(z,]0,T])\sim\sim\mathbb{S}^{1}\ ,
\end{array}
\end{equation}
denn 
die Restriktionen
$$\Psi(z,\mbox{id})|]0,T]\quad {\rm bzw.}\quad (\Psi(z,\mbox{id})|]0,T])^{-1}$$
sind stetig bez"uglich $\mathbf{T}(1)\cap ]0,T]$ und 
$\mathbf{T}(n)\cap \Psi(z,\mathbb{R})$ bzw. bez"uglich $\mathbf{T}(n)\cap \Psi(z,\mathbb{R})$ und $\mathbf{T}(1)\cap ]0,T]$.
\newline
Dass die relative Hom"oomorphie, wie wir mit der Implikation (\ref{eyditc}) bereits behaupteten, schw"acher als die unrelativierte 
Hom"oomorphie ist, wird uns im $\mathbb{R}^{3}$
anschaulich: Ist ${\rm K}\subset \mathbb{R}^{3}$ eine Punktmengenrepr"asentation des bekannten 
Kleeblattknotens\index{Kleeblattknoten} im $\mathbb{R}^{3}$, so gilt offensichtlich der Sachverhalt, den wir
als die Aussage
\begin{equation}\label{eyditg}
{\rm K}\sim\sim \mathbb{S}^{1}\ \land\ {\rm K}\not\sim \mathbb{S}^{1}
\end{equation}
notieren.
Unsere Untersuchungen werden ergeben, dass der kontinuierliche Determinismus,
den stetige Flussfunktionen gem"ass (\ref{kurzc})
objektivieren, 
die Alternative bedingt, die das Hauptthema dieser Abhandlung ist,
wenn kompakte Zustandsr"aume endlichdimensionaler reeller R"aume vorliegen:
Die 
abgeschlossenen H"ullen jeweiliger Trajektorien, die mit den jeweiligen
minimalen invarianten Mannigfaltigkeiten zusammenfallen,
sind dann entweder 
sensitive Attraktoren 
oder aber Fixpunkte oder geschlossene Trajektorien, solche Trajektorien also, die 
nach unserer nunmehr erl"auterten Redeweise 
zum Kreis $\mathbb{S}^{1}$ relativ hom"oomorph sind.
Deswegen kommt dem Fall, dass die minimalen invarianten Mannigfaltigkeiten
keine sensitive Attraktoren oder Fixpunkte, sondern zum Kreis $\mathbb{S}^{1}$ relativ hom"oomorph sind, 
einige Wichtigkeit zu, weswegen
wir es nicht unterliessen, die relative Hom"oomorphie im Exkurs der letzten 
Seiten
zu beschreiben und dieselbe gegen die unrelativierte Hom"oomorphie abzugrenzen,
die eine st"arkere Eigenschaft ist.\newline Die herausgestellte Differenz von
relativer Hom"oomorphie und unrelativierter Hom"oomorphie 
leitet uns naturgem"ass zu der Frage,
unter welchen Umst"anden die Kontinuit"at stetiger Flussfunktionen gem"ass (\ref{kurzc}) "uberhaupt 
zul"asst, dass jene herausgestellte Differenz
durch die minimalen invarianten Mannigfaltigkeiten realisiert sein kann, die
im Rahmen des durch jene Flussfunktionen objektivierten kontinuierlichen Determinismus m"oglich sind:
Unter welchen Umst"anden existieren f"ur stetige Flussfunktionen $\Psi$ gem"ass (\ref{kurzc}) "uberhaupt Zust"ande $x$, sodass
geschlossene Trajektorien $\Psi(x,\mathbb{R})$ auftreten, die zwar
zum Kreis $\mathbb{S}^{1}$ nicht
hom"oomorph sind, die zu demselben nichtsdestotrotz
relativ hom"oomorph sind?
Diese Frage k"onnen wir anschaulich kurzfassen als die Frage,
unter welchen Umst"anden der kontinuierliche Determinismus zul"asst,
dass sich eine Trajektorie verknotet.
\newline
Die folgende Frage k"onnen wir nicht so einfach beantworten, wie die Frage,
{\em ob} der kontinuierliche Determinismus zul"asst,
dass sich eine Trajektorie verknotet. Denn beispielsweise f"ur jene bereits vorgestellte 
Konkretisierung des Kleeblattknotens als eine 
Punktmenge ${\rm K}\subset \mathbb{R}^{3}$
gibt es ja offenbar
eine stetige Flussfunktion $\Psi_{{\rm K}}$, deren Zustandsraum
${\rm K}$ ist und die die einzige Trajektorie
$${\rm K}=\Psi_{{\rm K}}(x,\mathbb{R})$$
f"ur alle Zust"ande $x\in {\rm K}$ hat.
Gibt es aber eine Fortsetzung  
\begin{displaymath}
\begin{array}{c}
\Psi_{{\rm K}}^{+}:{\rm Z}\times \mathbb{R}\to {\rm Z}, \\
(z,t)\mapsto \Psi_{{\rm K}}^{+}(z,t) 
\end{array}
\end{displaymath} der Flussfunktion $\Psi_{{\rm K}}$, die ebenfalls
gem"ass (\ref{kurzc}) beschaffen ist und
f"ur die 
$${\rm K}=\Psi_{{\rm K}}^{+}(x,\mathbb{R})$$
f"ur alle Zust"ande $x\in {\rm K}$ gilt und 
deren Zustandsraum ${\rm Z}$ "uberdies eine offene Kugel $\mathbb{B}\subset{\rm Z}$ des $\mathbb{R}^{3}$
enth"alt, die ihrerseits die Repr"asentation des Kleeblattknotens ${\rm K}$
enth"alt?\newline 
Wir stossen also auf den nicht sch"arfer umrissenen Fragenkreis, 
unter welchen Umst"anden f"ur stetige Flussfunktionen $\Psi$ gem"ass (\ref{kurzc}) "uberhaupt Zust"ande $x$ existieren, sodass
geschlossene Trajektorien $\Psi(x,\mathbb{R})$ auftreten, die 
zum Kreis $\mathbb{S}^{1}$ nicht
hom"oomorph aber relativ hom"oomorph sind;
und indem wir die mit diesem Fragenkreis verbundene Thematik hier nicht eingehender behandeln,
stellen wir dieses Themenfeld als eine Aufgabe und als ein Arbeitsgebiet heraus.
Auf Letzteres wollen wir uns aber
nicht in dieser Abhandlung begeben, deren Thema eng lokalisiert ist. Letzteres ist ja \glqq ein gesuchter, dennoch bislang "ubersehener elementarer Satz\grqq ,
n"amlich der Satz 1.1, der elementare Quasiergodensatz, den wir nun formulieren wollen:
\subsection{Die Formulierung des elementaren Quasiergodensatzes}
Da die minimalen invarianten Mannigfaltigkeiten 
jeder trajektoriellen Partition Teilmengen 
des Zustandsraumes sind,
gilt f"ur kompakte Zustandsr"aume $\zeta$ kompakter trajektorieller Partitionen $[\Psi]$,
dass $m\in \mathcal{M}^{ 1}([\Psi])$ beschr"ankt ist.
Dabei nennen wir exakt trajektorielle Partitionen mit kompakten Zustandsr"aumen
kompakt.\index{kompakte trajektorielle Partition}
Wenn $[\Psi]$ eine kompakte trajektorielle $\mathcal{C}^{1}$-Partition ist,
ist $$\mathcal{M}^{ 0}([\Psi])=\mathcal{M}^{1/2} ([\Psi])\ ,$$
wie wir sp"ater zeigen werden.
F"ur kompakte trajektorielle Partitionen sind die Geschlossenheit ihrer Trajektorien oder die Mehrdimensionalit"at deren abgeschlossener
H"ullen einander ausschliessende Alternativen. Und nur die irrige Vorstellung, dass
die minimalen invarianten Mannigfaltigkeiten 
einfach die Trajektorien sind, l"asst den Sachverhalt dieser Alternative paradox erscheinen.
Aber die Konstruktion minimaler invarianter Mannigfaltigkeiten ist
keine, die Trajektorien lediglich auf eine andere Weise beschreibt. 
Wir erl"auterten die Differenz zwischen 
minimalen invarianten Mannigfaltigkeiten
und Trajektorien bereits und verweisen 
wiederholt
auf den Satz {\bf A.1} des Anhanges. Minimale invariante Mannigfaltigkeiten k"onnen h"oher als eindimensional sein.
Dies aber sahen wir nun schon ein:
Nur innerhalb minimaler invarianter Mannigfaltigkeiten kann Quasiergodik
vorliegen, weil jede Trajektorie exakt innerhalb einer solchen verl"auft.
Sei beispielsweise eine minimale invariante $\mathcal{C}^{1}$-Mannigfaltigkeit
im Fall einer Vielteilchen-Interpretation\index{Vielteilchen-Interpretation} 
einer jeweiligen, hinreichend kontinuierlichen Flussfunktion $\Psi$
nicht mit der entsprechenden Energiehyperfl"ache identisch, was
dann der Fall ist, wenn die Energie nicht die einzige Invariante dieser Flussfunktion $\Psi$ ist:
Dann kann keine Quasiergodik innerhalb der gesamten Energiehyperfl"ache vorliegen.
Die Quasiergodik innerhalb der gesamten Energiehyperfl"ache gibt es nur, wenn die 
Energie die einzige Invariante ist und
es zeichnet sich ab, dass dies nur in speziell konstruierten F"allen gegeben ist. 
Bei einem kompakten Zustandsraum m"usste, damit die Energie die einzige Invariante ist, die gesamte Energiehyperfl"ache\index{Energiehyperfl"ache}
ein sensitiver Attraktor\index{sensitiver Attraktor} sein, wie wir noch sehen werden. Wir verweisen auf den zweiten Abschnitt des folgenden Kapitels.
\newline
Best"unde die Differenz zwischen 
minimalen invarianten Mannigfaltigkeiten
und Trajektorien nicht, so w"are unsere folgende
Modellierung der Quasiergodenhypothese nicht zweckm"assig.
Die Behauptung des folgenden Satzes setzen 
wir n"amlich nun als die Objektivierung der Quasiergodenhypothese im Sinne P. und T. Ehrenfests an, die
die minimalinvasive Abwandlung von Boltzmanns 
unhaltbarer
Ergodenhypothese\index{Ergodenhypothese}
ist:
\newline
\newline
{\bf Elementarer Quasiergodensatz 1.1: }\index{elementarer Quasiergodensatz}
\newline
{\em Alle Trajektorien $\tau\in [\Psi]$ einer st"uckweise glatten minimalen invarianten Mannigfaltigkeit 
$\mu \in \mathcal{M}^{1/2} ([\Psi])$ einer kompakten trajektoriellen $\mathcal{C}^{1}$-Partition $[\Psi]$ liegen 
jedem Punkt der Mannigfaltigkeit $\mu$ beliebig nahe.
Es gilt f"ur alle $\tau\in [\Psi]$ und f"ur alle $\mu \in\ \mathcal{M}^{1/2} ([\Psi])$ die "Aquivalenz}
\begin{equation}\label{ergo}
\tau \subset \mu\ \Leftrightarrow\  \mathbf{cl} (\tau) = \mu\ . 
\end{equation}
\newline 
St"uckweise glatte minimale invariante Mannigfaltigkeiten $\mu \in \mathcal{M}^{1/2} ([\Psi])$
k"onnen wir als heterogene Niveauliniensysteme\index{Niveauliniensystem} von Tupeln st"uckweise differenzierbarer
Invarianter der trajektoriellen Partition $[\Psi]$ auffassen.
Die st"uckweise glatten minimalen invarianten Mannigfaltigkeiten $\mu \in \mathcal{M}^{1/2} ([\Psi])$
partionieren den Zustandsraum $\mathbf{P}_{2} \Psi$ per Konstruktionem, sodass es eine Parzialbehauptung des 
elementaren Quasiergodensatzes ist, dass
das Mengensystem
\begin{equation}\label{erddo}
[[\Psi]]:=\{ \mathbf{cl} (\tau):\tau\in [\Psi]\}\ \in\ \mathbf{part}(\mathbf{P}_{2} \Psi)
\end{equation}
den Zustandsraum $\mathbf{P}_{2} \Psi$ partioniert. Gerade diese Parzialbehauptung des 
elementaren Quasiergodensatzes ist die Aussage des Satzes von der Existenz der Zimmer 2.1.2,\index{Satz von der Existenz der Zimmer} die wir
im Abschnitt 2.1 beweisen. Diese Parzialbehauptung ist insofern eine hom"oomorphe Invariante, als f"ur jeden
Hom"oomorphismus $h$ mit der Definitionsmenge $\mathbf{P}_{1}h=\mathbf{P}_{2} \Psi$ die kotransformierte Aussage
\begin{equation}\label{erddoa}
[[h\circ\Psi]]_{\mathbf{T}(\mathbf{P}_{2}h)}:=\{ \mathbf{cl}_{\mathbf{T}(\mathbf{P}_{2}h)} (\tau):\tau\in [h\circ\Psi]\}\ \in\ \mathbf{part}(\mathbf{P}_{2} h)
\end{equation}
gilt, wobei $\mathbf{T}(\mathbf{P}_{2}h)$ eine Topologie der Wertemenge $\mathbf{P}_{2}h$ dieses  
Hom"oomorphismus $h$ sei; und f"ur alle Teilmengen $X\subset{\rm A}$ sei f"ur eine beliebige 
Topologie $\mathbf{T}({\rm A})$, die eine beliebige Menge ${\rm A}$ topologisiert
\begin{equation}\label{erddaa}
\mathbf{cl}_{\mathbf{T}({\rm A})}(X):=\bigcap\Bigl\{{\rm U}\in \mathbf{T}({\rm A}):X\subset{\rm U}\}\ ,
\end{equation}
sodass $\mathbf{cl}_{\mathbf{T}(\mathbf{P}_{2}h)}$ der transformierte H"ullenoperator ist.\index{transformierter H\"ullenoperator}
F"ur jede Menge ${\rm A}$ und f"ur jede Funktion $\Phi\in{\rm A}^{{\rm A}\times\mathbb{R}}$ und f"ur jede Topologie
$\mathbf{T}({\rm A})$
sind die Mengensysteme
\begin{equation}\label{wergo}
\begin{array}{c}
\lbrack\Psi\rbrack:=\{\Psi(x,\mathbb{R}):x\in{\rm A}\}\ ,\\
\lbrack\lbrack\Psi\rbrack\rbrack_{\mathbf{T}({\rm A})}:=\{\mathbf{cl}_{\mathbf{T}({\rm A})}(\tau):\tau\in\lbrack\Psi\rbrack\}
\end{array}
\end{equation}
festgelegt, wobei $\lbrack\lbrack\Psi\rbrack\rbrack$ einfach $[[\Psi]]_{\mathbf{T}({\rm A})}$ sei, wenn $\mathbf{T}({\rm A})$ die
nat"urliche Topologie desjenigen $\mathbb{R}^{n}$ nat"urlicher Dimension $n\in \mathbb{N}$ sei, in dem ${\rm A}$ liegt, sofern ${\rm A}$
so beschaffen ist, dass
dieser $\mathbb{R}^{n}$ existiert. Daher k"onnen wir
 \begin{equation}\label{vergo}
\begin{array}{c}
\lbrack\mbox{id}\rbrack:=\{\mbox{id}(x,\mathbb{R}):x\in\mathbf{P}_{2}\mbox{id}\}\ ,\\
\lbrack\lbrack\mbox{id}\rbrack\rbrack_{\mathbf{T}}:=\{\mathbf{cl}_{\mathbf{T}}(\tau):\tau\in\lbrack\mbox{id}\rbrack\}\ ,\\
\lbrack\lbrack\mbox{id}\rbrack\rbrack:=\{\mathbf{cl}_{\mathbf{T}(\mathbf{P}_{2}\mbox{id})}(\tau):\tau\in\lbrack\mbox{id}\rbrack\}
\end{array}
\end{equation}
als Operatoren oder als Abbildungen der Klasse derjeniger Funktionen $f$ auffassen, f"ur die es eine Menge ${\rm A}$ gibt, 
sodass
$f$ in der Menge ${\rm A}^{{\rm A}\times\mathbb{R}}$ ist. Dabei muss $\mathbf{T}$ eine jeweils angegebene Topologie sein, die 
die jeweilige Wertenmenge 
$\mathbf{P}_{2}\mbox{id}$ topologisiert; oder aber es muss $\mathbf{P}_{2}\mbox{id}$ in 
einem $\mathbb{R}^{n}$ liegen, dessen nat"urliche Topologie dann $\mathbf{T}(\mathbf{P}_{2}\mbox{id})$ sei. Im Fall einer
Heine-Descartesschen Kollektivierung $[\Psi]$ setzt sich die Zyklizit"atsbedingung an jeweilige  Trajektorien $\tau$, welche
der Zwischenwertsatz\index{Zwischenwertsatz} als das In-Einsfallen von Trajektorie $\tau$ und deren H"ulle $\mathbf{cl}(\tau)$ gem"ass der Unreichhaltigkeit
$$\mathbf{cl}(\tau)\setminus\tau=\emptyset$$
zu formulieren erlaubt, in der Form um, dass
$[[\Psi]]\setminus[\Psi]$ gerade das Mengensystem der entsprechenden nicht-trivialen Zimmer ist und dass genau dann 
$$[[\Psi]]\setminus[\Psi]=\emptyset$$
ist, wenn alle Trajektorien der Heine-Descartesschen Kollektivierung $[\Psi]$ Zyklen sind. 
\newline
Bei dieser Gelegenheit k"onnen wir sogleich kl"aren, was wir ganz allgemein als eine Flussfunktion 
und als eine Wellenfunktion\index{Wellenfunktion} auffassen:
Jede Abbildung der Klasse derjeniger Funktionen $f$, f"ur die es eine Menge ${\rm A}$ gibt, 
sodass
$f$ in der Menge ${\rm A}^{{\rm A}\times\mathbb{R}}$ ist, 
d.h., jede Funktion $f$, f"ur die
\begin{equation}\label{viiergo}
\mathbf{P}_{2}f=\mathbf{P}_{1}\mathbf{P}_{1}f\ \land\ \mathbf{P}_{2}\mathbf{P}_{1}f=\mathbb{R}
\end{equation}
gilt,
nennen wir eine Wellenfunktion.
Jede Wellenfunktion 
nennen wir
genau dann eine Flussfunktion, wenn
\begin{equation}\label{viergo}
[f]\in\mathbf{part}(\mathbf{P}_{2}f) 
\end{equation}
ist; und genau dann sagen wir, dass $f$ eine Flussfunktion binnen $\mathbf{P}_{2}f$ sei.\index{Flussfunktion}
\index{Flussfunktion im weiteren Sinn binnen einer Menge}
\newline
Es ist diese (\ref{erddoa}) gem"asse
topologische Invarianz
der Aussage des Satzes von der Existenz der Zimmer 2.1.2, die denselben allgemein und
zum wohl bedeutendsten Gesichtspunkt des elementaren Quasiergodensatzes macht. Suggestiv gesprochen: Auch, wenn wir
das durch eine kompakte trajektorielle $\mathcal{C}^{1}$-Partition $[\Psi]$ beschriebene 
Szenario mittels der Transformation desselben durch einen Hom"oomorphismus quasi zerkn"ullen,
so bleibt die 
von dem Satz von der Existenz der Zimmer behauptete Parzialbehauptung des 
elementaren Quasiergodensatzes, die
Partitivit"at gem"ass (\ref{erddoa}) g"ultig.
\newline
Offenbar wird allein schon durch den 
Satz von der Existenz der Zimmer 2.1.2 bereits die Ergodenfrage 
kontinuierlicher dreidimensionaler Billards gekl"art.
Unter einem kontinuierlichen dreidimensionalen Billard verstehen wir dabei ein  
Punktteilchensystem von ${\rm N}\in\mathbb{N}$ Punktteilchen, all deren Ortskoordinaten innerhalb eines
zusammenh"angenden Kompaktums ${\rm R}\subset\mathbb{R}^{3}$ verbleiben, an dessen Rand 
$\partial {\rm R}$
sie auf kontinuierliche und 
beschr"ankte Weise 
in dieses Kompaktum
${\rm R}$
zur"uckgestreut
werden. Dabei streuen die Punktteilchen auch aneinander, ebenfalls auf kontinuierliche und 
beschr"ankte Weise; was heissen
soll, dass sowohl die Streuvorg"ange aneinander als auch jene R"uckstreuprozesse so
ablaufen, dass
sich die Teilchengeschwindigkeiten stetig "andern und dass dabei deren   
Betr"age immer unterhalb einer f"ur das 
jeweilige Billard
vorgegebenen Grenze $V\in \mathbb{R}$ liegen.\index{kontinuierliches Billard} 
Die kontinuierlichen dreidimensionalen Billards
umfassen die kontinuierlichen zweidimensionalen Billards, die speziell so beschaffen sind, dass
die Ortskoordinaten aller Punktteilchen innerhalb eines
zusammenh"angenden Kompaktums ${\rm R}\subset\mathbb{R}^{3}$ sind, 
das in einer Ebene liegt. Wir bemerken hierbei, dass im Hinblick auf den 
Satz von der Existenz der Zimmer 2.1.2 die Komponentenzahl der 
Ortskoordinaten der Punktteilchen keine Rolle spielt.
Wir k"onnen genausogut ein $m$-dimensionales
kontinuierliches Billard
betrachten, bei dem ${\rm N}\in\mathbb{N}$ Quasiteilchen aneinander oder
am Rand $\partial {\rm R}$
eines
zusammenh"angenden Kompaktums ${\rm R}\subset\mathbb{R}^{m}$
auf kontinuierliche und 
beschr"ankte Weise streuen.\index{multidimensionales kontinuierliches Billard}
Die Quasiteilchen sind dabei durch $m$-Tupel $r\in {\rm R}\subset\mathbb{R}^{m}$ dargestellt und 
$m$
ist dabei eine nat"urliche Zahl. Diese Quasiteilchen sind dann nicht als physische Massenpunkte 
interpretiert.
\newline
Wie kl"art nun schon der 
Satz von der Existenz der Zimmer die Ergodenfrage 
kontinuierlicher multidimensionaler Billards? 
Der Zustandsraum $\zeta\subset\mathbb{R}^{2mN}$ eines kontinuierlichen $m$-dimensionalen Billards
von $N$ Quasiteilchen
ist 
als ein kompakter Zustandsraum formulierbar,
in dem der jeweilige Ortsraum, der im $\mathbb{R}^{mN}$ liegt,
eingebettet ist; etwa auf die Weise, dass 
$$(\mathbf{P}_{1},\mathbf{P}_{2},\dots\mathbf{P}_{mN})\zeta=\underbrace{{\rm R}\oplus{\rm R}\dots\oplus{\rm R}}_{N\ 
\mbox{Kokatenationen}}\subset\mathbb{R}^{mN}$$
der zum kartesischen Produkt ${\rm R}^{N}$ isomorphe Ortsraum ist.
Die den Zustandsraum $\zeta$ partionierenden Trajektorien sind
stetig. (Die Stetigkeit der Trajektorien soll dabei dies besagen: Wenn $\Psi\in \zeta^{\zeta\times \mathbb{R}}$ die
Flussfunktion des jeweiligen kontinuierlichen multidimensionalen Billards ist, sodass
f"ur jeden Zustand $x\in \zeta$ die Punktmenge
$\Psi(x,\mathbb{R})$ 
eine Trajektorie 
darstellt,
so ist jede Trajektorie $\Psi(x,\mathbb{R})$ die Wertemenge der stetigen Funktion $\Psi(x,\mbox{id})$,
die auf dem Zahlenstrahl $\mathbb{R}$ definiert ist.) Die mehr als eindimensionalen abgeschlossenen H"ullen, die mehr als eindimensionalen Zimmer,\index{Zimmer}
sind gem"ass dem Satz von deren Existenz 2.1.2
auf nicht triviale Weise quasiergodisch; und
Fixpunkte oder Zyklen sind auf triviale Weise quasiergodisch. 
\newline 
Zwar wird alleine durch den 
Satz 2.1.2 die Ergodenfrage multidimensionaler kontinuierlicher Billards
gekl"art,
die vielfach untersuchten, 
sogenannten Billards\index{Billard} aber, bei denen St"osse 
und Reflexionen an den Banden diskontinuierlich 
modelliert sind,
entziehen sich aber der   
unmittelbaren Anwendung des  
Satzes von der Existenz der Zimmer, 
der keine Aussage "uber die 
Ergodenfrage von Zustandsr"aumen macht, die durch Trajektorien 
partioniert sind, die nicht stetig sind; die also gleichsam
auseinanderreissen.
Wir d"urfen hier
stolz bemerken, dass unsere mathematische Formulierung des Quasiergodenproblemes 
zwar einerseits mehr leistet, als den Bedarf der Physik 
kausal abgeschlossener und finiter Systeme
zu decken, falls deren jeweilige trajektorielle Partitionen von 
$\mathcal{C}$-Kontinuit"at sind; d.h., so, dass
alle Trajektorien stetig sind. Andererseits erfordert schon die 
Behandlung der vielfach untersuchten 
Billards eine Ausdehnung des 
Satzes von der Existenz der Zimmer 2.1.2, die "uber die 
$\mathcal{C}$-Kontinuit"at der 
jeweiligen trajektoriellen Partitionen hinausgeht. Auf diese
Ausdehnung des 
Satzes von der Existenz der Zimmer 2.1.2 auf
trajektorielle Partitionen mit nur 
st"uckweise stetigen Trajektorien gehen wir in diesem Traktat nicht ein,
in dem wir stattdessen anderen Gesichtspunkten nachgehen.
\newline
Warum? Warum schreiten wir "uberhaupt zuerst "uber den Satz von der
Existenz der Zimmer 2.1.2 hinaus zum elementaren Quasiergodensatz 1.1 fort, ehe wir 2.1.2 f"ur 
st"uckweise stetige trajektorielle Partitionen generalisieren?
Und zeigen stattdessen die Identit"at der Zimmer mit
den st"uckweise glatten
minimalen invarianten Mannigfaltigkeiten 
$\mu \in \mathcal{M}^{1/2} ([\Psi])$?
Wir wollen uns nicht vorhalten, dass wir zwar 
mit dem Satz von der Existenz der Zimmer
die Ergodenfrage
abgeschlossener Systeme beantworten, dabei die 
jeweiligen Zimmer aber insofern nicht kennen,
als sich uns kein Weg zeigt, sie explizit zu bestimmen!
Gem"ass (\ref{dmimmo}) ist das Mengensystem minimaler invarianter Mannigfaltigkeiten
\begin{equation}\label{dmimoo} 
\mathcal{M}^{1/2} ([\Psi])=
\Bigl\{\lambda^{ -1}(\omega): (\lambda,\omega)\in{\rm Inv}^{1/2}([\Psi])\times\mathbb{R}\Bigr\}[3]\ .
\end{equation}
Eingestanden, die explizite Bestimmung der minimalen invarianten Mannigfaltigkeiten gem"ass dieser Identifizierung
ist dabei, wie wir erl"auterten, zun"achst wenig verfahrensgebend, doch womit uns findige Numeriker -- oder findige Analytiker -- "uberraschen,
k"onnen wir nicht absehen. 
\newline 
Die Identit"at (\ref{dmimoo}) und ihr Beweis im Abschnitt 2.3 mittels kinezentrischer Felder hat aber 
hervorgehobenermassen vor dem Hintergrund der Identifizierung nicht trivialer Zimmer mit sensitiven
Attraktoren im Abschnitt 2.2
eine 
weitreichende Bedeutung, die "uber die Interpretation innerhalb der Physik hinaus reicht:\newline
Dass innerhalb 
jeder Modellierung
eines deterministisch evolvierenden Szenario einer empirischen Wissenschaft durch einen kompakten 
Zustandsraum mit glatten Trajektorien
die m"ogliche Kenntnis "uber ein kausal abgeschlossenes und finites System damit ausgesch"opft ist,
dass dessen Invariante und eine aktuelle Lokalisierung dessen Zustandes
bekannt ist, ist f"urwahr nichts Neues. Das ist das Thema der Hamiltonischen Theorie.
Sie lehrt uns, dass die Kenntnis aller Invarianter und die Assymptote beliebig 
scharfer Zustandslokalisierung der maximale Kenntnisstand ist, weil dieser Kenntnisstand
der totale ist. 
\newline 
Dass dabei die aktuelle Lokalisierung eines Zustandes in einem sensitiven Attraktor
im Hinblick auf die Bestimmung einzelner Trajektorien, die dessen jeweilige Zustandsentwicklungen darstellen,
g"anzlich wertlos ist, dass
die einzelne Trajektorie in einem sensitiven Attraktor ohne 
empirische Relevanz ist und die einzelne Trajektorie nur noch eine ideelle Hilfslinie zur Konstruktion des jeweiligen sensitiven Attraktors ist,
ist die Bote der Chaostheorie.\index{Chaos} 
Das ist f"urwahr auch nicht neu. \newline 
Der Satz von der Existenz der Zimmer sagt nun aber mehr, wenn er uns zeigt, dass
die aktuelle Lokalisierung 
des Zustandes einer Darstellung eines kausal abgeschlossenen und finiten Systemes
{\em im generischen Fall} die Lokalisierung
in einem sensitiven Attraktor ist. Und das ist sehr wohl neu.\newline 
Es wurde zwar bislang durchaus schon vorgefunden, dass wir \glqq in der Realit"at\grqq  auf Schritt und Tritt 
auf
\glqq nichtlineare Systeme\grqq stossen und es in deren jeweiligen Zustandsr"aumen
sensitive Attraktoren {\em geben kann} und in der Regel auch gibt, in denen die nicht noch "arger dramatisierbare Form \glqq chaotischen 
Verhaltens\grqq  vorliegt. Diese Aussage, dass \glqq Chaos\grqq  insofern \glqq ganz normal\grqq  sei,
ist aber eine h"ochst unscharfe, wengleich eine zutreffend negative Einsch"atzung der 
Beschreibbarkeit der Welt. Hier nun wird eine weitere Dramatisierung dieser Beschreibbarkeitsnegation vorgelegt:
Es gibt diese mit dem sensitiven Attraktor verbundene dramatischste Form \glqq chaotischen 
Verhaltens\grqq  nicht nur \glqq normalerweise in der Realit"at\grqq: Diese dramatischste Form
ist die einzige Alternative zur Statik des Fixpunktes oder zur Zyklik der geschlossenen Trajektorie --
falls ein kompakter Zustandsraumes vorliegt.
Ist letztere nicht gegeben, dann liegt der Zerfall des Systemes vor, wie wir hier nicht ausf"uhren werden.
Insofern ist durch den elementaren Quasiergodensatz 1.1 die durch die Chaostheorie negierte  
Beschreibbarkeit der Welt sowohl extremalisiert als auch universalisiert. 
\newline 
Jeder
jeweilige sensitive Attraktor ist dabei
gem"ass (\ref{dmimoo}) durch die Invarianten vollst"andig festgelegt. Auch das ist also neu, dass sich gem"ass der Erg"anzung des Satzes 2.1.2. zum elementaren Quasiergodensatz 1.1
damit die empirische Relevanz der Trajektorien hin zu der {\em Relevanz der Invarianten} des Systemes verschiebt.
Das, was die Relevanz der Invarianten nicht zur totaleinzigen, zu einzig deren Relevanz macht,
ist lediglich, dass sogar noch in einem sensitiven Attraktor eine jeweilige Lokalisierung eines Zustandes 
eine Aussage "uber die Transeabilit"at dieser Lokalisierung in eine andere erlaubt.
Diese Transeabilit"at ist zwar mit Mitteln der Wahrscheinlichkeitsrechnung formulierbar. 
Diese Transeabilit"at "uberschreitet aber den Radius der herk"ommlichen Masstheorie in provokanter, in zun"achst geradezu 
paradoxer Weise: Die "Uberg"ange von einer jeweiligen Lokalisierung eines Zustandes eines 
sensitiven Attraktors 
in eine 
andere Lokalisierung sind n"amlich durch Trajektorienabschnitte objektiviert, die alle das Lebesgue-Mass null haben.   
\newline
Und was ist schliesslich dies, wie aus dem Beweis der Supplementierung des 
Satzes 2.1.2. zum elementaren Quasiergodensatz 1.1 im Abschnitt 2.3 hervorgeht,
dass die Invarianten des Systemes zu den erhebbaren Mittelwerten 
"uber den jeweiligen Attraktor "aquivalent sind?
Diese erhebbaren Mittelwerte eines $\nu$-dimensionalen Attraktors sind insofern 
zu den Invarianten des Systemes mit $n$ Freiheitsgraden "aquivalent, als  
$n-\nu$ erhebbare unbh"angige Mittelwerte genauso wie die 
$n-\nu$ Invarianten des jeweiligen Attraktors denselben festlegen.
Was ist nun schliesslich dies?
Dies ist nicht weniger als die endg"ultig gesch"arfte Formulierung der Boltzmannschen Hellsicht
seiner tragischen Ergodenhypothese.
\newline 
Es sei $f$ ein st"uckweise diffeomorpher Hom"oomorphismus des $\mathbb{R}^{ n}$ in sich, dessen
Wertemenge $f(\mathbb{R}^{ n})= \mathbb{R}^{ n}$ ist.
F"ur $f$
und f"ur jede kompakte trajektorielle Partition $[\Psi]$ gilt offensichtlich, dass $f([\Psi])$ 
eine kompakte trajektorielle Partition ist und, dass 
alle diese st"uckweise diffeomorphen Hom"oomorphismen $f$ differenzierbare Mannigfaltigkeiten
in dem Sinn erhalten, dass die "Aquivalenz 
$$\mu\in \mathcal{M}^{1/2} ([\Psi])\  \Leftrightarrow\ f(\mu)\in \mathcal{M}^{1/2} (f([\Psi])) $$
gilt. Wobei ausserdem die "Aquivalenz
$$\tau\in\ [\Psi]\  \Leftrightarrow\ f(\tau)\in\ f([\Psi])$$
wahr ist und, dass f"ur alle $\tau\in\ T\ ([\Psi]),\ \mu\in\ \mathcal{M}^{1/2} ([\Psi])$ die "Aquivalenz
$$\mathbf{cl} (\tau) = \mu\  \Leftrightarrow\ \mathbf{cl} (f (\tau)) = f(\mu)$$
gilt:
Diese st"uckweise diffeomorphen Hom"oomorphismen $f$
erhalten invariante minimale Mannigfaltigkeiten, Trajektorien und deren Abschl"usse. 
Es gilt also:
\newline
\newline
{\em Wenn wir zeigen k"onnen, dass f"ur alle kompakten trajektoriellen $\mathcal{C}^{1}$-Partitionen $[\Psi]$ gilt}
\begin{equation}\label{eines}
\{\mathbf{cl} (\tau):\tau\in[\Psi]\}=\mathcal{M}^{1/2} ([\Psi])
\end{equation}
{\em und zweitens, dass f"ur sie f"ur alle $\tau\in [\Psi],\ \mu\in\ \mathcal{M}^{1/2} ([\Psi])$ die "Aquivalenz}
\begin{equation}\label{iddito}
\mathbf{cl} (\tau) = \mu\  \Leftrightarrow\ \tau\in\ \mu
\end{equation}
{\em gilt, zeigen wir, dass der elementare Quasiergodensatz wahr ist.}\index{Reduktion kompakter trajektorieller Partitionen}
\newline
\newline
Diese beiden reduzierten, zentralen Aussagen werden wir nun
in den Abschnitten 2.1 und 2.3 beweisen. 
Der Abschnitt 2.2 ist mehr als ein aufschlussreiches Interludium. In diesem Abschnitt 
zeigen wir die Sensitivit"at f"ur nicht triviale Zimmer, welche der Satz 2.2.2 behauptet.
Dessen Korollar, der Trivialit"atssatz\index{Trivialit\"atssatz} ist es, der den Zweck unsere Konzeption, mit Hilfe der Komanenz und Immanenz zu argumentieren,
deutlich macht. Der Trivialit"atssatz grenzt den Satz von der Existenz der Zimmer\index{Satz von der Existenz der Zimmer} 2.1.2 
ab
gegen die relativ triviale 
Partitivit"atsbehauptung, die f"ur insensitive Phasenfl"usse gilt. 
\newline
Das Fazit des Satzes 2.2.2 ist es, dass kontinuierliche dynamische Systeme entweder keinen kompakten Zustandsraum haben,
was dem Systemzerfall\footnote{Der Fall, dass Trajektorien in 
einen im Unendlichen liegenden Zustand laufen, jedoch aus demselben zur"uckkehren, ist transformatorisch reduzierbar auf die 
Beschreibung in einem entsprechend transformierten kompakten Zustandsraum. Nur, wenn die Trajektorien in 
einen im Unendlichen liegenden Zustand laufen und dabei nicht mehr in beschr"ankbare Bereiche zur"uckkehren,
ist diese transformative Reduktion auf den Fall nicht mehr m"oglich, in dem die Pr"amissen des Satzes von der Existenz der Zimmer
vorliegen. Dieser irreduzible Fall eines nicht kompakten Zustandsraumes beschreibt den Systemzerfall.} entspricht; oder aber -- grob gesagt -- sie sind in dem Sinn stabil,
dass deren Zustandsraum kompakt ist und dann liegt aber die Universalit"at
sensitiver Attraktoren insofern vor, als 
f"ur diese kontinuierlichen dynamischen Systeme die Alternativengegebenheit gilt,
die dem Hamletschen Fragen
\begin{center}
{\bf sensitiver Attraktor oder Zyklus oder Fixpunkt -- das ist hier die Frage!}
\end{center}
zugrundeliegt.\newline\newline
{\small Wer bis hierher mitkam, dem ist deutlich geworden, dass diese Abhandlung noch diese
Frage offen\index{offene Frage} l"asst:
N"amlich die Frage, ob es "uberhaupt nicht-triviale Zimmer f"ur Heine-Descartessche Kollektivierungen gibt.
Ja, diese Frage bleibt offen, was uns wenig verwundern sollte.\newline
Auch nach Abschluss dieser Abhandlung 
ist die Frage, ob es "uberhaupt nicht-triviale Zimmer gibt, nicht in dem Sinn reduzierbar,
dass aus der Insensitvit"at jeder Heine-Descartesschen Partition, die nur Zyklen umfasst,
die Existenz nicht-triviale Zimmer folgt. Die Insensitvit"at jeder Heine-Descartesschen 
Partition, die nur Zyklen umfasst, g"alte es dabei erst zu zeigen.
Gesetzt, dies gel"ange, so folgte damit aber die Existenz nicht-triviale Zimmer
auch noch nicht:
\newline
Wenn es eine Heine-Descartessche Kollektivierung gibt, die eine Partition einer kompakten Mannigfaltigkeit in lauter 
Zyklen oder Fixpunkte ist und
welche dabei
nicht insensitiv eichbar ist, (d.h., wenn jede Partition $[\Psi]$ einer kompakten Mannigfaltigkeit in lauter 
Zyklen oder Fixpunkte, die lokal parallel verlaufen,
so beschaffen ist, dass diese Partition $[\Psi]$ nicht insensitiv eichbar ist,) 
dann ist es m"oglich, dass es keine nicht-trivialen Zimmer f"ur Heine-Descartessche Kollektivierungen gibt;
f"ur Heine-Descartessche Kollektivierungen genausowenig, wie f"ur insensitive Kollektivierungen.
Wir sind dabei geneigt, zu sagen, 
dass aber der andere Fall, dass es keine solche Partition $[\Psi]$ gibt, offensichtlich zutreffe. 
Dieser andere Fall schliesst nicht aus,
dass keine nicht-trivialen Zimmer existieren.
Wenn jede solche Heine-Descartessche Partition $[\Psi]$
in Zyklen oder Fixpunkte
so beschaffen ist, dass es eine Flussfunktion $\hat{\Psi}$ gibt, f"ur die $[\Psi]=[\hat{\Psi}]$ ist,
wobei $\hat{\Psi}$ insensitiv ist und nur stetige Fl"usse hat,
dann ist n"amlich nur gezeigt, dass es nicht-triviale Zimmer f"ur Heine-Descartessche Kollektivierungen geben muss --
sofern es Heine-Descartessche Kollektivierungen kompakter Mannigfaltigkeiten gibt, die --
sensitiv sind... Ist das teuflisch? zyklisch? zirkul"os?\newline
Wir kommen um den konstruktiven Nachweis nicht-trivialer Zimmer letztlich wohl nicht herum. 
Dieser Nachweis ist
erbracht, wenn
wir f"ur ein einziges konkretes dynamisches System, dessen Trajektorien eine Heine-Descartessche Kollektivierung
bilden, zeigen, dass es nicht lauter zyklische Trajektorien hat.}

\chapter{Der Beweis des elementaren Quasiergodensatzes}
\begin{flushright}
Allen Kraftvollen: Der Stein rollt!\\
\end{flushright}
\section{Erster Schritt: Der reelle Determinismus impliziert, dass Zimmer existieren.}\label{secao}
Was sind Zimmer? -- Was Zimmer hier innnerhab dieser Abhandlung sind, das haben wir bereits im letzten Unterabschnitt gesagt: Zimmer sind die abgeschlossenen H"ullen von Trajektorien.
Sie sind aber nicht etwa die abgeschlossenen H"ullen irgendwelcher eindimensionaler Mannigfaltigkeiten.
Nein, sie sind die abgeschlossenen H"ullen 
$$\mathbf{cl}_{\mathbf{T}(\mathbf{P}_{2}h)} (\tau)\in [[h\circ\Psi]]_{\mathbf{T}(\mathbf{P}_{2}h)}$$
der hom"oomorphen Bilder $h(\tau)$ der Trajektorien $\tau\in [\Psi]$ einer
trajektoriellen Partition $[\Psi]$ f"ur einen Hom"oomorphismus $h$ bez"uglich
der Topologien $\mathbf{T}(\mathbf{P}_{1}h)$ und $\mathbf{T}(\mathbf{P}_{2}h)$. 
Wir sehen, dass
es Zimmer auch in viel allgemeineren R"aumen als in endlichdimensionalen und reellen R"aumen gibt. Ob sich
die Benennung der Zimmer\index{Zimmer} als
solche verfestigen wird, oder aber, ob letztere bloss eine Arbeitsbenennung innerhalb dieses Traktates
bleibt, das sie
mit sensitiven Attraktoren
identifizieren wird,
sei dahingestellt.
\footnote{Nicht nur die Theorie des Determinismus metrischer R"aume ist dem Autor bereits wohlvertraut.
Letztere sprengt aber den hier gesteckten Rahmen. In der allgemeinen Theorie tut sich die Differenz zwischen generalisierten Zimmern und generalisierten Attraktoren auf und der Begriff der Zimmer und der Begriff der Attraktoren verzweigt sich,
was die eigene Benennung der in diesem Sinn autonomen generalisierten Zimmer rechtfertigt.
Begriffsverzweigung dieser Art gibt es bei so mancher 
Generalisierung, wie alte wie junge alte Hasen wissen: Wenn eine Generalisierung den Bereich verl"asst, innerhalb welchem eine Symmetrie
gegeben ist, 
dann dissoziieren auch die begrifflichen Koinzdenzen, welche 
diese Symmetrie voraussetzen.
Die pr"agenerelle Symmetrie kann dabei eine im Stillschweigen vorausgesetzte sein,
wenn dieselbe erst ausserhalb jenes Bereiches, erst im allgemeineren Rahmen 
als eine relative 
formulierbar ist, die nur im Bezug auf die allgemeinen Konstrukte
der verallgemeinerten Theorie konstituierbar ist. 
Dann geht die Generalisierung mit der Entdeckung der Symmetrie des
jeweiligen pr"agenerellen Generalisierungskeimes einher.
Die eine Seite der Medaille ist die jeweilige verallgemeinerungserzeugte Begriffsdissoziation, deren andere ist die
Entdeckung der jeweiligen
verallgemeinerungserzeugten Symmetrie des Pr"agenerellen.}
\newline\newline
Wir haben jene zwei ganz am 
Schluss des ersten Abschnittes exponierten Aufgaben gestellt, diejenigen Behauptungen zu 
beweisen, die (\ref{iddito}) und (\ref{eines}) darstellen. 
Die Exponierung der "Aquivalenz (\ref{iddito}) und der Identit"at (\ref{eines}) ist das Ergebnis der Konzentrierung der Quasiergodenfrage auf 
die entprechenden
beiden Beweisaufgaben.
Die Zentralit"at der "Aquivalenz (\ref{iddito}) und der Identit"at (\ref{eines}) besteht gerade darin,
dass die Konzentrierbarkeit der Quasiergodenfrage auf (\ref{iddito}) und (\ref{eines})
gegeben ist.\newline
Die
Identit"at (\ref{eines}) erweitert die "Aquivalenz (\ref{iddito}) zu der Aussage, dass f"ur jede trajektorielle Partition $[\Psi]$ 
und f"ur jede deren minimaler invarianter Mannigfaltigkeiten
$\mu \in \mathcal{M}^{ 1/2} ([\Psi])$ und f"ur jede der Trajektorien $\tau\in [\Psi]$ der trajektoriellen Partition $[\Psi]$ 
sogar die Identit"at
\begin{equation}\label{eddito}
\mathbf{cl} (\tau) = (\mathcal{M}^{1/2}  ([\Psi]))_{\tau}
\end{equation}
wahr ist.
Ohne, dass wir in diesem enggesteckten Rahmen die logische Architektur der Theorie dynamischer Systeme im Hinblick auf die Quasiergodenproblematik 
darlegen wollen, um innerhalb jener 
herausgearbeiteten logischen Architektur
die Position der Aussagen (\ref{eines}) und (\ref{iddito}) zu bestimmen:
Wir widmeten den ersten Abschnitt dieser Zielsetzung,
darzustellen, 
dass die "Aquivalenz (\ref{iddito}) die Quasiergodenfrage schon im Wesentlichen beantwortet.
Wenn wir jene beiden Beweisaufgaben erf"ullen, (\ref{eines}) und (\ref{iddito}) zu zeigen, beschreiben wir die Geometrie der trajektoriellen Partitionen. Geometrie nennen
wir vereinfachterweise die Behandlung von Punktmengen in endlichdimensionalen und reellen R"aumen.\newline
Es gibt eine Verallgemeinerung von trajektoriellen Partitionen f"ur metrische R"aume, in denen generalisierte
Trajektorien existieren. 
Diese Verallgemeinerungen trajektorieller Partitionen k"onnen je nach der Reichhaltigkeit der ihnen zugrundeliegenden metrischen R"aume 
allerdings auch recht triviale Konstruktionen sein. Heine-Descartessche trajektorielle Partitionen sind komanent\index{komanente trajektorielle Partition} und global immanent:\index{global immanente trajektorielle Partition}
\newline
Wir finden n"amlich folgende beiden  
auf der Grundlage der Metrisiertheit des Zustandsraumes 
verallgemeinerbaren Eigenschaften der
trajektoriellen Partitionen, mit denen es folgende Bewandnis hat:
Diese beiden Eigenschaften von trajektoriellen Partitionen sind die Grundlage einer abstrakteren und
axiomatischen Fassung der Ergodentheorie, die sich weit 
"uber die Ergodentheorie endlichdimensionaler und reeller Zustandsr"aume hinaus fortsetzen l"asst, n"amlich
sogar "uber metrische Zustandsr"aume hinaus.
Und schon alleine deshalb sollen wir jene beiden Eigenschaften herausstellen, auch wenn wir hier auf die
generalisierte Ergodentheorie nicht eingehen.
\newline
Wir m"ussen aber geradezu jene besagten beiden  
verallgemeinerbaren Eigenschaften aus folgenden Grund pr"asentieren: 
Ihre kombinierte Betrachtung leitet n"amlich zur zu beweisenden Einsicht,
dass die abgeschlossenen H"ullen der Trajektorien 
einer jeweiligen trajektoriellen Partition
deren Zustandsraum partionieren. 
\newline
\newline
{\em Wir sind hier bem"uht, nicht alleine f"ur 
"uberragende Experten zu schreiben. Deshalb widmen wir mancher Betonung Raum:
Und vielleicht ist selbst der Fachmann nicht dagegen gefeit, beim fl"uchtigen Blick
die zu zeigende Partitivit"at 
$$[[\Psi]]\in\mathbf{part}(\mathbf{P}_{2}\Psi)$$
f"ur {\bf Heine-Descartessche Kollektivierungen} $[\Psi]$ gem"ass der 
Punkte} (1)-(3) {\em der Einleitung und die 
gleiche Partitivit"at
f"ur {\bf insensitive Kollektivierungen}\index{insensitive Kollektivierung} $[\Psi]$ gem"ass}
(\ref{knurrze}) {\em 
auseinanderzuhalten. Dann muss jener unsere Argumentation als umst"andlich 
und als befremdlich empfinden.\newline
F"ur insensitiv eichbare Flussfunktionen\index{insensitiv eichbare Flussfunktion} ist der Beweis
jener Partitivit"at einfach und bekannt. Dem weniger Versierten darf die
Aufgabe zugemutet werden, diesen einfach Beweis zu formulieren: Dazu
braucht es  
die Konzeption von Komanenz und Immanenz freilich nicht.
Die Aufgabe indess, jene  
Partitivit"at f"ur Heine-Descartessche Kollektivierungen zu zeigen, wird
vermutlich unter der mit der Konzeption von Komanenz und Immanenz verbundenen Perspektive
einer eleganten Argumentation erschlossen.} 
\newline
\newline
Die Eigenschaft jeder gem"ass (\ref{kurza})-(\ref{kurzc}) beschaffenen
Flussfunktion $\Psi$ binnen einer als Zustandsraum aufgefassten 
Teilmenge $\zeta\subset \mathbb{R}^{n}$, dass es
f"ur alle Paare $(t,\varepsilon)\in \mathbb{R}\times\mathbb{R}^{+}$ eine positive reelle Zahl $\delta$ von der Art gibt, dass f"ur alle 
Paare von Zust"anden $(x,y)\in\zeta$ die Implikation 
\begin{equation}
||x  - y|| <  \delta\ \Rightarrow\ ||\Psi(x,t) - \Psi(y,t)|| < \varepsilon
\end{equation}
wahr ist,
nennen wir genau dann, wenn 
die G"ultigkeit dieser Implikation
gegeben ist, die Komanenz 
der Flussfunktion $\Psi$.  
Die Komanenz ist demnach eine spezielle Form der Stetigkeit.\footnote{Und zwar eine,
die wir f"ur alle diejenigen Funktionen $\Phi$
auf naheliegende Weise verallgemeinern k"onnen,
f"ur die es Menge ${\rm A}$ gibt, 
f"ur die $\Phi\in{\rm A}^{{\rm A}\times\mathbb{R}}$ gilt, wobei ${\rm A}$
die durch eine auf ${\rm A}$ 
erkl"arte Metrik $d_{{\rm A}}$ mit der durch diese Metrik induzierten Topologie $\mathbf{T}(d_{{\rm A}})$ topologisiert
ist. Wir k"onnen diese Stetigkeitsform aber auch f"ur alle Funktionen $\Phi\in{\rm A}^{{\rm A}\times\mathbb{R}}$
generalisieren, f"ur die es Menge ${\rm A}$ gibt, die durch eine nicht notwendigerweise metrische Topologie 
$\mathbf{T}({\rm A})$ topologisiert ist. Sowohl dieser naheliegenden metrischen Formulierung
der Komanenz als auch deren weiteren Generalisierung gehen wir hier nicht nach, weil diese jenseits der 
programmatischen 
Elementarit"at dieses Traktates liegt.}
Dass wir diese als Komanenz bezeichnen, liegt daran, dass die Trajektorien
der Kollektivierung $[\Psi]$ einer komanenten Flussfunktion $\Psi$ in dem Sinn
beieinander bleiben, ko-manieren, den diese Definition der 
Komanenz pr"azisiert.
Ferner ordnet diese Definition 
die Komanenz zwar zun"achst einer jeweiligen Flussfunktion 
und nicht 
etwa der
Kollektivierung $[\Psi]$ der jeweiligen Flussfunktion zu.\index{Komanenz eines Trajektorienkollektives}
Wir k"onnen aber auch jeder 
Kollektivierung $[\Psi]$ einer jeweiligen Flussfunktion $\Psi$ binnen $\zeta\subset \mathbb{R}^{n}$
gegebenenfalls
deren Komanenz 
oder aber deren Nicht-Komanenz
als Eigenschaft zuordnen, die
die Kollektivierung $[\Psi]$ entweder hat oder aber nicht. Und in diesem Gebrauch ist
die Komanenz die Eigenschaft eines Mengensystemes von Trajektorien und
nicht die Eigenschaft einer Abbildung, sodass sie dann nur indirekt als eine Form der Stetigkeit aufgefasst werden kann.
Das, der
Kollektivierung $[\Psi]$ die Komanenz zuzuordnen,  
ist deshalb m"oglich, weil f"ur alle Flussfunktionen $\Xi_{1}$ und $\Xi_{2}$, f"ur die
$$[\Xi_{1}]=[\Xi_{2}]=[\Psi]$$
gilt,
die "Aquivalenz
\begin{equation}\label{ovaeva}
\begin{array}{c}
\Bigl(\forall\ (t,\varepsilon) \in  \mathbb{R}\times\mathbb{R}^{ +}\ 
\exists\delta \in \mathbb{R}^{ +}\ \forall\ (x, y)\in\bigcup [\Psi]\\
||x  - y|| <  \delta\ \Rightarrow\ ||\Xi_{1}(x,t) - \Xi_{1}(y,t)|| < \varepsilon\Bigr)\\
\Leftrightarrow \\
\Bigl(\forall\ (t,\varepsilon) \in  \mathbb{R}\times\mathbb{R}^{ +}\ 
\exists\delta \in \mathbb{R}^{ +}\ \forall\ (x, y)\in\bigcup [\Psi]\\
||x  - y|| <  \delta\ \Rightarrow\ ||\Xi_{2}(x,t) - \Xi_{2}(y,t)|| < \varepsilon\Bigr) 
\end{array}
\end{equation}
wahr ist:\index{Komanenz eines Trajektorienkollektives}\index{Komanenz}
Eine Flussfunktion ist genau dann komanent, wenn ihre Kollektivierung komanent ist und
ein Kollektiv ist genau dann komanent, wenn es eine komanente
Flussfunktion gibt, deren Kollektivierung jenes Kollektiv ist.
\newline
Der Begriff der Komanenz ist f"ur jedwede reelle Flussfunktion
entworfen. Den Begriff sowohl der punktweisen als auch der globalen Immanenz hingegen entwerfen 
wir nur f"ur spezielle Flussfunktionen, n"amlich f"ur normale Flussfunktionen.\index{normale Flussfunktion}
Jede Flussfunktion $\Xi$ binnen einer 
durch eine jeweilige Topologie $\mathbf{T}$
topologisierten Menge ${\rm X}$ bezeichnen wir als
bez"uglich der Topologie $\mathbf{T}$ genau dann normal, wenn diese Flussfunktion
$\Xi$ in ihrer zweiten Ver"anderlichen stetig ist bez"uglich der 
nat"urlichen Topologie des Zahlenstrahles $\mathbf{T}(1)$ und bez"uglich jener
Topologie $\mathbf{T}$, die den jeweiligen Zustandsraum ${\rm X}$ topologisiert.
Das heisst, dass genau dann, wenn $\Xi$ bez"uglich der Topologie $\mathbf{T}$ normal ist,
f"ur alle $x\in {\rm X}$
\begin{equation}\label{inormkl}
\Xi(x,\mbox{id})\in\mathcal{C}_{\mathbf{T}(1),\mathbf{T}}(\mathbb{R},{\rm X})
\end{equation}
gilt.\index{Normalit\"at einer Flussfunktion bez\"uglich einer Topologie} 
Wir k"onnen den Begriff insensitiver Eichbarkeit f"ur jede Zustandsraumtopologisierung\index{insensitiv eichbare Flussfunktion} 
analog zu (\ref{knurrze}) verfassen.
Die Normalit"at einer Flussfunktion impliziert nicht, dass sie insensitiv eichbar ist.
Dementsprechend meinen wir genau dann, wenn wir eine Flussfunktion $\Psi$ binnen 
einer Teilmenge
$\zeta\subset \mathbb{R}^{n}$
eines reellen Raumes der Dimension $n\in\mathbb{N}$ unattributiert als normal bezeichnen,
dass $\Psi$ bez"uglich der Topologie $\mathbf{T}(n)$ normal ist, wobei 
$\mathbf{T}(n)$ die nat"urliche Topologie des $\mathbb{R}^{n}$ ist.
\newline
Als die
punktweise Immanenz\index{punktweise Immanenz einer Flussfunktion} einer normalen, d.h., in ihrer
zweiten Ver"anderlichen stetigen, 
Flussfunktion $\Psi$ binnen $\zeta\subset \mathbb{R}^{n}$
bezeichnen wir 
gegebenenfalls
deren Eigenschaft,
dass 
f"ur jede Zahl $\varepsilon \in \mathbb{R}^{+}$ und 
einen Zustand
$x\in\zeta$ eine Zahl $T\in\mathbb{R}^{+}$ von der Art existiert,
dass die Menge
$$\mathbb{B}_{ \varepsilon} (0) + \Psi(x,[-T, T])$$
die gesamte Trajektorie $\Psi(x, \mathbb{R})$ enth"alt.
Und zwar sei dies die punktweise Immanenz der Flussfunktion $\Psi$ im 
Zustand $x$.\newline
Die Menge $\mathbb{B}_{ \varepsilon} (0) + \Psi(x,[-T, T])$
k"onnen wir uns als einen Schlauch 
vom Durchmesser $2\varepsilon$ 
vorstellen, der von zwei Kugelkappen ebendesselben Durchmessers verschlossen ist, 
und dessen Achse hierbei der Trajektorienabschnitt $\Psi(x,[-T, T])$ ist -- falls $\varepsilon$ oder
$T$ hinreichend klein sind.
Was wir uns als den imaginativen Sachverhalt vorstellen,
dass dieser Schlauch um den Trajektorienabschnitt 
$\Psi(x,[-T, T])\subset \Psi(x,\mathbb{R})$ 
die gesamte Trajektorie $\Psi(x, (\mathbb{R})$
enth"alt, ist in der Inklusion
\begin{equation}\label{inkl}
\mathbb{B}_{ \varepsilon} (0) + \Psi(x,[-T, T])\ \supset\ \Psi(x,\mathbb{R})
\end{equation}
formuliert. 
Als
die globale Immanenz einer normalen Flussfunktion $\Psi$ binnen $\zeta\subset \mathbb{R}^{n}$
\index{globale Immanenz einer Flussfunktion} bezeichnen wir gegebenenfalls genau diesen Sachverhalt, dass die folgende Aussage
f"ur $\Psi$ gilt: F"ur jede positive reelle Zahl $\varepsilon \in \mathbb{R}^{ +}$ 
existiert eine andere positive reelle Zahl $T \in \mathbb{R}^{ +}$ von der Art,
dass f"ur alle Zust"ande $ x \in \bigcup [\Psi]$ die Inklusion 
(\ref{inkl}) gilt.
\newline
\newline
Dass die punktweise Immanenz 
einer jeden 
normalen Flussfunktion $\Psi$ binnen $\zeta\subset \mathbb{R}^{n}$
gegeben ist, 
deren Zustandsraum $\zeta$ beschr"ankt ist,
sehen wir leicht:
\newline
\newline
Weil $\zeta$ beschr"ankt ist, ist 
die abgeschlossene H"ulle der Trajektorie
$\mathbf{cl} (\Psi(x, \mathbb{R}))$ kompakt f"ur jeden Zustand
$x\in \zeta$.
F"ur jede positive reelle Zahl $\varepsilon$
gibt es eine Folge $\{t_{ j}\}_{ j \in\mathbb{N}} \in \mathbb{R}^{ \mathbb{N}}$ von der Art,
dass die Vereinigung 
\begin{equation}
\bigcup\Bigl\{ \mathbb{B}_{ \varepsilon} (\Psi(x,t_{ j}))\ :\ j \in \mathbb{N}\Bigr\}
\end{equation}
dieses Kompaktum $\mathbf{cl} (\Psi(x, \mathbb{R}))$ "uberdeckt. Der Satz von Heine-Borel\index{Satz von Heine-Borel}
sagt uns daher, dass eine
endliche Teilmenge $I \subset \mathbb{N}$ existiert, f"ur die die Inklusion 
\begin{equation}
\bigcup\Bigl\{ \mathbb{B}_{ \varepsilon} (\Psi(x,t_{ j}))\ :\ j \in I\Bigr\}  \supset\ 
\mathbf{cl}(\Psi(x, \mathbb{R}))
\end{equation}
gilt. Und daher existiert
auch ein Intervall $[-t (\varepsilon),  t (\varepsilon) ]$ von der Art,
dass f"ur alle $j\in I$ die eindimensionale Lokalisierung
\begin{equation}
t_{ j}\ \in\ [-t (\varepsilon), t  (\varepsilon)] 
\end{equation}
wahr ist. 
Wir k"onnen die 
hier via einer Flussfunktion $\Psi$ definierten  
Trajektorien $\Psi(x, \mathbb{R})$ aus dem
Zusammenhang mit einer Flussfunktion $\Psi$ befreien:
Jede 
stetige Abbildung $\vartheta$ des Zahlenstrahles in eine Menge $\zeta\subset\mathbb{R}^{n}$
des $\mathbb{R}^{n}$ mit endlicher Dimension $n\in\mathbb{N}$ hat offenbar
eine Wertemenge $\mathbf{P}_{2}\vartheta=\vartheta(\mathbb{R})$, die insofern eine Verallgemeinerung einer 
Trajektorie ist: 
Jede Trajektorie
ist zwar ein derartiges stetiges Zahlenstrahlbild $\vartheta(\mathbb{R})$. Es ist
aber nicht jedes solche Zahlenstrahlbild $\vartheta(\mathbb{R})$
eine Trajektorie. Denn das stetige Zahlenstrahlbild $\vartheta(\mathbb{R})$
kann 
im Gegensatz zu einer Trajektorie
so beschaffen sein, dass es sich selbst durchdringt.
Wegen der Stetigkeit der Abbildung $\vartheta$ gibt es nur endlich viele
Punkte der Wertemenge $\vartheta(\mathbb{R})$, in denen sich das Zahlenstrahlbild $\vartheta(\mathbb{R})$ selbst durchdringt,
falls dasselbe beschr"ankt ist. Daher ist es dann als eine Projektion
einer Trajektorie eines h"oher- jedoch endlichdimensionalen
Raumes $\mathbb{R}^{m}$ auffassbar. Vor diesem Hintergrund nennen wir
alle der beschriebenen stetigen Zahlenstrahlbilder $\vartheta(\mathbb{R})$ -- auch die nicht beschr"ankten --
projektive Trajektorien.
Jede beschr"ankte projektive Trajektorie\index{projektive Trajektorie} $\vartheta(\mathbb{R})\subset\mathbb{R}^{n}$ ist
offenbar immanent in diesem Sinn,
dass f"ur alle positiven und reellen Zahlen $\varepsilon$ und f"ur alle ihre Punkte $x\in \vartheta(\mathbb{R})$ eine
Art Schlauch $\mathbb{B}_{ \varepsilon} (0) + \vartheta(\mbox{I})$ um einen Abschnitt
$\vartheta(\mbox{I})\ni x$ existiert, welcher
die gesamte projektive Trajektorie $\vartheta(\mathbb{R})$ enth"alt, wobei der 
Trajektorienabschnitt $\vartheta(\mbox{I})$ endliche L"ange hat, weil $\mbox{I}\subset\mathbb{R}$ ein abgeschlossenes 
Intervall sein soll.
Wir erkennen die f"ur alle projektiven Trajektorien $\vartheta(\mathbb{R})$ g"ultige
"Aquivalenz von deren jeweiliger Immanenz im beschriebenen Sinn und deren jeweiliger Beschr"anktheit. 
\newline
Die Beschr"anktheit des Zustandsraumes $\zeta$
jeder normalen Flussfunktion $\Psi$ binnen $\zeta$ ist also
ein hinreichendes, aber offenbar kein notwendiges Kriterium daf"ur, dass
die stetige Flussfunktion $\Psi$ 
-- oder gleichermassen ihre Kollektivierung $[\Psi]$ --  
global immanent ist: Denn die globale Immanenz einer Flussfunktion $\Psi$
ist die punktweise Immanenz aller Trajektorien ihrer Kollektivierung $[\Psi]$,
wie die vorangegangene Betrachtung beschr"ankter 
projektiver Trajektorie deutlich macht,
sodass 
folgendes Pendant
zu der "Aquivalenz (\ref{ovaeva}) formulierbar ist:
F"ur alle normalen Flussfunktionen $\Xi_{1}$ und $\Xi_{2}$, f"ur die
$$[\Xi_{1}]=[\Xi_{2}]=[\Psi]$$
gilt, gilt
die "Aquivalenz\index{Immanenz eines Trajektorienkollektives}\index{Immanenz}
\begin{equation}\label{evaeva}
\begin{array}{c}
\Bigl(\forall\ (x,\varepsilon) \in \bigcup[\Psi] \times\mathbb{R}^{ +}\ 
\exists T \in \mathbb{R}^{ +}\\
\mathbb{B}_{ \varepsilon} (0) + \Xi_{1}(x,[-T, T])\ \supset\ \Xi_{1}(x,\mathbb{R})\Bigr)\\
\Leftrightarrow \\
\Bigl(\forall\ (x,\varepsilon) \in \bigcup[\Psi] \times\mathbb{R}^{ +}\ 
\exists T \in \mathbb{R}^{ +}\\
\mathbb{B}_{ \varepsilon} (0) + \Xi_{2}(x,[-T, T])\ \supset\ \Xi_{2}(x,\mathbb{R})\Bigr)\ .
\end{array}
\end{equation}
Eine normale Flussfunktion ist genau dann global immanent, wenn ihre Kollektivierung global immanent ist und
ein Kollektiv ist genau dann global immanent, wenn es eine global immanente
Flussfunktion gibt, deren Kollektivierung es ist.
\newline
Ist die Flussfunktion $\Psi$ global immanent, so existiert eine Abbildung 
\begin{equation}
\begin{array}{c}
a_{\Psi}:\ \zeta\times \mathbb{R}^{ +}\ \to \ \mathbb{R}_{ 0}^{ +}\ ,\\
(x, \varepsilon)\ \mapsto\ a_{\Psi} (x, \varepsilon)\ ,
\end{array}
\end{equation}
wobei wir die Zahl
\begin{equation}
 a_{\Psi} (x, \varepsilon):=\inf\Bigl\{t \in \mathbb{R}_{ 0}^{ +}:\ \Psi(x,[-t, t ]) + 
\mathbb{B}_{ \varepsilon} (0) \subset\ \Psi(x, \mathbb{R})\Bigr\}
\end{equation}
naheliegender Weise als die Immanenzzeit f"ur\index{Immanenzzeit} jedes jeweilige Paar $(x,\varepsilon)\in\zeta\times \mathbb{R}^{ +}$ 
bezeichnen.
Wenn die Abbildung 
$a_{\Psi} (\mbox{id}, \varepsilon)$ stetig ist, so existiert, wenn der Zustandsraum $\zeta$ zudem kompakt ist, sogar
die Zahl
\begin{equation}
A_{\Psi} (\varepsilon)\ :=\ \max\Bigl\{a_{\Psi} (x, \varepsilon)\ : x \in \zeta\ \Bigr\}\ ,
\end{equation}
eine $\varepsilon$-Immanenzzeit.
Und wenn $a_{\Psi} (\mbox{id}, \varepsilon)$ f"ur alle  $\varepsilon \in \mathbb{R}^{ +}$ stetig ist, so existiert demnach sogar 
eine pauschale Immanenzfunktion\index{pauschale Immanenzfunktion}
\begin{equation}\begin{array}{c}
A_{ \Psi} : \mathbb{R}^{ +} \to \mathbb{R}^{ +}\ ,\\ 
\varepsilon \mapsto A_{ \Psi} (\varepsilon)\ .
\end{array}
\end{equation}
Es gelingt uns, zu zeigen,
dass f"ur eine $\mathcal{C}^{1}$-Flussfunktion $\Psi$,
deren Zustandsraum $\zeta$ beschr"ankt ist,
die Funktion $a_{\Psi}\ (\mbox{id}, \varepsilon)$ stetig ist f"ur alle 
$\varepsilon \in \mathbb{R}^{ +}$,
wenn wir von der Komanenz der Flussfunktion $\Psi$ ausgehen k"onnen, deren 
Ableitung nach der
zweiten Ver"anderlichen $\partial_{2}\Psi$ stetig ist. F"ur sie
gibt es die Abbildung
\begin{equation} 
\begin{array}{c}
\hat{\Psi}:\zeta\ \to\ \mathbb{R}^{ n}\ ,\\ 
x\mapsto \hat{\Psi} (x):=\partial_{2}\Psi(x,0)\ ,
\end{array}
\end{equation}
sodass auch die Abbildung
\begin{equation} 
\begin{array}{c}
\zeta\times\zeta\ \to\ \mathbb{R}^{ +}_{ 0}\ ,\\ (x, y) \mapsto || \hat{\Psi} (x) - \hat{\Psi} (y) ||
\end{array}
\end{equation}
stetig ist und deren Maximum
\begin{equation}
\Theta (\varepsilon)\ :=\ \max \Bigl\{||\hat{\Psi} (x) - \hat{\Psi} (y)||\ :\ x, y \in \zeta\ \wedge\ ||x-y||\ \leq\ \varepsilon\Bigr\}
\end{equation}
f"ur alle $\varepsilon \in\mathbb{R}^{ +}$ existiert. Daher gibt
es 
eine entsprechende Abbildung $\Theta$ mit der Definitionsmenge $\mathbf{P}_{1}\Theta=\mathbb{R}^{ +}$, deren 
Wertemenge im positiven Zahlenstrahl
$\mathbb{R}^{ +}$ liegt, und deren
jeweilige Werte $\Theta (\varepsilon)$ sind. Diese Funktion
$\Theta$ ist offensichtlich stetig.
Und offenbar w"achst sie monoton, wobei $\Theta (\varepsilon)$ f"ur $\varepsilon\to 0$ verschwindet. $\Theta$ bleibt aber f"ur 
$$\varepsilon\ \geq \sup\Bigl\{||x-y||:\ x,\ y\ \in\ \zeta\Bigr\}$$ konstant,
weshalb $\Theta$ von einer streng monotonen differenzierbaren Funktion auf dem Intervall 
$$[0,\ \sup \{||x-y||: x, y\ \in\ \zeta\}]$$ von oben beschr"ankt wird, deren Ableitung 
auf ihrer Wertemenge den Maximalwert annimmt.
Also  gibt es einen Flussexponenten\index{Flussexponent}  $\kappa (\Psi) \in \mathbb{R}^{ +}$ als die kleinste
positive reelle Zahl, f"ur die f"ur alle $\varepsilon \in \mathbb{R}^{ +}$
\begin{equation}
\Theta (\varepsilon)\ \leq\ \kappa (\Psi)\ \varepsilon 
\end{equation}
gilt.
Also ist f"ur alle $x, y \in \zeta\ ,t \in \mathbb{R}^{ +}$ die Absch"atzung
$$||\Psi(x,t) - \Psi(y,t)||$$
$$\leq\ ||x - y||\ +\ \int_{0}^{t}\ ||\hat{\Psi}(x,s) - \hat{\Psi}(y,s)||\ ds$$
$$\leq\ ||x - y||\ +\ \kappa (\Psi)\ \int_{0}^{t}\ ||\Psi(x,s) -\Psi(y,s)||\ ds$$
wahr. 
Die folgende Annahme beinhaltet offenbar einen Widerspruch:
Die Annahme,
dass es f"ur eine stetige Funktion $\alpha$ mit der Definitionsmenge $\mathbf{P}_{1}\alpha=\mathbb{R}^{ +}$, deren 
Wertemenge im positiven Zahlenstrahl
$\mathbb{R}^{ +}$ liegt, eine positive reelle Zahl $T$ gibt, f"ur die f"ur alle $t \in\ [0,\ T]$ die Gleichung
$$\alpha (t) = \alpha (0)\ +\ \kappa (\Psi)\ \int_{0}^{t}\ \alpha (s)\ ds$$
gilt und dass dabei einerseits die Beschr"ankung
$$\alpha (0)\ \geq\ ||x\ -\ y||$$
eingehalten ist,
aber andererseits eine Zahl $t^{ \star}$ in dem Intervall $[0, T]$ existiert, f"ur die die Beschr"ankung
$$\alpha (t^{ \star})\ \leq\ ||\Psi(x,t^{ \star}) - \Psi(y,t^{ \star})||$$
zutrifft. 
Also ist f"ur alle $t \in \mathbb{R}^{ +}, x, y \in \zeta$ die Ungleichung
\begin{equation}\label{vovo}
||\Psi(x,t) - \Psi(y,t)||\ \leq\ ||x - y||\ \exp(\kappa(\Psi)t)
\end{equation}
erf"ullt.
\newline
Es sei $f$ ein Hom"oomorphismus zwischen endlichdimensionalen reellen R"aumen.
Wenn dieser Hom"oomorphismus auf dem Zustandsraum $\bigcup \Gamma$ einer trajektoriellen $\mathcal{C}^{1}$-Partition $\Gamma$ eines Kompaktums 
eines endlichdimensionalen reellen Raumes
Lipschitz-stetig ist,
wenn also eine Schranke $K \in \mathbb{R}^{ +}$ existiert mit
\begin{equation}
K\ :=\ \sup \left\{\frac{||f\ (x)-\ f\ (y)||}{||x\ -\ y||}\ :\ x, y\ \in\ \bigcup \Gamma\right\}\ ,
\end{equation}
so ist das Bild $f (\Gamma)=\{f(\tau):\tau\in\Gamma\}$ 
ein Kollektiv, das
komanent ist.
Wenn $\Psi$ eine stetige Flussfunktion ist, deren Kollektivierung $[\Psi]$ mit der trajektoriellen Partition
$\Gamma$ "ubereinstimmt und wenn
$\Psi$ dabei
den
Flussexponenten $\kappa(\Psi)\in\mathbb{R}^{+}$ hat, so gibt
es zu dem Bildkollektiv $f(\Gamma)$
also die Flussfunktion $f\circ\Psi$, deren Kollektivierung 
$$[f\circ\Psi]=f([\Psi])=f (\Gamma)$$
dieses Bildkollektiv $f(\Gamma)$ 
ist. 
Die Flussfunktion $f\circ\Psi$ 
hat dabei den Flussexponenten\index{Flussexponent}
\begin{equation} 
\log(K)\ +\ \kappa(\Psi) = \kappa(f\circ\Psi)\ .
\end{equation}
Zeigen wir nun die\newline\newline
{\bf Bemerkung 2.1.1: Die Unvereinbarkeit der Komanenz mit\newline 
der Existenz unstetiger Immanenzfunktionen}\index{Unvereinbarkeit unstetiger Immanenzfunktionen mit der Komanenz}\newline
{\em Die Immanenzfunktion $a_{\Psi}$ jeder komanenten und punktweise 
immanenten Flussfunktion $\Psi$ 
ist stetig.}
\newline
\newline
{\bf Beweis:}\newline
Es sei $$\{z_{ j}\}_{ j \in \mathbb{N}}\ \in\ \zeta^{ \mathbb{N}}$$ eine gegen $z \in \zeta$ 
konvergente Punktfolge.
Existierte eine positive reelle Zahl $\tau \in \mathbb{R}^{ +}$ von der Art,
dass f"ur eine positive reelle Zahl $\delta \in \mathbb{R}^{ +}$ die Ungleichung
$$a_{\Psi} (z_{j}, \delta)\ \leq\ \ a_{\Psi} (z,\delta)- \tau$$
f"ur alle $j \in \mathbb{N}$ g"alte, so g"abe es zum einen eine reelle Zahl
$$t(z) \in \mathbb{R} \setminus \bigcup\{a_{\Psi} (z_{ j}, \delta) \cdot [-1,1]:j \in \mathbb{N}\}$$
und zum anderen eine positive reelle Zahl $\eta$, f"ur die f"ur alle $j\in \mathbb{N}$
$$\Psi(z,t(z))\ \not\in\ 
\Psi\Bigl(z,\bigcup_{j\in\mathbb{N}} a_{\Psi} (z_{j},\delta) \cdot [-1, 1]\Bigr)+ \mathbb{B}_{\delta+3\eta} (0)$$
g"alte. 
Wegen der Komanenz der Flussfunktion $\Psi$ g"abe es dabei eine nat"urliche Zahl $j^{\star}\in\mathbb{N}$ von der Art,
dass auch die Ungleichung
\begin{equation}\label{eisner}
||\Psi(z_{ j}, t_{ z})\ -\ \Psi(z,t_{ z})||\ <\ \eta
\end{equation}
f"ur alle $j > j^{ \star}$ zutr"afe, sodass f"ur alle $j > j^{ \star}$ die Inklusionenkette
$$\Psi\Bigl(z,\bigcup_{j\in\mathbb{N}} a_{\Psi} (z_{j},\delta) \cdot [-1, 1]\Bigr)+ \mathbb{B}_{\delta+3\eta} (0)\supset$$
$$\Psi\Bigl(z_{j}, a_{\Psi} (z_{j},\delta) \cdot [-1, 1]\Bigr)+ \mathbb{B}_{\delta+\eta} (0)\supset$$
$$\Psi(z_{j},\mathbb{R})+ \mathbb{B}_{\eta} (0)$$
g"alte. Der Zustand $\Psi(z,t(z))$ k"onnte demnach nicht nur kein Element
der 
Trajektorie $\Psi(z_{j},\mathbb{R})$ sein, wenn $j > j^{ \star}$ ist; er k"onnte dar"uber hinaus
auch nicht in der Menge $\Psi(z_{j},\mathbb{R})+ \mathbb{B}_{\eta} (0)$ sein, was
im Widerspruch zu der Ungleichung (\ref{eisner}) steht.
\newline
Existierte hingegen eine positive reelle Zahl $\tau \in \mathbb{R}^{ +}$ von der Art,
dass f"ur eine positive reelle Zahl $\delta \in \mathbb{R}^{ +}$ die Ungleichung
$$a_{\Psi} (z_{j}, \delta)\ \geq\ \ a_{\Psi} (z,\delta)+ \tau$$
f"ur alle $j \in \mathbb{N}$ g"alte, so g"abe es reelle Zahlen
\begin{equation}\label{kurta}
t(z_{j}) \in \bigcap\{a_{\Psi} (z_{ j}, \delta) \cdot [-1,1]:j \in \mathbb{N}\}\setminus(a_{\Psi} (z, \delta) \cdot[-1,1])
\end{equation}
und eine positive reelle Zahl $\eta$, f"ur die f"ur alle $j\in \mathbb{N}$
\begin{equation}\label{kurtb}
\Psi(z_{ j},t(z_{j}))\ \not\in\ \Psi(z,a_{\Psi} (z, \delta)\cdot [-1, 1])+ \mathbb{B}_{\delta+2\eta} (0)
\end{equation}
g"alte, wobei es wegen der Komanenz der Flussfunktion $\Psi$ eine nat"urliche Zahl $j^{\star}\in\mathbb{N}$ von der Art g"abe,
dass f"ur alle $j > j^{ \star}$ die Ungleichung (\ref{eisner}) zutr"afe.
Wegen derselben g"alte daher die Inklusionenkette
$$\Psi(z,a_{\Psi} (z, \delta)\cdot [-1, 1])+ \mathbb{B}_{\delta+2\eta} (0)\supset$$
$$\Psi(z,\mathbb{R})+\ \mathbb{B}_{2\eta} (0)\supset$$
$$\Psi(z,t(z_{j})\cdot [-1, 1])+\ \mathbb{B}_{2\eta} (0)\supset$$
$$\Psi(z_{j},t(z_{j})\cdot [-1, 1])\ni\Psi(z_{ j},t(z_{j}))$$
f"ur alle $j > j^{ \star}$, was aber der Annahme der Existenz der Zahlen
$t(z_{j})$, f"ur welche die Aussagen (\ref{kurta}) und (\ref{kurtb}) gelten, widerspricht: Wenn die
Zustandsfolge $\{z_{ j}\}_{ j \in \mathbb{N}}$ gegen $z$ konvergiert, konvergiert auch die Zahlenfolge
$\{a_{\Psi} (z_{j}, \delta)\}_{ j \in \mathbb{N}}$ gegen die Zahl $a_{\Psi} (z, \delta)$.
\newline
{\bf q.e.d.}
\newline
\newline
Offensichtlich k"onnen komanente und punktweise 
immanente Flussfunktionen $\Psi$ durchaus insensitive Kollektivierungen haben.\index{insensitive Kollektivierung}
Wenn der Zustandsraum $\zeta$ kompakt ist, gibt es also eine pauschale Immanenzfunktion
$A_{\Psi}$
von der Art,
dass die Inklusion
\begin{equation}
\Psi(x,A_{\Psi}(\delta) \cdot [-1,1])\ +\ \mathbb{B}_{ \delta} (0)\ \supset\ \Psi(x, \mathbb{R})
\end{equation}
f"ur alle Zust"ande $x\in \zeta$ wahr ist.
Neben der Immanenzfunktion $A_{\Psi}$ existiert die Funktion 
$b_{ \Psi}$ auf ihrer Definitionsmenge $\mathbf{P}_{1}b_{ \Psi}=]0,\infty[\times]0,\infty[$ 
mit den jeweiligen Werten
$$b_{ \Psi} (\delta , t) :=\ \sup\ \Bigl\{||x\ -\ y||\ :
x, y \in \zeta\quad \wedge\ $$
\begin{equation}
s \in [-t,\ t]\Rightarrow ||\Psi(x,s)-\Psi(y,s)||\ \leq \delta\Bigr\}\ \in\ \mathbb{R}^{ +} ,
\end{equation}
die wir naheliegender Weise die Komanenzfunktion der 
Flussfunktion\index{Komanenzfunktion einer Flussfunktion} 
$\Psi$ nennen.
Dar"uber hinaus existiert die Funktion 
$B_{ \Psi}$ auf $\mathbb{R}^{ +}$ mit der Definitionsmenge $\mathbf{P}_{1}B_{ \Psi}=\mathbb{R}^{ +}$, deren 
Wertemenge im positiven Zahlenstrahl
$\mathbb{R}^{ +}$ liegt und deren 
jeweilige Werte
f"ur alle Zahlen $\delta$ ihrer Definitionsmenge
\begin{equation}
B_{ \Psi}\ (\delta)\ :=\ b_{ \Psi}\ (\delta/3,\ A_{\Psi}\ (\delta/3))
\end{equation}
seien. Diese abgewandelte Komanenzfunktion $B_{ \Psi}$ ist nun aber so beschaffen, 
dass f"ur alle Zust"ande $x, y\in \zeta$ und f"ur alle positiven reellen Zahlen $\delta$ die Implikation
$$||x - y||\ \leq\ B_{ \Psi} (\delta)\quad \wedge\quad
t\ \in\ A_{\Psi}\ (\delta/3)\ \cdot\ [-1,\ 1]\ \Rightarrow$$ 
\begin{equation}
||\Psi(x,t)\ - \ \Psi(y,t)||\ \leq\delta/3
\end{equation}
wahr ist.
Also sind die beiden Inklusionen 
$$\Psi(x, \mathbb{R}) \subset \mathbb{B}_{ \delta/3} (0) + \Psi(x,A_{\Psi} (\delta/3) \cdot [-1,\ 1])$$
und 
$$\Psi(y,\mathbb{R}) \subset \mathbb{B}_{ \delta/3} (0) + \Psi(y,A_{\Psi} (\delta/3) \cdot\ [-1,\ 1])$$
wahr, weshalb
die Absch"atzung
$$\sup \Bigl\{ \inf \Bigl\{||\Psi(x,s) - \Psi(y,t)|| : t\in \mathbb{R}\Bigr\} : s\in \mathbb{R}\Bigr\}\ \leq\delta$$  
f"ur alle $x, y \in \zeta$ mit
$$||x - y|| \leq B_{ \Psi} (\delta)$$
gilt.
Daher gilt aber auch die Implikation 
$$y \in \mathbf{cl} (\Psi(x, \mathbb{R})) \Rightarrow \Psi(y,\mathbb{R})  \subset \mathbf{cl} (\Psi(x, \mathbb{R}))\  ,$$
also sogar die Implikation
$$y \in \mathbf{cl} (\Psi(x, \mathbb{R})) \Rightarrow \mathbf{cl} (\Psi(y,\mathbb{R})) = \mathbf{cl} (\Psi(x, \mathbb{R}))\ .$$ 
Demnach haben wir in der Relation
\begin{equation}\label{zurxy}
N_{ [\Psi]}\ :=\ \Bigl\{(x,y) \in \zeta\times \zeta: y \in \mathbf{cl} (\Psi(x, \mathbb{R}))\Bigr\}
\end{equation}
eine in $\zeta\times \zeta$ enthaltene "Aquivalenzrelation des Zustandsraumes $\zeta$ gefunden.
Und exakt
diese "Aquivalenzrelation (\ref{zurxy}) nennen wir
die N"aherelation der Flussfunktion $\Psi$.\index{N\"ahe-Relation einer Flussfunktion}
Exakt 
die Elemente ihrer Quotientenmenge\index{Notationskonvention}
\begin{equation}\label{trapppp}
\Xi\ \in\ [ \zeta : N_{ [\Psi]}]\ =:\ [[\Psi]]
\end{equation}
nennen wir die 
Zimmer\index{Zimmer} 
der Flussfunktion $\Psi$. Offensichtlich gilt dabei die "Aquivalenz
\begin{equation}
x \in \Xi \in [[\Psi]]\ \Leftrightarrow\ \Xi = \mathbf{cl} (\Psi(x, \mathbb{R}))\ .
\end{equation}
Exakt 
jede Partition $[[\Psi]]$ des Zustandsraumes in Zimmer $\Xi \in [[\Psi]]$
nennen wir gegebenenfalls 
die nat"urliche Partition der jeweiligen Flussfunktion $\Psi$. 
\index{nat\"urliche Partition}\newline
Der Befund, dass die N"aherelation $N_{ [\Psi]}$ existiert und daher deren Quotientenmenge
\begin{equation}
[[\Psi]]=\Bigl\{ \mathbf{cl} (\Psi(x, \mathbb{R})):x\in\mathbf{P}_{2}\Psi\Bigr\}\ \in\ \mathbf{part}(\mathbf{P}_{2}\Psi)
\end{equation}
den Zustandsraum $\mathbf{P}_{2}\Psi$ partioniert, zeigt sich dem Auge der Anschauung offenbar
als eine Gegebenheit, die eigentlich durch die Beschaffenheit alleine des Kollektives $[\Psi]$ 
bedingt sein muss 
und nicht etwa bloss 
durch die Beschaffenheit
spezieller Flussfunktionen, deren Kollektivierung $[\Psi]$ ist.
Die Besinnung auf die "Aquivalenzen (\ref{ovaeva}) und (\ref{evaeva})
begr"undet den Eindruck, der sich f"ur unsere Anschauung ergibt. Es ist daher gerechtfertigt, dass wir die
N"aherelation $N_{ [\Psi]}$
mit dem Kollektiv $[\Psi]$ indizieren, statt
mit einer Flussfunktion.
Die Indizierung mit einer Flussfunktion
haben wir aber bei der 
punktweisen Imanenzfunktion $a_{\Psi}$, der pauschalen Imanenzfunktion $A_{\Psi}$, der Konanenzfunktion $b_{\Psi}$
und schliesslich der Begleitfunktion $B_{\Psi}$ 
zu nehmen, weil
f"ur ein komanent-immanent-beschr"anktes Trajektorienkollektiv
$[\Psi]$ und f"ur
zwei normale Flussfunktionen $\Xi_{1}$ und $\Xi_{2}\not=\Xi_{1}$, die komanent-immanent sind
und f"ur die
$$[\Xi_{1}]=[\Xi_{2}]=[\Psi]$$
gilt, im allgemeinen die Differenzen
\begin{equation}\label{evaeva}
\begin{array}{c}
a_{\Xi_{1}}\not=a_{\Xi_{2}}\ ,\\
A_{\Xi_{1}}\not=A_{\Xi_{2}}\ ,\\
b_{\Xi_{1}}\not=b_{\Xi_{2}}\ ,\\
B_{\Xi_{1}}\not=B_{\Xi_{2}}
\end{array}
\end{equation}
vorliegen.\newline
Zimmer
partionieren also den Zustandsraum jedes komanent-immanenten Trajektorienkollektives,
dessen Zustandsraum $\zeta\subset\mathbb{R}^{n}$ eine beschr"ankte Teilmenge eines 
endlichdimensionalen und reellen Raumes $\mathbf{R}^{n}$ ist,
dessen Dimension $n\in\mathbb{N}$ endlich ist.
Wir fragen uns da nat"urlicherweise, welche ergodentheoretische Verh"altnisse in R"aumen vorliegen,
die nicht endlichdimensional sind und ausserdem, ob die beschr"ankten 
Teilmengen $\zeta\subset\mathbb{R}^{n}$, die
Zustandsr"aume komanent-immanenter Trajektorienkollektive sind, kompakt sein m"ussen.\newline
Der beschr"ankte Zustandsraum $\zeta\subset\mathbb{R}^{n}$ eines komanent-immanenten Trajektorienkollektives
$[\Psi]$ ist
zwar die Vereinigung 
$$\zeta=\bigcup[[\Psi]]\ ,$$
d.h., derjeniger Kompakta $\Xi\in[[\Psi]]$, die die Zimmer sind.
Aber auf dem beschr"ankten und offenen Zustandsraum
$]0,1[$ ist beispielsweise $\{\{x\}:x\in]0,1[\}$ ein komanent-immanentes Trajektorienkollektiv.
\newline
Der Frage nach den ergodentheoretischen Verh"altnissen in R"aumen, deren 
Dimension nicht endlich ist,
gehen wir hier nicht nach. 
Diese Frage
ist dabei eine Spezifizierung der weit allgemeineren Frage
nach Verallgemeinerungen der Zimmer f"ur irgendwelche
topologische R"aume.
Und zu dieser sehr allgemeinen Frage k"onnen wir allerdings schon jetzt,
ausgehend von der Beobachtung, dass die nat"urliche 
Partition eine topologische Invariante ist, 
dies sagen:
Die nat"urliche 
Partition eines komanent-immanenten Trajektorienkollektives $[\Psi]$,
dessen Zustandsraum $\zeta\subset\mathbb{R}^{n}$ beschr"ankt ist, ist insofern eine topologische Invariante,
als f"ur jeden
topologischen Raum $(\bigcup\mathbf{T},\mathbf{T})$ 
und f"ur jeden Hom"oomorphismus $f$, der zwischen demselben 
und dem Zustandsraum $\zeta\subset\mathbb{R}^{n}$ besteht, mit
\begin{equation}\label{avaeva}
f\in\mathcal{C}_{\mathbf{T}(n),\mathbf{T}}(\zeta,\bigcup\mathbf{T})
\end{equation}
offenbar die G"ultigkeit der Implikation
\begin{equation}\label{bvaeva}
\Xi\in[[\Psi]]\Rightarrow f(\Xi)\in [[f\circ\Psi]]\in\mathbf{part}(\bigcup[[f\circ\Psi]])
\end{equation}
bedingt ist, wobei 
$$f(\Xi)=\{f(x):x\in\Xi\}$$
ein jeweiliges Bildzimmer ist, wenn $\Xi$ ein Zimmer ist .
Denn sowohl die Partitivit"at eines Mengensystemes als auch
die Abgeschlossenheit einer Menge ist eine topologische Invariante.
Damit kommen wir zu einer sehr allgemeinen 
ergodentheoretischen
Aussage und passen unsere Begriffe deren Generalit"at an,
indem wir folgende Vereinbarung
treffen:
F"ur jeden topologischen Raum $(\bigcup\mathbf{T},\mathbf{T})$ nennen
wir jedes Mengensystem 
$$\Theta\subset\bigcup\mathbf{T}$$
genau dann ein Trajektorienkollektiv,
wenn es einen 
gem"ass (\ref{avaeva}) beschaffenen
Hom"oomorphismus $f$ gibt,
ferner eine nat"urliche Zahl $n$ und einen
Zustandsraum $\zeta\subset\mathbb{R}^{n}$ einer Flussfunktion $\Psi$
binnen desselben,
die so beschaffen sind, dass dieses Mengensystem
\begin{equation}\label{cvaeva}
\Theta=\{f(\Psi(x,\mathbb{R}):x\in \zeta\}=\Bigl\{\{f(y):y\in\Psi(x,\mathbb{R})\}:x\in \zeta\Bigr\}
\end{equation}
ist. Genau dann, wenn dabei $\zeta$
beschr"ankt ist, nennen wir das Trajektorienkollektiv $\Theta$ beschr"ankt.
Genau dann, wenn es eine komanent-immanente Flussfunktion $\Psi$ gibt, f"ur die
(\ref{cvaeva}) zutrifft, nennen
wir das
Trajektorienkollektiv $\Theta$ komanent-immanent. Die 
Elemente der Partition $f(\Xi)\in [[f\circ\Psi]]$ in (\ref{bvaeva}) nennen wir dabei die 
Zimmer der nat"urlichen Partition $[[f\circ\Psi]]$. 
Jedes solche komanent-immanent-beschr"ankte Trajektorienkollektiv ist ein hom"oomorphes Bild 
einer 
trajektoriellen $\mathcal{C}^{1}$-Partition eines Kompaktums eines endlichdimensionalen reellen Raumes. 
\newline
Zimmer sind ihrer Natur nach kompakt und sie liegen im jeweiligen Zustandsraum $\xi$ separiert.
Sie enthalten keine singul"are Trajektorie $\{x\}\subset\xi$, ohne dieselbe zu sein,
denn singul"are Trajektorien sind offensichtlich selber Zimmer. Also haben sie nach Brouwers Fixpunktsatz\index{Brouwerscher Fixpunktsatz} 
nie das Geschlecht eins.
\newline
Bekanntlich gilt f"ur $A,B \subset \mathbb{R}^{ n}$ im allgemeinen {\em nicht},
dass die Implikation
$$\mathbf{cl} (A) \cap B \neq \emptyset\ \Rightarrow\ \mathbf{cl} (B) \cap A \neq \emptyset\ $$
wahr ist,
auch nicht, wenn $A$ und $B$ speziell Trajektorien sind und wenn dabei die Dimension $n>1$ ist.
Gerade aber das ist anders, wenn 
$A$ und $B$ spezielle Teilmengen eines endlichdimensionalen und reellen 
Raumes sind, n"amlich die
Trajektorien eines speziellen, n"amlich eines immanent-komanenten Trajektorienkollektives
in einem beschr"ankten Zustandsraum.
\newline\newline
Wir heben unser Resultat heraus:
\newline
\newline
{\bf 2.1.2 Der Satz von der Existenz der Zimmer:}\index{Satz von der Existenz der Zimmer} \newline
{\em Der Zustandsraum
eines jeden beschr"ankten, komanent-immanenten Trajektorienkollektives
ist partioniert in die abgeschlossenen H"ullen seiner Trajektorien.}
\newline
\newline 
Trivial, aber wichtig ist nun, dass f"ur beschr"ankte Trajektorien $\tau \subset \mathbb{R}^{ n}$ keineswegs
genau dann
$$\mathbf{cl} (\tau) = \tau$$ 
ist,
wenn die beschr"ankte Trajektorie $\tau$ geschlossen ist, d.h., wenn sie zu dem Kreis $\mathbb{S}^{1}$ relativ hom"oomorph ist.\index{relative Hom\"oomorphie}
Diese "Aquivalenz von 
Geschlossenheit und Abgeschlossenheit
gilt aber,
wenn $\tau$ ein Element einer 
trajektoriellen $\mathcal{C}^{1}$-Partition eines Kompaktums eines endlichdimensionalen reellen Raumes ist.
Die Wichtigkeit dieser "Aquivalenz von Geschlossenheit und Kompaktheit einer beschr"ankten Trajektorie einer 
trajektoriellen $\mathcal{C}^{1}$-Partitionen eines Kompaktums eines endlichdimensionalen reellen Raumes
ist, dass alle nicht geschlossenen Trajektorien einer solchen trajektoriellen Partition
in nicht-trivialen Zimmern $\mathbf{cl} (\tau)$ mit
\begin{equation}
\mathbf{cl} (\tau)\ \setminus\ \tau\ \neq\ \emptyset
\end{equation}
liegen. Dabei gibt es bekanntlich
trajektorielle $\mathcal{C}^{1}$-Partitionen eines Kompaktums eines endlichdimensionalen reellen Raumes, deren
Elemente nicht geschlossene Trajektorien sind.
\newline\newline
Wir nennen f"ur jede Flussfunktion $\Psi$ mit der Komanenzfunktion $b_{\Psi}$ und mit der
Immanenzfunktion $A_{\Psi}$ exakt die Funktion
\begin{equation}\label{nerke}
B_{\Psi}\ :=\ b_{\Psi}\ \Bigl(\mbox{id}/3,\ A_{\Psi}\ (\mbox{id}/3)\Bigr)\ \in\ (\mathbb{R}^{ +})^{\ \mathbb{R}^{ +}}
\end{equation} die
Begleitfunktion der Flussfunktion $\Psi$.\index{Begleitfunktion einer Flussfunktion}\newline \newline
Diese Benennung liegt nahe,
weil $B_{\Psi}$ eine Art von 
Begleitung formuliert.
Die Begleitfunktion $B_{\Psi}$ quantifiziert 
n"amlich diejenige Form von Begleitung, die
objektivierbar ist als die Existenzaussage
\begin{equation}\label{erket}
\begin{array}{c}
\forall\ x,y\ \in\ \lbrack\bigcup \Psi\rbrack,\ t \in\ \mathbb{R}\ \exists\  t(x,y,t)\in\ A_{\Psi}(\delta/3)\cdot[-1,1]:\\
||x - y||\ <\ B_{\Psi} (\delta)\ \Rightarrow\ ||\Psi(x,t) - \Psi(y,t(x,y,t))||\ <\ \delta\ ,
\end{array}
\end{equation}
die von der Aussage zu unterscheiden ist,
dass f"ur alle $x,y \in \lbrack\bigcup \Psi\rbrack,\ t \in\ \mathbb{R}$ einfach
die Implikation
\begin{equation}
||x - y||\ <\ B_{\Psi} (\delta)\ \Rightarrow\ ||\Psi(x,t) - \Psi(y,t)||\ <\ \delta
\end{equation}
gilt,
die im Allgemeinen {\em nicht} wahr ist: Im Vorausblick auf die Sensitivit"at nicht trivialer Zimmer betonen und erl"autern wir
ausf"uhrlich, welche Bewandnis es mit der Begleitfunktion $B_{\Psi}$ der Flussfunktion
$\Psi$ hat: Wenn ein Zustand $y\in\zeta\subset\mathbb{R}^{n}$ des beschr"ankten Zustandsraumes $\zeta$ eines endlichdimensionalen und reellen Raumes $\mathbb{R}^{n}$ 
von einem anderen Zustand $x\in\zeta$ weniger als $B_{\Psi}(\delta)$ entfernt ist, wobei $\delta$ eine positive reelle Zahl und
$\Psi$ eine komanent-immanente Flussfunktion ist und $B_{\Psi}$ deren 
Begleitfunktion,\index{Begleitfunktion einer Flussfunktion} so
bleiben die jeweiligen Zust"ande $\Psi(y,\pm t)$ solange in einer Kugel des Radius $\delta$ um $\Psi(x,\pm t)$,
solange $t\in [0,A_{\Psi}(\delta/3)[$ ist, wobei $A_{\Psi}$ die pauschale Komanenzfunktion der
Flussfunktion $\Psi$ ist. Es gilt dann die Aussage
\begin{displaymath}
\mathbb{B}_{\delta/3}(0)\ +\ \Psi(x,A_{\Psi}(\delta/3)\cdot[-1,1])\ \supset\ \downarrow [[\Psi]]_{\{x\}}\ \land
\end{displaymath}
\begin{equation}\label{klerob}
\mathbb{B}_{\delta/3}(0)\ +\ \Psi(y,A_{\Psi}(\delta/3)\cdot[-1,1])\ \supset\ \downarrow [[\Psi]]_{\{y\}}\ \land
\end{equation}
\begin{displaymath}
t\in  [0,A_{\Psi}(\delta/3)]\Rightarrow||\Psi(x,\pm t)-\Psi(y,\pm t)||<\delta/3\ .
\end{displaymath}
Der sich 
in den Zustand $\Psi(y,t)$ entwickelnde
Zustand $y$ begleitet also die {\em gesamte Trajektorie} $\Psi(x,\mathbb{R})$
insofern, als es f"ur alle Zust"ande $\tilde{x}\in\Psi(x,\mathbb{R})$
einen Zustand $\Psi(y,t(\tilde{x}))\in \Psi(y,A_{\Psi}(\delta/3)\cdot[-1,1])$ gibt, dessen Distanz
$||\tilde{x}-\Psi(y,t(\tilde{x}))||$ kleiner als $\delta$ ist.
Wegen (\ref{klerob}) gilt die Implikation
\begin{equation}\label{kleroc}
\begin{array}{c}
||x-y||<B_{\Psi}(\delta)\Rightarrow\\
\mathbb{B}_{\delta}(0)\ +\ \Psi(x,A_{\Psi}(\delta/3)\cdot[-1,1])\ \supset\ \lbrack\lbrack\Psi\rbrack\rbrack_{\{x\}}\ \land\\
\mathbb{B}_{\delta}(0)\ +\ \Psi(x,A_{\Psi}(\delta/3)\cdot[-1,1])\ \supset\ \lbrack\lbrack\Psi\rbrack\rbrack_{\{y\}}\ .
\end{array}
\end{equation}
F"ur $t\in [0,A_{\Psi}(\delta/3)[$ ist der von der Trajektorie $\Psi(y,\mathbb{R})$
durchlaufene Zustand $\Psi(y,\pm t)$ in der Menge
$$\mathbb{B}_{\delta}(0)\ +\ \Psi(x,t\cdot[-1,1])\ ,$$
die zwar f"ur einen hinreichend kleinen Radius $\delta$ zigarrenf"ormig ist, die aber im allgemeinen recht unf"ormig ist.
Der Sachverhalt, dass die Aussage (\ref{klerob}) f"ur alle positiven reellen Zahlen $\delta$
gilt, impliziert, dass es das Zimmer $\downarrow [[\Psi]]_{\{x\}}$ gibt; und dieser Sachverhalt hat als Konsequenz, dass, falls $y\in \downarrow [[\Psi]]_{\{x\}}$
und $||x-y||<B_{\Psi}(\delta)$ zutrifft, 
die Implikation
\begin{equation}\label{kleroc}
\begin{array}{c}
y\in\mathbf{cl}(\Psi(x,\mathbb{R}))\Rightarrow\\
\lbrack\lbrack\Psi\rbrack\rbrack_{\{x\}}=\mathbf{cl}(\Psi(x,\mathbb{R}))=\mathbf{cl}(\Psi(y,\mathbb{R}))=
\lbrack\lbrack\Psi\rbrack\rbrack_{\{y\}}
\end{array}
\end{equation}
gilt. Die wichtige 
Konsequenz jenes Sachverhaltes ist also die Existenz der Zimmer, 
die der Implikation (\ref{kleroc}) gleichwertig ist. Falls die Bedingung $$y\not\in\Psi(x,\mathbb{R})$$
eingehalten ist, verl"asst der von der Trajektorie $\Psi(y,\mathbb{R})$
durchlaufene Zustand $\Psi(y,\pm t)$ aber gegebenenfalles die Kugel
$$\mathbb{B}_{\delta/3}(\Psi(x,t))$$
wenn $t>A_{\Psi}(\delta/3)$ ist,
gleich, ob der Zustand $y\not\in\Psi(x,\mathbb{R})$ im Zimmer $\downarrow [[\Psi]]_{\{x\}}$ ist, in dem der Zustand $x$ ist,
oder nicht.
Wenn dem nicht so w"are, w"are die dem Satz 2.2.2 gem"asse Sensitivit"at nicht trivialer Zimmer nicht m"oglich,
wie wir sehen werden.
\newline  
Die Aussage (\ref{erket}) ist "aquivalent zu folgender Aussage:
F"ur alle 
Zustandsraumpunkte $x, y \in \bigcup [\Psi]$ und f"ur alle positiven reellen Zahlen $\delta\in \mathbb{R}^{+}$
sind die folgenden Implikationen und "Aquivalenzen
$$\Bigl(\mathbb{B}_{ B_{\Psi} (\delta)} (0)\ +  \Psi(x,\mathbb{R})\Bigr)\  \cap\ \Psi(y,\mathbb{R})\ \neq\ \emptyset\  \Rightarrow$$
$$\mathbb{B}_{ \delta} (0) + \Psi(x,\mathbb{R})\ \supset\ \Psi(y,\mathbb{R})\ ,$$
$$\mathbb{B}_{ \delta} (0) + \Psi(x,\mathbb{R})\ \supset\ \Psi(y,\mathbb{R})\
 \Leftrightarrow\ \mathbb{B}_{ \delta} (0) + \Psi(y,\mathbb{R})\ 
\supset\ \Psi(x,\mathbb{R})\ ,$$
$$\mathbb{B}_{ \delta} (0) + \Psi(y,\mathbb{R})\ 
\supset\ \Psi(x,\mathbb{R})\
\Rightarrow$$
$$\mathbb{B}_{ \delta} (0) + \mathbf{cl} (x\ (\mathbb{R}))\ \supset\ \mathbf{cl} (y\ (\mathbb{R}))$$
wahre Aussagen.
Dies impliziert, dass f"ur alle $x,y\in \bigcup [\Psi]$, $\delta \in \mathbb{R}^{ +}$ die Implikation
\begin{equation}\label{hausdogge}
||x - y||\ <\ B_{\Psi} (\delta)\ \Rightarrow\  \mbox{D}\ \Bigl(\mathbf{cl} (x\ (
\mathbb{R})),\ \mathbf{cl} (y\ (\mathbb{R}))\Bigr)\ <\ \delta
\end{equation}
gilt. Dabei ist
$\mbox{D}$ die Hausdorff-Distanz\index{Hausdorff-Distanz},
eine Funktion mit der Definitionsmenge $2^{\mathbb{R}^{ n}}\times 2^{\mathbb{R}^{ n}}$,  
deren jeweiliger Wert f"ur alle Teilmengen $A, B \subset\mathbb{R}^{ n}$
\begin{equation}\label{hausdog}
\begin{array}{c}
\mbox{D} (A, B)\ :=\inf \Bigl\{\delta\in[0,\infty[:\mathbb{B}_{ \delta} (0) +A\supset B\ \land 
\mathbb{B}_{ \delta} (0) +B\supset A \Bigr\}
\end{array}
\end{equation}
ist. Die Wertemenge der Hausdorff-Distanz ist die Menge $\mathbb{R}\cup\{\infty\}$,
weil das Infimum der leeren Menge mit $\infty$ identifiziert ist. Ferner finden wir die 
"Aquivalenz
\begin{equation}\label{kleroa}
\mathbf{cl}(\Psi(x,\mathbb{R}))=\mathbf{cl}(\Psi(y,\mathbb{R}))\Leftrightarrow \mbox{D}(\Psi(x,\mathbb{R}),\Psi(y,\mathbb{R}))=0\ .
\end{equation} 
Die Restriktion $\mbox{D}|\mathbf{C}(n)$ ist eine Metrik, wenn
\begin{equation}
\mathbf{C}(n):=\mathbf{C}(\mathbf{T}(n))
\end{equation}
das Mengensystem aller abgeschlossenen Mengen ist, die im $\mathbb{R}^{ n}$ liegen;
wobei $\mathbf{C}$ der universelle, 
auf der Klasse der Topologien definierte Operator sei,
der jeder jeweiligen Topologie $\mathbf{T}$ dasjenige Mengensystem $\mathbf{C}(\mathbf{T})$ zuordnet, das das Mengensystem aller bez"uglich $\mathbf{T}$
abgeschlossenen Mengen des zu $\mathbf{T}$ geh"orenden topologischen Raumes ist. 
Der jeweilige
endlichdimensionale und reelle Raum
$\mathbb{R}^{ n}$ sei dabei der, in dem auch alle Trajektorien des
gerade betrachteten Kollektives 
$[\Psi]$ enthalten sind.
Auch wenn die 
nicht restringierte
Hausdorff-Distanz D keine Metrik ist, ist
die Implikation (\ref{hausdogge}) nichtsdestotrotz insofern eine Stetigkeitsaussage, als die 
restringierte Hausdorff-Distanz $\mbox{D}|\mathbf{C}(n)$
sehr wohl eine Metrik ist.\footnote{ 
Eine eingehendere Darstellung der Hausdorff-Distanz
findet man in dem Lehrbuch \cite{topo}.}
Es gilt also als Konsequenz der 
beschriebenen
Begleitung\index{Begleitung} eines komanent-immanenten Trajektorienkollektives $[\Psi]$
mit einem abgeschlossenen Zustandsraum 
$$\mathbf{cl}(\bigcup [\Psi])=\bigcup [\Psi]$$
dies:
Genau dann, wenn eine Folge von Zust"anden 
$\{x_{ j}\}_{ j\in \mathbb{N}} \in (\bigcup [\Psi])^{ \mathbb{N}}$ gegen den Zustand $x$ konvergiert, der
dann ja im abgeschlossenen Zustandsraum $\bigcup [\Psi]$ ist, wobei $\mathbf{cl} (\Psi(x,\mathbb{R}))$ ein Zimmer ist, 
konvergiert eine Zimmernfolge\index{Zimmernfolge}
$\{\mathbf{cl} (x_{ j} (\mathbb{R}))\}_{ j\in \mathbb{N}}$ bez"uglich der Hausdorff-Metrik 
$\mbox{D}|\mathbf{C}(n)$
gegen ein Zimmer der 
nat"urlichen Partition $[[\Psi]]$. Wir haben die 
folgende Bemerkung bereits bewiesen:
\newline
\newline
{\bf 2.1.3 Bemerkung:}\newline {\bf Die Vollst"andigkeit des Raumes von Zimmern}\index{Vollst\"andigkeit des Raumes von Zimmern}\newline
{\em Ist $[\Psi]$ ein komanent-immanentes Trajektorienkollektiv,
so ist
das Paar}  
$$([[\Psi]], \mbox{D})$$
{\em genau dann ein vollst"andiger metrischer Raum, wenn der Zustandsraum 
$$\mathbf{cl}(\bigcup [\Psi])=\bigcup [\Psi]$$
abgeschlossenen ist.} $\mbox{D}$ {\em ist dabei die Hausdorff-Distanz, die von der euklidischen Metrik 
des $\mathbb{R}^{n}$, in dem das Trajektorienkollektiv $[\Psi]$ dargestellt ist,
abgeleitet ist.}
\newline\newline
Dazu, dass $([[\Psi]], \mbox{D})$ ein vollst"andiger metrischer Raum ist, braucht das komanent-immanente Trajektorienkollektiv $[\Psi]$ also nicht notwendigerweise einen
beschr"ankten Zustandsraum 
$\bigcup [\Psi]$
zu haben.\newline
Wir erg"anzen die Bemerkung 2.1.3 ausserdem noch um die folgenden Hinweise:
Der metrische Raum
$(\mathbf{C}(n), \mbox{D})$ 
ist zwar vollst"andig, er
hat aber bekanntlich im generischen Fall keine abz"ahlbare Basis, sodass
er 
im generischen Fall kein polnischer Raum ist.\index{polnischer Raum}
(Zum Beispiel ist die Flussfunktion
$\Phi\in[0,1]^{[0,1]\times \mathbb{R}}$ mit
$$\Phi(x,\mathbb{R})=\{x\}$$
f"ur alle $x\in[0,1]$ so beschaffen, dass der dieser Flussfunktion $\Phi$ entsprechende Raum der Zimmer 
$$([[\Phi]],\ \mbox{D})=\Bigl(\Bigl\{\{x\}:x\in[0,1]\Bigr\},\ |\mathbf{P}_{1}\downarrow-\mathbf{P}_{2}\downarrow|\Bigr)$$
ist. Wobei jede Basis dieses vollst"andigen metrischen Raumes dessen Tr"agermenge $\{\{x\}:x\in[0,1]\}$
enth"alt, deren Kardinalit"at $\aleph_{1}$, die M"achtigkeit des Kontinuums, 
ist. Der Operator 
$\downarrow$ bezeichnet dabei das Element jeder jeweiligen einelementigen Menge.) 
Daher k"onnen
auf der Grundlage 
des im allgemeinen nicht polnischen Raumes $([[\Psi]], \mbox{D})$
stochastische Prozesse,
die die "ubergeordnete Zufallsdynamik der Zimmer beschreiben,
nicht in der "ublichen Weise eingerichtet werden.
Eine derartige 
"ubergeordnete Zufallsdynamik der Zimmer, die
von der deterministischen Evolution, die die jeweiligen Zimmer pr"agt, 
unabh"angig ist,
modellierte
dabei die externe statistische St"orung eines deterministischen Systemes,
die dem angedeuteten Unabh"angigkeitskriterium gen"ugt.
\newline
Weil die Begleitfunktion monoton w"achst und sie im Ursprung verschwindet, gilt
folgende Feststellung, die wir sp"ater noch gebrauchen:
Ist $[\Psi]$ eine 
trajektorielle $\mathcal{C}^{1}$-Partition eines Kompaktums eines endlichdimensionalen reellen Raumes,
so ist
$$A:\ (\bigcup [\Psi])^{\ 2}\ \to\ \mathbb{R}^{ +}_{ 0}\ ,$$
\begin{equation}\label{vonvon}
(x, y)\ \mapsto\ A\ (x,y)\ :=\ \max\Bigl\{
||\Psi(x,t)-\Psi(y,t)||:\ t\in\ \mathbb{R}\Bigr\}
\end{equation}
eine stetige Funktion; sie existiert wegen der Kompaktheit des Zustandsraumes $\bigcup[\Psi]$ und wegen den Stetigkeiten der 
jeweiligen
Funktionen 
$\Psi(x,\mbox{id})$ f"ur alle Zust"ande $x$ des 
Zustandsraumes.
Also ist auch das Paar $(\bigcup [\Psi], A)$ ein metrischer Raum, der vollst"andig ist, weil 
der Zustandsraum 
$\bigcup [\Psi]$ 
der trajektoriellen $\mathcal{C}^{1}$-Partition $[\Psi]$
kompakt ist.\newline\newline
Sehen wir nun von der "Ubertragbarkeit der Zimmer in hom"oomorphe Bilder
endlichdimensionaler und reeller Zustandsr"aume einmal ab, die offenbar in die generalisierte Theorie der Zimmer f"uhrt.
\section{Unverhofftes Wiedersehen: Zimmer sind sensitive Attraktoren.}
\subsection{Was sind sensitive Attraktoren?}
Wiedererkennen macht Freude:
Um den Bezug zu der in der Literatur "ublichen Redeweise ganz unmittelbar herzustellen,
in der vielfach auch diskrete dynamische Systeme betrachtet werden,
beziehen wir uns auf den sogenannten Phasenfluss.\index{Phasenfluss}
F"ur jede 
Flussfunktion $\Psi$ mit $\Psi(\mbox{id},0)=\mbox{id}$
und f"ur jede reelle Zahl $t$ sei
die Abbildung 
\begin{equation}\label{kurzza}
\begin{array}{c} 
\Psi(\mbox{id},t):\zeta\to\zeta\ ,\\
z\mapsto\Psi(z,t)
\end{array}
\end{equation}
und deren Notation\index{Notationskonvention} 
als 
\begin{equation}\label{kurzz}
\Psi^{t}:=\Psi(\mbox{id},t)
\end{equation}
festgelegt. Diese Notation (\ref{kurzz}) ist nicht un"ublich und 
dadurch motiviert,
dass f"ur alle reelle Zahl $t_{1}$ und $t_{2}$
\begin{equation}
\begin{array}{c} 
\Psi^{t_{1}+t_{2}}=\Psi(\mbox{id},t_{1}+t_{2})\\
=\Psi(\Psi(\mbox{id},t_{1}),t_{2})=\Psi^{t_{1}}\circ\Psi^{t_{2}}
\end{array}
\end{equation}
gilt, falls die Gruppe $(\{\Psi^{t}:t\in\mathbb{R}\},\circ)$
Abelsch ist. Wir betrachten Flussfunktionen, deren Kollektivierung Heine-Descartesch ist.\index{Heine-Descartessche Kollektivierung}
Solche Flussfunktionen m"ussen allerdings keineswegs Abelsche Gruppen des Phasenflusses haben.
Die so notierte Abbildung $\Psi^{t}$ ist ein Beispiel f"ur 
einen Phasenfluss
in dem Zustandsraum $\zeta$, der nun hier
gerade im 
endlichdimensionalen und reellen Raum $\mathbb{R}^{n}$ liegt. In der 
Literatur ist auch die
sogenannte Gruppe des Phasenflusses
\begin{equation}\label{kurzzb}
\Bigl((\Psi^{t})^{\mathbb{Z}},\circ\Bigr)
\end{equation}
eine Konstruktion, die viel
verwendet wird, die aber auch manchmal selber als 
Phasenfluss bezeichnet wird. 
Wir nennen hier hingegen exakt
jede Gruppe 
\begin{equation}\label{kurzzc}
(\xi^{\mathbb{Z}},\circ)
\end{equation}
eine Phasenflussgruppe, wenn
$\xi$ eine
Bijektion
einer Menge $Z$ auf sich ist,
wobei $\circ$ die Komposition von Abbildungen sei und $Z$ als ein
allgemeiner Zustandsraum\index{Phasenflussgruppe}
aufgefasst wird.\footnote{Bei manchen Autoren gilt 
die Phasenflussgruppe
aber auch als die Objektivierung 
des Begriffes des dynamischen Systemes, das mit einer 
jeweiligen Phasenflussgruppe identifiziert wird.}
Die Notation der Menge $\xi^{\mathbb{Z}}$ ist dabei 
erl"auterungsbed"urftig,
obwohl unter ihr vielleicht kaum die Menge der Abbildungen der 
ganzen Zahlen $\mathbb{Z}$ in die Menge 
$\xi\subset Z\times Z$, die die Bijektion $\xi$ ja ist, verstanden wird.
Jene Menge der Abbildungen steht nichtsdestotrotz
als $\xi^{\mathbb{Z}}$ eigentlich da.
Und daher sagen wir hier ausdr"ucklich, dass $\xi^{\mathbb{Z}}$ 
folgendermassen festgelegt sei:
Die Schreibweisen
\begin{equation}\label{kurzzd} 
\begin{array}{c} 
\xi^{0}:=\mbox{id}\ ,\\
\xi^{1}:=\xi\ ,\\
k\in\mathbb{Z}\Rightarrow\xi^{k+1}:=\xi\circ\xi^{k}\\
\end{array} 
\end{equation}
sind fraglos etabliert, sodass
gem"ass ebenfalls durchaus "ublicher 
Notationssystematik
\begin{equation}\label{kurzze} 
\xi^{\mathbb{Z}}:=\{\xi^{k}:k\in\mathbb{Z}\}
\end{equation}
die Menge aller sukzessiven Kompositionen der Bijektionen
$\xi$ oder $\xi^{-1}$ der Menge $Z$ auf sich selbst ist. Diese Lesart ist es, die nun hier gelte.
Und dementsprechend ist
$$(\Psi^{t})^{\mathbb{Z}}=\Psi^{t\mathbb{Z}}$$ und auch
\begin{equation}\label{kurzzg} 
\Psi^{\mathbb{R}}:=\{\Psi^{t}:t\in\mathbb{R}\} 
\end{equation}
eine gerechtfertigte Schreibweise.
Die dem Paar $(\Psi,t)$ zugeordnete Phasenflussgruppe
(\ref{kurzzb}) beschreibt demnach ein diskretes dynamisches System
f"ur die spezielle durch den Takt $t\in\mathbb{R}$ gegebene Taktung. 
Ein dynamisches System ist zu seiner Gruppe des Phasenflusses $(\Xi,\circ)$ 
"aquivalent. 
F"ur jedes Element $\xi\in\Xi$ ist die 
f"ur alle $\xi\in\Xi$ gleiche
Definitionsmenge 
$$\downarrow\mathbf{P}_{1}\Xi=\mathbf{P}_{1}\xi$$ der
Zustandsraum des durch $(\Xi,\circ)$ gegebenen dynamisches Systemes,
wobei f"ur jede einelementige Menge $\Xi$
\begin{equation}\label{pontie}
\downarrow \Xi:\in \Xi
\end{equation}
das Element derselben sei.
Wir nennen exakt jedes Paar
\begin{equation}\label{kurzzf} 
\Bigl((\Xi,\circ),\mathbf{T}(\downarrow\mathbf{P}_{1}\Xi)\Bigr)
\end{equation}
ein topologisches dynamisches System oder ein Smale-System,\index{Smale-System}
wenn dabei
$(\downarrow\mathbf{P}_{1}\Xi,\mathbf{T}(\downarrow\mathbf{P}_{1}\Xi))$ ein topologischer Raum ist.
Die Benennung 
topologischer dynamischer Systeme 
nach S. Smale liegt deshalb nahe, weil es
zuvorderst S. Smale war, der in den sp"aten sechziger
Jahren des zwanzigsten Jahrhunderts
z.B. in seiner klassischen Arbeit \cite{male}
die Verbindung zwischen der Theorie
dynamischer Systeme und der mengentheoretischen Topologie 
etablierte.   
\newline 
F"ur jedes Element $\xi\in\Xi$ ist ferner $(\xi^{\mathbb{Z}},\circ)$ eine Untergruppe 
der Phasenflussgruppe
$(\Xi,\circ)$ 
und 
diese Untergruppe $(\xi^{\mathbb{Z}},\circ)$ ist das
"Aquivalent eines eigenen dynamischen Systemes, das aber seinen Zustandsraum 
$\mathbf{P}_{1}\xi=\downarrow\mathbf{P}_{1}\Xi$
mit dem 
durch $(\Xi,\circ)$ charakterisierten dynamischen System gemein hat.\newline
Betrachten wir die Phasenflussgruppe 
$(\xi^{\mathbb{Z}},\circ)$ und
die Fixmengen\footnote{Bei \cite{garr} sind 
sogenannte Attraktoren senso latto\index{Attraktoren senso latto} vorgestellt, die
im wesentlichen die Fixmengen der Abbildung $\xi$ sind.} der Abbildung $\xi$: 
Offensichtlich ist f"ur jede Teilmenge $\tilde{A}\subset\mathbf{P}_{1}\xi$ die Identit"at
$$\tilde{A}=\bigcap_{z\in\mathbb{Z}}\xi^{z}\ (\tilde{A})$$
genau dann\footnote{W"are $\tilde{A}$ eine Menge, f"ur die $\tilde{A}\setminus\xi(\tilde{A})\not=\emptyset$ w"are, dann
l"age $x\in\tilde{A}\setminus\xi(\tilde{A})$ so, dass 
$\xi(x)\not\in \xi^{2}(\tilde{A})$ w"are, sodass $x$ nicht in dem Schnitt der rechten Seite w"are.
Ein Element $x\in \xi(\tilde{A})\setminus\tilde{A}$ hingegen existierte genau dann, wenn
$y:=\xi^{-1}(x)\in \tilde{A}\setminus\xi^{-1}(\tilde{A})$ w"are,
was heisst, dass
$\xi^{-1}(y)\not\in \xi^{-2}(\tilde{A})$ w"are.}
erf"ullt, wenn
$$\xi(\tilde{A})= \tilde{A}$$
gilt, wobei das Mengensystem aller Fixmengen der Abbildung $\xi$
\begin{equation}\label{aqttaa}
\mathbf{T}(\xi):=\{\tilde{A}\subset\mathbf{P}_{1}\xi:\xi(\tilde{A})= \tilde{A}\}
\end{equation}
selber eine $\mbox{T}_{2}$-Topologie des Zustandsraumes ist.
Dieses Mengensystem ist
gerade die Topologie, die wir als die invariante Topologie der Abbildung $\xi$
bezeichnen, die eine Bijektion ihrer Definitionsmenge 
auf dieselbe ist. Exakt eine Bijektion einer Menge 
auf dieselbe nennen wir einen Autobolismus.\index{Autobolismus}\index{invariante Topologie eines Autobolismus}
Die invariante Topologie $\mathbf{T}(\xi)$ jedes Autobolismus $\xi$ erf"ullt das
Trennungsaxiom.\index{Trennungsaxiom}
Eine
Fixmenge 
$A\in\mathbf{T}(\xi)$ des zu $(\xi^{\mathbb{Z}},\circ)$ "aquivalenten dynamischen Systemes,
das die erste Komponente eines Smale-Systemes 
$((\xi^{\mathbb{Z}},\circ),\mathbf{T})$
ist,
heisst genau dann ein
Attraktor des Smale-Systemes $((\xi^{\mathbb{Z}},\circ),\mathbf{T})$, wenn 
es
f"ur alle offenen Teilmengen $a,b\in \mathbf{T}\cap A$ eine ganze Zahl $z$ von der Art gibt,
dass die Reichhaltigkeitsaussage
\begin{equation}\label{durch}
\xi^{z} (a)\cap b\neq\emptyset
\end{equation}
gilt und wenn diese Fixmenge 
$A\in\mathbf{T}(\xi)$ ausserdem bez"uglich der Topologie $\mathbf{T}$ des Smale-Systemes 
$((\xi^{\mathbb{Z}},\circ),\mathbf{T})$ quasikompakt ist.
Vorzugsweise wurden bislang die Attraktoren 
reeller und endlichdimensionaler Zustandsr"aume f"ur Smale-Systeme betrachtet,
deren Topologie die nat"urliche Topologie des jeweiligen Zustandsraumes ist. Deren Attraktoren sind also Kompakta
bez"uglich der nat"urlichen Topologie des jeweiligen Zustandsraumes.
\newline
Das Kriterium
(\ref{durch}) objektiviert also die Koh"arenz innerhalb einer besonderen 
Fixmenge der Abbildung $\xi$, der
ein Attraktor ist. Deswegen nennen wir exakt die Erf"ulltheit dieses Kriteriums
(\ref{durch})
durch eine
Fixmenge 
$A\in\mathbf{T}(\xi)$
die 
Koh"arenz dieser Fixmenge\index{Koh\"arenz-Kriterium} und das Kriterium
(\ref{durch}) selbst das Koh"arenz-Kriterium, das ein Attraktor erf"ullen muss,
um ein solcher zu sein.\footnote{Das Kriterium f"ur sogenannte 
topologische Transitivit"at ist eine Verallgemeinerung des diskret formulierten Koh"arenz-Kriteriums.
(\ref{durch}). Wir werden aber ausserhalb dieser Abhandlung alsbald abstrakte und allgemeine Koh"arenz-Kriterien 
vorstellen.} 
Die Begriffsbildung des 
Attraktors 
kommt dabei offenbar eigentlich auch ohne die Bezugnahme auf eine jeweilige 
festgelegte Zustandsraumtopologie aus.
Attraktoren sind insofern wie 
die jeweiligen 
Fixmengen auch f"ur dynamische Systeme konzipiert, 
die keine Smale-Systeme sind.\index{Smale-System}
In der Literatur hat sich aber die zus"atzliche Forderung der 
Quasikompaktheit eines Attraktors bereits teilweise eingeb"urgert.
\newline
Da offenbar f"ur alle offenen Teilmengen $a,b\subset \chi\in[[\Psi]]$ eine ganze Zahl von der Art existiert,
dass die Reichhaltigkeitsaussage
\begin{equation}\label{durchh}
(\Psi^{t})^{z} (a)\cap b\neq\emptyset
\end{equation}
wahr ist f"ur alle $t\in \mathbb{R}\setminus\{0\}$, ist jedes nicht triviale Zimmer der nat"urlichen Partition
$[[\Psi]]$ ein Attraktor des zu $((\Psi^{t})^{\mathbb{Z}},\circ)$ "aquivalenten dynamischen Systemes,
falls die Implikation
\begin{equation}\label{curch}
\tau\in[\Psi]\Rightarrow\chi\setminus \tau\not=\emptyset
\end{equation}
gilt,
d.h., falls $\chi$ kein triviales Zimmer ist. 
In dem Fall jedoch, dass $\chi$ ein triviales Zimmer ist und das f"ur die 
Nicht-Trivialit"at eines Zimmers notwendige und hinreichende Kriterium
(\ref{curch}) nicht gilt, ist $\chi$ nach der erst noch zu zeigenden Bemerkung A.3
entweder eine singul"are Trajektorie oder aber
eine geschlossene Trajektorie, die eine positive Periode $T_{\chi}\in\mathbb{R}^{+}$ hat.
Wenn $\chi$ ein mit einer singul"aren Trajektorie identisches Zimmer ist,
so gilt die Reichhaltigkeitsaussage (\ref{durchh}) f"ur alle $t\in \mathbb{R}$ auf "ahnliche Weise wie f"ur nicht triviale Zimmer,
jedoch
nicht ganz genauso:
Wenn $\chi$ eine geschlossene Trajektorie mit der Periode $T_{\chi}\in\mathbb{R}^{+}$ ist,
ist $\chi$ 
f"ur jede Taktung $t\in \mathbb{R}\setminus T_{\chi}\mathbb{Z}$
ein Attraktor des zu $((\Psi^{t})^{\mathbb{Z}},\circ)$ "aquivalenten dynamischen Systemes.
\newline
Zimmer sind also auf die beschriebene Weise Attraktoren:
Ausser in dem Fall ihrer nicht singul"aren Geschlossenheit sind sie 
Attraktoren des Smale-Systemes
$$\Bigl((\Psi^\mathbb{R},\circ),\mathbf{T}(\dim\Psi)\Bigr)\ .$$
Ausserdem sind Zimmer insofern elementar,
als es keine kleineren Attraktoren gibt:
Ist die Teilmenge $A\subset\bigcup[[\Psi]]$ ein Attraktor des komanent-immanenten Trajektorienkollektives
$[\Psi]$ in dem Sinn, dass f"ur die sie repr"asentierende Flussfunktion $\Psi$
(\ref{durchh}) gilt, so gilt die Implikation
\begin{equation}
A\subset\chi\in[[\Psi]]\Rightarrow \chi\subset A \ .
\end{equation}
Wir erw"ahnen noch, dass die Eigenschaft eines Zimmers aus $[[\Psi]]$,
ein Attraktor 
des Smale-Systemes
$((\Psi^\mathbb{R},\circ),\mathbf{T}(\dim\Psi))$
zu sein, f"ur alle Flussfunktionen $\Psi$ gilt,
die das komanent-immanente Trajektorienkollektiv $[\Psi]$ repr"asentieren:
\newline
Der Begriff des Zimmers ist eine Abstraktion des Begriffes 
des Attraktors. Ersterer ist n"amlich eine Abstraktion insofern,
als er die definitorische Bindung abstreift,
die der Begriff des Attraktors 
eines dynamischen Systemes noch hat:
Der Begriff des Attraktors ist
an eine spezielle Repr"asentation eines dynamischen Systemes
in Gestalt einer  
Flussfunktion $\Psi$ gebunden. Der Begriff des Zimmers ist hingegen insofern gleichermassen abstrakter und  
geometrisch, als er auf der Gegebenheit trajektorieller Partitionen des Zustandsraumes fusst.\newline 
Es gibt viel Literatur zu Attraktoren und es wurden auch viele Attraktoren studiert. Besonderes Interesse kam 
dabei den sogenannten sensitiven Attraktoren\index{sensitiver Attraktor} zu: 
Die \glqq Sensitivit"at\grqq  der sensitiven Attraktoren ist der Audruck des Erstaunens "uber das Ph"anomen, 
dass die evolution"are Unterscheidbarkeit von Zust"anden
gegeben sein kann, obwohl deren  
Unterschiedenheit nicht aufl"osbar ist.
\newline
Dieses Sensitivit"atsph"anomen\index{Sensitivit\"atsph\"anomen} pr"asentiert sich dabei in logischer
Diversit"at: Es lassen sich verschiedene, keineswegs logisch gleichwertige Erscheinungsbilder der
Sensitivit"at dynamischer Systeme differenziert formulieren.
Diejenige Form der Sensitivit"at, die wir als den Repr"asentanten der etablierten Hauptformen der   
Sensitivit"at dynamischer Systeme auffassen k"onnen, ist die folgende Erscheinung,
die es
genauso, wie es sie f"ur     
Flussfunktionen gibt,
auch allgemeiner f"ur Wellenfunktionen gibt.
Dabei verstehen wir ja gem"ass (\ref{viiergo}) unter einer 
Wellenfunktion\index{Wellenfunktion} jede
Abbildung $\Xi$ aus einer Menge
\begin{equation}\label{calzb}
\mathbf{P}_{2}\Xi^{\mathbf{P}_{1}\mathbf{P}_{1}\Xi\times\mathbf{R}}\ ,
\end{equation}
falls deren Wertemenge 
$\mathbf{P}_{2}\Xi\subset\mathbf{P}_{1}\mathbf{P}_{1}\Xi$ 
im ersten kartesischen Faktor $\mathbf{P}_{1}\mathbf{P}_{1}\Xi$
ihrer Definitionsmenge $\mathbf{P}_{1}\Xi$ enthalten ist, der in einem endlichdimensionalen und reellen Raum liege:
Auch f"ur diese Wellenfunktion $\Xi$ kann deren Sensitivit"at in dem Zustand $x\in \mathbf{P}_{1}\mathbf{P}_{1}\Xi$ auf die Weise
vorliegen, die sich durch die Hilfskonstruktion der Menge
\begin{equation}\label{aalzb}
\begin{array}{c}
M_{\Xi}(x):=
\Bigl\{\delta\in\mathbb{R}^{+}:\ \forall\ \varepsilon\in\mathbb{R}^{+}
\ \exists\ y\in \mathbb{B}_{\varepsilon}(x)\cap\mathbf{P}_{1}\mathbf{P}_{1}\Xi\ \forall T\in\mathbb{R}^{+}\ \exists\ 
t\in ]T,\infty[\\
||\Xi(x,t)-\Xi(y,t)||\ >\ \delta\Bigr\}
\end{array}
\end{equation}
so formulieren l"asst:
Genau dann, wenn die
Menge $M_{\Xi}(x)$ nicht leer ist, nennen wir die Wellenfunktion $\Xi$
im Zustand $x\in \mathbf{P}_{1}\mathbf{P}_{1}\Xi$
sensitiv.\index{Sensitivit\"at einer Wellenfunktion}
Ferner ist die Sensitivit"at der Wellenfunktion $\Xi$ durch 
die Funktion
\begin{equation}\label{aalzbb}
\begin{array}{c}
\Delta_{\Xi}:\mathbf{P}_{1}\mathbf{P}_{1}\Xi\to\mathbb{R}^{+}\cup\{-\infty\}\ ,\\
x\mapsto\Delta_{\Xi}(x):=\sup M_{\Xi}(x)
\end{array}
\end{equation}
quantitativ bewertbar, exakt welche wir fortan 
f"ur jede Wellenfunktion $\Xi$
als deren Aufl"osungsfeld\index{Aufl\"osungsfeld einer Wellenfunktion} bezeichnen wollen. 
Genau dann, wenn die Wellenfunktion $\Xi$ im Zustand 
$x\in \mathbf{P}_{1}\mathbf{P}_{1}\Xi$ nicht sensitiv ist, ist ihr Aufl"osungsfeld in $x$ singul"ar,
n"amlich 
das mit $-\infty$ identifizierte Supremum der leeren Menge. 
Diese allgemein formulierte Sensitivit"at umfasst damit
auch eine der Sensitivit"atsformen einer jeweiligen Flussfunktion $\Xi$
in einem jeweiligen Zustand; und gerade diese 
Sensitivit"atsform einer jeweiligen Flussfunktion
bezeichnen wir als deren zustandsweise
Sensitivit"at im jeweiligen Zustand.\index{zustandsweise Sensitivit\"at einer Flussfunktion}  
\newline
Nehmen wir an, uns werde eine Flussfunktion $\Psi$ vorgelegt, 
deren Zustandsraum $\chi$ in einem endlichdimensionalen und reellen Raum liege
und die "uberdies komanent sei. Solche Flussfunktionen gibt es sehr wohl. 
Komanenz und Sensitivit"at der Flussfunktion $\Psi$ in dem Sinn, dass 
$$\Delta_{\Psi}(\mathbf{P}_{2}\Psi)\not=\{-\infty\}$$
ist, und dass $\Psi$ dabei komanent ist, schliessen einander nicht aus!
F"ur $\Psi$ gelte
die Aussage 
\begin{equation}\label{saalzb}
\begin{array}{c}
\exists x\in \chi\ \forall\ \delta\in\mathbb{R}^{+}\exists\varepsilon\in\mathbb{R}^{+}
\ \forall\ y\in \mathbb{B}_{\varepsilon}(x)\ \exists T\in\mathbb{R}^{+}:\\
t\in ]T,\infty[\Rightarrow
||\Psi(x,t)-\Psi(y,t)||\ <\ \delta\ , 
\end{array}
\end{equation} 
die behauptet,
dass es einen Zustand $x\in \chi$ gibt,
in dem $\Psi$ nicht 
zustandsweise Sensitivit"at zeige.
Dann gibt 
es also erstens die Funktion
\begin{equation}\label{saalzc}
\begin{array}{c}
\varepsilon_{x}:\mathbb{R}^{+}\to\mathbb{R}^{+},\\
\delta\mapsto\varepsilon_{x}(\delta)
\end{array}
\end{equation}
und zu derselben zweitens die Funktion
\begin{equation}\label{saalzd}
\begin{array}{c}
{\rm T}:\mathbb{R}^{+}\times\mathbb{B}_{\varepsilon_{x}(\delta)}(x)\cap \chi \to\mathbb{R}^{+},\\
(\delta,y)\mapsto{\rm T}(\delta,y)\ ,
\end{array}
\end{equation}
deren jeweiliger 
Wert f"ur alle $y\in \mathbb{B}_{\varepsilon_{x}(\delta)}(x)$ das Infimum
\begin{equation}
{\rm T}(\delta,y):=\inf\Bigl\{T\in\mathbb{R}^{+}:t>T\Rightarrow||\Psi(x,t)-\Psi(y,t)|| < \delta\Bigr\}
\end{equation}
sei. 
Nehmen wir des weiteren an, dass diese Funktion 
${\rm T}$
stetig sei:
Dann ist das Supremum
\begin{equation}\label{salzdd}
{\rm T}^{\star}(\delta):=
\sup {\rm T}(\delta,\mathbb{B}_{\varepsilon_{x}(\delta)}(x))\ \in\ \mathbb{R}^{+}
\end{equation}
eine positive reelle Zahl,
f"ur die die Implikation
$$(y,t)\in \mathbb{B}_{\varepsilon_{x}(\delta)}(x)\times ]{\rm T}^{\star}(\delta),\infty[
\ \Rightarrow$$ 
$$||\Psi(x,t)-\Psi(y,t)||\ <\ \delta$$ 
gilt.
Wegen der Komanenz der Flussfunktion $\Psi$ gilt
aber auch andererseits, dass es eine Funktion 
\begin{equation}\label{salzf}
\begin{array}{c}
\varepsilon^{x}:\mathbb{R}^{+}\times\mathbb{R}^{+} \to\mathbb{R}^{+},\\
(\delta,T)\mapsto\varepsilon^{x}(\delta,T)
\end{array}
\end{equation}
gibt, 
f"ur die
f"ur alle Paare $(\delta,T)\in\mathbb{R}^{+}\times\mathbb{R}^{+}$ die Implikation
$$(y,t)\in \mathbb{B}_{\varepsilon^{x}(\delta,T)}(x)\times [0,{\rm T}^{\star}(\delta)]\ \Rightarrow$$ 
$$||\Psi(x,t)-\Psi(y,t)||\ <\ \delta$$
wahr ist, sodass f"ur alle $(\delta,T,t)\in\times\mathbb{R}^{+}\times\mathbb{R}^{+}\times\mathbb{R}^{+}$ die Implikation
$$||y-x||<\min\{\varepsilon^{x}(\delta,{\rm T}^{\star}(\delta)),\varepsilon_{x}(\delta)\}\ \Rightarrow$$
$$||\Psi(x,t)-\Psi(y,t)||\ <\ \delta$$
gilt. 
Diese Flussfunktion
$\Psi$ ist also, wenn die Annahme (\ref{saalzb})
zutrifft, so beschaffen, dass die durch $x$ verlaufende Trajektorie $\Psi(x,\mathbb{R})$ mit ihrer
abgeschlossenen H"ulle
$\mathbf{cl}(\Psi(x,\mathbb{R}))$ identisch ist,
wie sich leicht zeigen l"asst.
Ist $\chi$ ein nicht triviales Zimmer, so kann daher die Annahme (\ref{saalzb})
nicht gelten.
Diese Beobachtung legt uns geradezu aufdringlich nahe, zu versuchen, die
zustandsweise Sensitivit"at der nicht trivialen Zimmer
komanent-imanenter
Flussfunktionen binnen beschr"ankter und endlichdimensionaler und reeller
Zustandsr"aume zu zeigen.
Und zwar, indem 
wir die Stetigkeit 
der jeweiligen Entsprechnung zu der Funktion ${\rm T}$ gem"ass
(\ref{saalzd})
beweisen. Dieser Weg ist vielleicht gangbar. Er ist jedoch steinig und dieser Weg liegt im G"ahnen der Gefahr,
dass die seiner Idee folgende Argumentation nicht l"uckenlos ist.
\subsection{Die zustandsweise Sensitivit"at nicht trivialer Zimmer}
Um innerhalb eines allgemeineren Rahmens einzusehen,
dass nicht triviale Zimmer die zustandsweise Sensitivit"at haben,
betrachten wir wieder 
nicht speziell eine Flussfunktion, sondern 
eine Wellenfunktion $\Xi$, deren Wertemenge
$\mathbf{P}_{2}\Xi\subset\mathbf{P}_{1}\mathbf{P}_{1}\Xi$, wie vorhin,
im zweiten kartesischen Faktor $\mathbf{P}_{1}\mathbf{P}_{1}\Xi$
ihrer Definitionsmenge $\mathbf{P}_{1}\Xi$ enthalten ist, der in einem endlichdimensionalen und reellen Raum liege.
$\Xi$ sei aber diesesmal eine stetige Wellenfunktion.
Zu dieser $\mathcal{C}$-Wellenfunktion $\Xi$ gibt es die durch sie bestimmte sekund"are Funktion in jedem Punkt 
$x\in\mathbf{P}_{1}\mathbf{P}_{1}\Xi$
\begin{equation}\label{salzd}
\begin{array}{c}
{\rm T}[\Xi,x]:\mathbb{R}^{+}\times\mathbf{P}_{2}\Xi \to\mathbb{R}^{+}\cup\{\infty\},\\
(\delta,y)\mapsto{\rm T}[\Xi,x](\delta,y)\ ,
\end{array}
\end{equation}
deren jeweiliger 
Wert f"ur alle $y\in\mathbf{P}_{1}\mathbf{P}_{1}\Xi$ das Infimum
\begin{equation}
{\rm T}[\Xi,x](\delta,y):=\inf\Bigl\{T\in\mathbb{R}^{+}:t>T\Rightarrow||\Xi(x,t)-\Xi(y,t)|| < \delta\Bigr\}
\end{equation}
sei. Wegen der Stetigkeit der Wellenfunktion $\Xi$ gilt daher die Gleichung
\begin{equation}
||\Xi(x,{\rm T}[\Xi,x](\delta,y))-\Xi(y,{\rm T}[\Xi,x](\delta,y))||\ =\ \delta\ ,
\end{equation}
die den Wert 
${\rm T}[\Xi,x](y,\delta)$ 
als absolutes Maximum der Restriktion 
$$(||\Xi(x,\mbox{id})-\Xi(y,\mbox{id})||)|[{\rm T}[\Xi,x](\delta,y),\infty[ $$
auf das Intervall $[{\rm T}[\Xi,x](\delta,y),\infty[$ kennzeichnet.
Dieser Wert
${\rm T}[\Xi,x](y,\delta)$
ist eine Stelle des Zahlenstrahles, in der die Funktion
$\Xi_{x,y}$ monoton
f"allt, wobei $\Xi_{a,b}$
f"ur jedes Paar $(a,b)\in\mathbf{P}_{1}\mathbf{P}_{1}\Xi\times\mathbf{P}_{1}\mathbf{P}_{1}\Xi$ 
die auf dem Zahlenstrahl als 
$$\Xi_{a,b}:=||\Xi(x,\mbox{id})-\Xi(y,\mbox{id})||$$ 
erkl"arte, nicht negative, reelle Funktion sei.
Einfachheitshalber
betrachten wir nun eine Wellenfunktion $\Phi$, die ansonsten wie $\Xi$ beschaffen sei,
ausser, dass sie "uberdies eine $\mathcal{C}^{0,1}$-Wellenfunktion\footnote{$\mathcal{C}^{0,1/2}$-Wellenfunktion seien
exakt diejenigen Wellenfunktionen, die in der ersten Ver"anderlichen
stetig und in der zweiten Ver"anderlichen st"uckweise stetig differenzierbar sind und
$\mathcal{C}^{0,-1/2}$-Wellenfunktion seien dementsprechend
exakt diejenigen Wellenfunktionen, die in der ersten Ver"anderlichen
stetig und in der zweiten Ver"anderlichen st"uckweise stetig sind.} sei,
womit bezeichnet sei, dass $\Phi$ in der ersten Ver"anderlichen
stetig und in der zweiten Ver"anderlichen stetig differenzierbar sei.
Dann ist das streng monotone Fallen der $\Xi_{x,y}$ entsprechenden 
Funktion $\Phi_{x,y}$ "ubersichtlich zu erfassen. Es ist
genau dann der Fall, wenn das ebenfalls auf dem ganzen Zahlenstrahl existierende 
Produkt dieser  
Funktion $\Phi_{x,y}$ mit ihrer Ableitung
$$\Phi_{x,y}\Phi_{x,y}'$$
an der Stelle ${\rm T}[\Phi,x](y,\delta)$ einen negativen Wert
annimmt. Dass $\Phi_{x,y}\Phi_{x,y}'$ an dieser 
Stelle ${\rm T}[\Phi,x](y,\delta)$ dabei einen positiven Wert
hat, widerspricht der Konstruktion des Infimums ${\rm T}[\Phi,x](y,\delta)$.
Es ist zwar die Positivit"at 
$$\Phi_{x,y}\Phi_{x,y}'|[{\rm T}[\Phi,x](\delta,y),\infty[>0$$
nicht m"oglich, aber
es kann sein, dass
der kritische Fall vorliegt, dass 
\begin{equation}\label{slzcra}
\Phi_{x,y}\Phi_{x,y}'({\rm T}[\Phi,x](y,\delta))=0
\end{equation}
ist. (In ${\rm T}[\Phi,x](y,\delta)$
ist dann aber dennoch nicht notwendigerweise ein lokales Maximum der 
Funktion $\Phi_{x,y}$.)
Betrachten wir nur den unkritischen Fall, dass 
\begin{equation}\label{slzcrb}
\Phi_{x,y}\Phi_{x,y}'({\rm T}[\Phi,x](y,\delta))<0 
\end{equation}
ist. Es sei gem"ass unserer Notation (\ref{not1}) von Auswahlen durch eine Menge
$$J\in\mathbf{T}(1)_{\{{\rm T}[\Phi,x](y,\delta)\}}$$
eine hinreichend kleine Umgebung des Punktes ${\rm T}[\Phi,x](y,\delta)$, sodass
$||\Phi(x,t)-\Phi(y+h_{k},\mbox{id})|||J$ die Restriktion
der Funktion $||\Phi(x,t)-\Phi(y+h_{k},\mbox{id})||$ auf jene offene Teilmenge $J\subset\mathbb{R}$, mithin selber eine Funktion ist, und
der Wert
$$(||\Phi(x,t)-\Phi(y+h_{k},\mbox{id})|||J)^{-1}(\delta)$$
eine reelle Zahl, die in $J$ ist. Dieser Wert
verschiebt sich in diesem Fall gegen"uber ${\rm T}[\Phi,x](y,\delta)$ beliebig wenig, wenn 
$k$ ein hinreichend hoher Folgenindex 
der Nullfolge $$\{h_{k}\}_{k\in\mathbb{N}}\in\Bigl(\mathbb{R}^{\dim\Phi}\Bigr)^{\mathbb{N}}$$ ist, die
so beschaffen sei, dass die Aussage
\begin{equation}\label{salzbd}
\begin{array}{c}
k\in\mathbb{N}\Rightarrow y+h_{k}\in\mathbf{P}_{2}\Phi\ \land\\
\lim_{k\to\infty}||h_{k}||=0
\end{array}
\end{equation}
gelte, wobei $\dim\Phi$ die Komponentenzahl der $\mathcal{C}^{0,1}$-Wellenfunktion $\Phi$ sei.
\newline
Selbst in diesem unkritischen Fall der besagten strengen Monotonie 
gem"ass (\ref{slzcrb})
ist die Funktion ${\rm T}[\Phi,x](\delta,\mbox{id})$ dennoch nicht notwendigerweise
in einer 
Umgebung des Punktes $x\in \mathbf{P}_{1}\mathbf{P}_{1}\Phi$
stetig:
Die Funktion
${\rm T}[\Phi](y,\mbox{id})$ ist offenbar nicht stetig,
wenn es eine streng monoton wachsende Folge 
$$\{t_{k}\}_{k\in\mathbb{N}}\in ([{\rm T}[\Phi,x](\delta,y),\infty[)^{\mathbb{N}}$$ gibt, 
f"ur die die Grenzwertaussage
\begin{equation}\label{salzdd}
\lim_{k\to\infty} \Phi(y,t_{k})=
\delta
\end{equation}
gilt. Und dass wir dabei diese M"oglichkeit nicht ausschliessen k"onnen, belegen
entsprechende Beispiele, deren Existenz evident ist. Zwei Sacherhalte sind auseinanderzuhalten:
Die Tatsache, dass es 
in dem unkritischen Fall der besagten strengen Monotonie 
gem"ass (\ref{slzcrb}) 
eine Umgebung $J$ des Punktes ${\rm T}[\Phi,x](y,\delta)$ gibt, f"ur die 
der Grenzwert
\begin{equation}\label{lzsdd}
\lim_{k\to\infty}{\rm T}[\Phi,x](y,\delta)-(||\Phi(x,t)-\Phi(y+h_{k},\mbox{id})|||J)^{-1}(\delta)=0
\end{equation}
ist, differiert von der Stetigkeitsbehauptung, dass
\begin{equation}\label{lzsda}
\lim_{k\to\infty}{\rm T}[\Phi,x](y,\delta)-{\rm T}[\Phi,x](y+h_{k},\delta)=0
\end{equation}
ist.
Die entprechende Sachverhaltsdifferenz zwischen
(\ref{lzsdd}) und (\ref{lzsda})
besteht offenbar auch f"ur eine $\mathcal{C}$-Wellenfunktion $\Xi$ noch so hoher
Kontinuit"at:
Wir k"onnen 
gegebenenfalls auch aus der 
noch so hohen
Kontinuit"at der Wellenfunktion $\Xi$ alleine noch
keinen Schluss auf die Stetigkeit der Funktion ${\rm T}[\Xi,x](\mbox{id},\delta)$
ziehen. In der Funktion ${\rm T}[\Xi,x](\mbox{id},\delta)$ stossen wir also auf ein 
Epikonstrukt zu einer  
Wellenfunktion $\Xi$, 
das insofern bemerkenswert ist, als
$\Xi$
von noch so hoher Kontinuit"at sein kann und dennoch ist ${\rm T}[\Xi,x](\mbox{id},\delta)$ nicht stetig:
Die potenzielle
Differenz zwischen dem Sachverhalt 
(\ref{lzsdd}) und der Stetigkeit (\ref{lzsda}) ist offenbar eine Frage der
Kontinuit"at 
des zu der
Wellenfunktion $\Xi$ konstruierten sekund"aren Epikonstruktes ${\rm T}[\Xi,x](\mbox{id},\delta)$, die
aber keine alleinige Frage der
Kontinuit"at der prim"aren Wellenfunktion $\Xi$ ist.
Dieses Epikonstrukt ${\rm T}[\Xi,x](\mbox{id},\delta)$ kann 
dennoch -- selbst in dem unkritischen Fall der (\ref{slzcrb}) entsprechenden strengen Monotonie --
sprunghaft sein.\footnote{Das wundert nicht nur den Fachmann wenig. Nichtsdestotrotz halten wir die 
gew"ahlte Benennung f"ur motiviert.} 
Demnach ist die folgende Benennung motiviert:
F"ur jede $\mathcal{C}$-Wellenfunktion $\Xi$, deren Wertemenge 
$\mathbf{P}_{2}\Xi\subset\mathbf{P}_{1}\mathbf{P}_{1}\Xi$ 
im zweiten kartesischen Faktor $\mathbf{P}_{1}\mathbf{P}_{1}\Xi$
ihrer Definitionsmenge $\mathbf{P}_{1}\Xi$ enthalten ist, der in einem endlichdimensionalen und reellen Raum liegt und
f"ur jeden Punkt $x\in \mathbf{P}_{1}\mathbf{P}_{1}\Xi$ nennen wir 
exakt die 
Funktion ${\rm T}[\Xi,x]$ die atque-fecit-saltus - Funktion
zu der Wellenfunktion $\Xi$ an der Stelle $x\in \mathbf{P}_{1}\mathbf{P}_{1}\Xi$.\index{atque-fecit-saltus - Funktion}
\newline
Die Stetigkeit einer 
jeweiligen atque-fecit-saltus - Funktion
${\rm T}[\Xi,x]$ in ihrer ersten Ver"anderlichen wird
nicht durch die 
Kontinuit"at der Wellenfunktion $\Xi$ alleine
bestimmt. Sie wird n"amlich von der Kollektivierungsstruktur 
der Wellenfunktion $\Xi$ in hohem Grad mitbestimmt.
Dies zeigt die folgende 
\newline\newline
{\bf Bemerkung 2.2.1:\newline 
Die Unkritische Stetigkeit einer atque-fecit-saltus-Funktion}\index{Stetigkeit von atque-fecit-saltus-Funktionen}\newline
{\em Es sei $\Xi$
eine Wellenfunktion, die speziell eine komanent-immanente
$\mathcal{C}^{0,1}$-Flussfunktion $\Xi$ ist, deren Zustandsraum $\mathbf{P}_{2}\Xi$
beschr"ankt ist und $\chi\in[[\Xi]]$ ein Zimmer.
Die Spezifizierung der zu $\Xi$ geh"orenden atque-fecit-saltus - Funktion ${\rm T}[\Xi,x](\delta,\mbox{id})$ 
f"ur $\delta\in\mathbb{R}^{+}$ und f"ur einen Zustand $x\in \chi$ ist
auf
der Teilmenge 
\begin{equation}\label{shipa}
{\rm U}_{\delta}^{\Xi}(x):=\Bigl\{y\in\mathbb{B}_{{\rm T}[\Xi,x]
(\delta,y)}(x)\cap\chi:\Xi_{x,y}\Xi_{x,y}'({\rm T}[\Xi,x](y,\delta))\not=0\Bigr\}
\end{equation}
stetig.}
\newline\newline
{\bf Beweis:}\newline
Das Kompaktum $\chi\in[[\Xi]]$ ist ein Zimmer.
Deswegen gelingt es uns zu zeigen, dass die Funktion
${\rm T}[\Xi,x](\mbox{id},\delta)|\chi$ f"ur die positive reelle Zahl $\delta$
und den Zustand $x\in \chi$
im Punkt $y\in \chi\in[[\Xi]]$ in dem Fall tats"achlich stetig ist, der der
in dem Sinn unkritische Fall ist, dass in $y\in \chi\in[[\Xi]]$
\begin{displaymath}
\Xi_{x,y}\Xi_{x,y}'({\rm T}[\Xi,x](y,\delta))<0 
\end{displaymath}
gilt:
\newline
Weil $y$ in dem selben Zimmer $\chi$ ist, in dem auch $x$ ist,
schneidet die Menge
$$\Xi(y,[{\rm T}[\Xi,x](\delta,y),\infty[)$$
jede Kugel
$$\mathbb{B}_{\delta-\vartheta}(x)$$
um $x$, gleich, wie klein deren Radius $\delta-\vartheta\in\mathbb{R}^{+}$ ist.
F"ur jede positive reelle Zahl $\vartheta\in ]0,\delta[$
muss es also die positive reelle Zahl $\Theta(y,\vartheta)$
geben, f"ur die
f"ur alle $t\in[\Theta(y,\vartheta),\infty[$
die Ungleichung
$$\Xi_{x,y}(t)\ <\ \delta-\vartheta$$
erf"ullt sein muss, was mit dem 
Szenario nicht vereinbar ist, das dem Zutreffen
der Grenzwertaussage (\ref{salzdd}) entspricht.
Die entsprechende Zahl $\Theta(y_{k},\vartheta)$ gibt es dabei auch f"ur alle Folgenindizes 
$k\in \mathbb{N}$ einer gegen $y$ konvergenten Folge
$$\{y_{k}\}_{k\in\mathbb{N}}\in \chi^{\mathbb{N}}\ .$$
Es gilt sowohl die Implikation
$$t>{\rm T}[\Xi,x](\delta,y)\Rightarrow\Xi_{x,y}(t)\ <\ \delta$$
als auch die Implikation
$$t>\Theta(y,\vartheta)\Rightarrow\Xi_{x,y}(t)\ <\ \delta-\vartheta$$
und daher weicht der Funktionswert ${\rm T}[\Xi,x](\delta,y_{k})$
von dem Funktionswert ${\rm T}[\Xi,x](\delta,y)$
beliebig wenig ab, wenn $||y-y_{k}||$ hinreichend klein ist.
\footnote{Denn es gibt die Zahl $\alpha$
im Intervall $$]0,{\rm T}[\Xi,x](\delta,y)-\Theta(y,\vartheta)[\ ,$$ f"ur die
das Maximum 
$$\max\Bigl\{\delta-\Xi_{x,y}(t):t\in[{\rm T}[\Xi,x](\delta,y)-\alpha,\Theta(y,\vartheta)]\Bigl\}$$
$$=\delta-\Xi_{x,y}({\rm T}[\Xi,x](\delta,y)-\alpha)> 0$$
positiv ist.
Also gibt es eine nat"urliche Zahl $k^{\star}\in \mathbb{N}$, f"ur die
f"ur $k>k^{\star}$ die Norm
$||y-y_{k}||$ so klein ist,
dass f"ur alle
$t\in [{\rm T}[\Xi,x](\delta,y)-\alpha,\Theta(y,\vartheta)]\cup[\Theta(y,\vartheta),\Theta(y+h_{k},\vartheta)]$ die Ungleichung 
$$\Xi_{x,y_{k}}(t)<\delta$$
gilt.}
Also ist die Funktion ${\rm T}[\Xi,x](\delta,\mbox{id})$ auf
der gem"ass (\ref{shipa}) festgelegten Teilmenge 
${\rm U}_{\delta}^{\Xi}(x)$ 
der relativen Kugel $\mathbb{B}_{{\rm T}[\Xi,x](\delta,y)}(x)\cap\chi$ stetig.\newline
{\bf q.e.d.}
\newline\newline
Bedenken wir nun noch einmal die Feststellung, die wir f"ur die komanente Flussfunktion $\Psi$ trafen, 
ehe wir dann anschliessend die atque-fecit-saltus - Funktionen einf"uhrten.
Diese Feststellung war, dass, wenn die Annahme (\ref{saalzb}) f"ur die Flussfunktion
$\Psi$ in $x\in\mathbf{P}_{2}\Psi$ zutrifft,
kein nicht triviales Zimmer $\chi\in[[\Psi]]$ den Zustand $x$
als Element hat:
Wir erkennen dann, dass unsere kurze Untersuchung der atque-fecit-saltus - Funktionen andeutet,
dass das zustandsweise Sensitivit"atsph"anomen in einem Zustand $x$ bei 
Flussfunktionen $\Psi$, deren Kollektivierung eine nat"urliche Partion in Zimmer
bestimmt,
als singul"ares Verhalten verstanden werden kann;
n"amlich als das in dem Sinn 
singul"are Verhalten 
der zu der jeweiligen Flussfunktion $\Psi$ gegebenen atque-fecit-saltus - Funktion ${\rm T}[\Psi,x]$, die exakt im Fall
der zustandsweisen Sensitivit"at der Flussfunktion $\Psi$ im Zustand $x$ nicht beschr"ankbar ist.
\newline
Wir f"ugen dieser Bemerkung 2.2.1 noch die tiviale Perspektive hinzu, dass wir die Behauptung des folgenden Satzes 2.2.2 zeigen,
falls wir
f"ur jeden kritischen Zustand
$$y^{\star}\in\mathbb{B}_{\varepsilon_{x}(\delta)}(x) \setminus{\rm U}_{\delta}^{\psi}(x)$$
und f"ur jede gegen denselben konvergente Folge
$$\{y_{k}\}_{k\in\mathbb{N}}\in ({\rm U}_{\delta}^{\Psi}(x))^{\mathbb{N}}$$
die Grenzwertaussage
$$\lim_{k\to\infty}{\rm T}[\Xi](y_{k},\delta)={\rm T}[\Xi](y^{\star},\delta)$$
zeigen k"onnen. Diesen Weg m"ussen wir aber nicht einschlagen.
\newline\newline
Wir finden,
dass nicht triviale Zimmer eines komanent-immanenten Trajektorienkollektives 
in Sinne des folgenden Satzes 
sensitive Attraktoren sind:
\newline\newline
{\bf Satz 2.2.2: Die Sensitivit"at nicht trivialer Zimmer:}\index{Satz von der Sensitivit\"at nicht trivialer Zimmer}\newline
{\em Alle nicht trivialen Zimmer $\chi\in[[\Psi]]$, die Zimmer der nat"urlichen Partition
eines beschr"ankten, komanent-immanenten Trajektorienkollektives $[\Psi]$ sind, sind f"ur alle
es repr"asentierenden Flussfunktionen $\Psi$ 
sensitive Attraktoren in dem Sinn, dass es eine positive reelle Zahl $Q(\chi,\delta)$ von der Art gibt,
dass
f"ur alle $t\in]0,Q(\chi,\delta)[$ und f"ur alle nicht trivialen Zimmer $\chi\in[[\Psi]]$}
\begin{equation}\label{salz}
\begin{array}{c}
\forall x\in \chi\exists\delta\in\mathbb{R}^{+}\forall\varepsilon\in\mathbb{R}^{+}
\exists y\in \mathbb{B}_{\varepsilon}(x),\ z\in\mathbb{Z}:\\
||(\Psi^{t})^{z}(x)-(\Psi^{t})^{z}(y)||\ >\ \delta\ 
\end{array}
\end{equation} 
{\em gilt.}\newline
\newline
\newline
{\bf Beweis:}\newline
Es sei $\chi\in[[\Psi]]$ ein nicht triviales Zimmer,
das wegen der Beschr"anktheit des komanent-immanenten Trajektorienkollektives $[\Psi]$ kompakt ist.
\newline 
Unser Beweis ist indirekt und er besteht vor allem darin, zu zeigen, dass
die Annahme, dass $\chi$ kein sensitiver Attraktor ist, im Widerspruch 
zur Nichttrivialit"at des Zimmers $\chi$ steht, von dem wir bereits wissen, dass es ein Attraktor ist:
\newline 
Nehmen wir an, es gelte die Aussage 
\begin{equation}\label{salzb}
\begin{array}{c}
\exists x\in \chi\ \forall\ \delta_{1}\in\mathbb{R}^{+}\exists\varepsilon\in\mathbb{R}^{+}
\ \forall\ y\in \mathbb{B}_{\varepsilon}(x)\ \exists T\in\mathbb{R}^{+}:\\
t\in ]T,\infty[\Rightarrow
||\Psi(x,t)-\Psi(y,t)||\ <\ \delta_{1}\ . 
\end{array}
\end{equation}
Dann gibt 
es also erstens die Funktion
\begin{equation}\label{salzc}
\begin{array}{c}
\varepsilon_{x}:\mathbb{R}^{+}\to\mathbb{R}^{+},\\
\delta\mapsto\varepsilon_{x}(\delta)
\end{array}
\end{equation}
und zu derselben zweitens die Funktion
\begin{equation}\label{salzd}
\begin{array}{c}
{\rm T}\langle\varepsilon_{x}\rangle:\mathbb{R}^{+}\times\mathbb{B}_{\varepsilon_{x}(\delta)}(x) \to\mathbb{R}^{+},\\
(\delta,y)\mapsto{\rm T}\langle\varepsilon_{x}\rangle(\delta,y)\ ,
\end{array}
\end{equation}
deren jeweiliger 
Wert f"ur alle $y\in \mathbb{B}_{\varepsilon_{x}(\delta)}(x)$ das Infimum
\begin{equation}\label{wthisb}
{\rm T}\langle\varepsilon_{x}\rangle(\delta,y):=\inf\Bigl\{T\in\mathbb{R}^{+}:t>T\Rightarrow||\Psi(x,t)-\Psi(y,t)|| < \delta\Bigr\}
\end{equation}
sei. F"ur alle $(z,\varepsilon)\in \chi\times\mathbb{R}^{+}$ sei dabei w"ahrend dieses Beweises mit 
$\mathbb{B}_{\varepsilon}(z)$
die offene relative Kugel
$$\mathbb{B}_{\varepsilon}(z):=\{\tilde{z}\in\chi:||z-\tilde{z}||<\varepsilon\}$$
bezeichnet. Die
Funktion ${\rm T}\langle\varepsilon_{x}\rangle$ ist eine atque-fecit-saltus - Funktion\index{atque-fecit-saltus - Funktion} und mit deren
Hilfe zu argumentieren, erfordert Vorsichtigkeit, wie wir bereits in den Vorbemerkungen zu dem
Satz 2.2.2 herausstellten. 
\newline
Die Argumentation mittels der atque-fecit-saltus - Funktion ${\rm T}\langle\varepsilon_{x}\rangle$
halten wir daher gerne kurz und dies ist uns dadurch erm"oglicht,
dass $y$ in dem selben Zimmer ist, in dem auch $x$ ist.
Per constructionem 
schneidet die Menge
$$\Psi(y,[{\rm T}\langle\varepsilon_{x}\rangle(\delta_{1},y),\infty[)$$
auch die Kugel
$$\mathbb{B}_{\varepsilon_{x}(\delta_{2})}(x)$$ 
f"ur die positive reelle Zahl 
$$\delta_{2}=\frac{1}{2}\delta_{1}\ .$$
Also muss es auch die positive reelle Zahl $\Theta(y,\delta_{2})$
geben, f"ur die
f"ur alle $t\in[\Theta(y,\delta_{2}),\infty[$
die Ungleichung
$$||\Psi(x,t)-\Psi(y,t)||\ <\ \delta_{2}$$
erf"ullt sein muss. Wir sehen, dass wir die gleiche Schlussweise immer wieder
wiederholen k"onnen: Die 
Menge
$\Psi(y,[{\rm T}\langle\varepsilon_{x}\rangle(\delta_{j},y),\infty[)$
schneidet
die Kugel
$$\mathbb{B}_{\varepsilon_{x}(\delta_{j+1})}(x)$$ 
f"ur die positive reelle Zahl 
$$\delta_{j+1}:=\frac{1}{2}\delta_{j}$$
f"ur alle $j\in\mathbb{N}$.
Demnach fliesst die Kugel
$\mathbb{B}_{\varepsilon_{x}(\delta_{1})}(x)$ in dem Sinn
in die 
Kugel
$\mathbb{B}_{\varepsilon_{x}(\delta_{j})}(x)$
f"ur alle 
$j\in\mathbb{N}$,
dass f"ur alle 
$j\in\mathbb{N}$ und $z\in\mathbb{B}_{\varepsilon_{x}(\delta_{1})}(x)$
die Aussage
$$\exists t\in\mathbb{R}^{+}:\ \Psi(z,[t,\infty[)\ \subset\
\Psi(\mathbb{B}_{\varepsilon_{x}(\delta_{j})}(x),[0,\infty[)$$
gilt. Die gleiche Argumentationsweise k"onnen wir genausogut f"ur die 
entwicklungsinvertierte, mit Hilfe der Festlegung (\ref{ausfxx}) formulierte
Flussfunktion\index{entwicklungsinvertierte Flussfunktion}
$$\Psi^{-}:=\Psi(\mathbf{P}_{1},-\mathbf{P}_{2})$$
wie f"ur die urspr"ungliche Flussfunktion $\Psi$ praktizieren:
Jeder der Punkte $z\in\mathbb{B}_{\varepsilon_{x}(\delta_{1})}(x)$
muss daher beliebig nahe bei $x$ liegen, was aber mit 
der positiv definiten Funktion $\varepsilon_{x}$ 
gem"ass (\ref{salzc})
per constructionem 
nicht vereinbar ist.  
\newline
Daher kann die Aussage (\ref{salzb}) nicht gelten, wenn
$\chi\in[[\Psi]]$ ein nicht triviales Zimmer ist.
\newline
Der Rest ist nun einfach: F"ur jedes 
nicht triviale Zimmer $\chi\in[[\Psi]]$ gilt also
\begin{equation}\label{sazb}
\begin{array}{c}
\forall x\in \chi\ \exists\ \delta\in\mathbb{R}^{+}\forall\varepsilon\in\mathbb{R}^{+}
\ \exists\ y\in \mathbb{B}_{\varepsilon}(x)\ \forall T\in\mathbb{R}^{+}\ \exists\ t\in ]T,\infty[:\\
||\Psi(x,t)-\Psi(y,t)||\ >\ \delta\ . 
\end{array}
\end{equation} 
F"ur jeden Zustand $x$ des nicht trivialen Zimmers $\chi$
und
f"ur jede positive reelle
Zahl $\tau$ gibt es eine weitere positive reelle
Zahl $\delta\in\mathbb{R}^{+}$ von der Art, dass f"ur jedes $\varepsilon\in\mathbb{R}^{+}$
ein beliebig nahe bei $x$ liegender Zustand $y\in \mathbb{B}_{\varepsilon}(x)$ existiert,
f"ur den 
f"ur
jede nat"urliche Zahl $m$ 
eine streng monoton wachsende Folge 
$$\{t(j,m,\tau)\}_{j\in \mathbb{N}}\in [m\tau,\infty[^{\mathbb{N}}$$
existiert,
f"ur die f"ur alle $j\in \mathbb{N}$
$$||\Psi(x,t(j,m,\tau))-\Psi(y,t(j,m,\tau))||\ >\ \delta$$
gilt. 
Da $\chi$ kompakt und $\Psi$ stetig ist,
gibt es eine gr"osste Zahl $q(\chi,\delta/4)\in\mathbb{R}^{+}$,
f"ur die f"ur alle $(z,t)\in\chi\times \mathbb{R}$ die Ungleichung 
$$\sup||\Psi(z,t+[-q(\chi,\delta/4),q(\chi,\delta/4)])-\Psi(z,t+[-q(\chi,\delta/4),q(\chi,\delta/4)])||\ \leq\ \delta/4$$
wahr ist, sodass f"ur alle 
$$t\in \bigcup\Bigl\{t(j,m,\tau)+[-q(\chi,\delta/4),q(\chi\delta/4)]:j\in \mathbb{N}\Bigr\}$$
die Absch"atzung
$$||\Psi(x,t)-\Psi(y,t)||\ >\ \delta/4$$
gilt. 
Wenn $\tau\leq q(\chi,\delta/4)$ ist, dann gilt also die Aussage
\begin{displaymath}
\begin{array}{c}
\exists\delta/4\in\mathbb{R}^{+}\forall\varepsilon\in\mathbb{R}^{+}
\exists y\in \mathbb{B}_{\varepsilon}(x),\ z\in\mathbb{N}:\\
||\Psi^{\tau})^{z}(x)-\Psi^{\tau})^{z}(y)||\ >\ \delta/4\ .
\end{array}
\end{displaymath} 
{\bf q.e.d.}
\newline
\newline
Es folgt damit unmittelbar das 
\newline
\newline
{\bf Korolllar 2.2.3: Trivialit"atssatz\index{Trivialit\"atssatz}}\newline
{\em Alle Elemente jeder insensitiven Heine-Descartessche Kollektivierung $[\Psi]$ sind Zyklen.} 
\newline
\newline
{\bf Beweis:}
\newline
Die zustandsweise Sensitivit"at jeder Flussfunktion $\Psi$, f"ur die
$$\Delta_{\Psi}(\mathbf{P}_{2}\Psi)\not=\{-\infty\}$$
gilt, widerspricht deren
Asensitivit"at in dem Sinn, dass gem"ass (\ref{knurrze}) f"ur alle $t\in \mathbb{R}$
$$(\Psi^{t})^{-1}\mathbf{cl}\Psi^{t}=\mbox{id}$$
ist.
{\bf q.e.d.}
\newline
\newline
Wir f"ugen dem Satz 2.2.2 die folgende 
Vermutung hinzu, die wir die Vermutung maximaler Sensitivit"at\index{Vermutung maximaler Sensitivit"at} nennen: 
\newline
\newline
{\em Dar"uber hinaus k"onnte 
der Aufl"osungsfeldwert jedes Zimmers $\chi$ als dessen jeweilige Aufl"osungskonstante 
$\Delta(\chi)$
in dem Sinn existieren, dass 
das 
Aufl"osungsfeld\index{Aufl\"osungsfeld einer Wellenfunktion} 
der auf das jeweilige Zimmer $\chi$ restringierten Wellenfunktion $\Psi|\chi\times\mathbb{R}$
konstant ist: Es nimmt den singul"aren Wert $-\infty$ f"ur jedes triviale Zimmer an, keiner
dessen Zust"ande
sensitiv ist.
F"ur 
jedes nicht triviale Zimmer $\chi$ hingegen k"onnte als jeweilige 
Aufl"osungskonstante\index{Aufl\"osungskonstante eines Zimmers}
f"ur alle $x\in\chi$
die Zahl
\begin{equation}
\Delta(\chi):=\Delta_{\Psi|\chi\times\mathbb{R}}(x)=||\chi||\in \mathbb{R}^{+}\ ,
\end{equation} 
vorliegen und einfach mit 
dem Durchmesser des betrachteten nicht trivialen Zimmers $\chi$ "ubereinstimmen.}
\newline
\newline
Ferner bemerken wir:
Wenn ein Takt $t>Q(\chi,\delta)$ 
gr"osser als die entsprechende Konstante $Q(\chi,\delta)$ 
ist, so gibt es eine nat"urliche Zahl 
$p$, f"ur die $t/p\leq Q(\chi,\delta)$ gilt. 
Die Aussage (\ref{salz})
des Satzes 2.2.2 gilt also auch f"ur den Takt $\tau>Q(\chi,\delta)$, wenn
die 
f"ur einen jeweiligen Zustand $y\in \mathbb{B}_{\varepsilon}(x)$ und f"ur das Tripel
$$(x,\delta,\varepsilon)\in \chi\times\mathbb{R}^{+}\times\mathbb{R}^{+}$$
spezifizierte 
Menge  
$$\{z\in\mathbb{N}: ||\Psi^{t/p})^{z}(x)-\Psi^{t/p})^{z}(y)||\ >\ \delta\}\ ,$$
deren Kardinalit"at die Kardinalit"at $\mathbf{card}(\mathbb{N})$ ist,
nicht endlich viele Vielfache der nat"urlichen Zahl $p$ hat,
wenn also
$$\mathbf{card}\Bigl(\{z\in\mathbb{N}: ||\Psi^{t/p})^{z}(x)-\Psi^{t/p})^{z}(y)||\ >\ \delta\}\ \cap\ p \mathbb{N}\Bigr)=\mathbf{card}(\mathbb{N})$$
ist.
Dann ist $\chi$ ein sensitiver Attraktor auch f"ur den Takt $\tau>Q(\chi,\delta)$.
\newline
Beispiele f"ur 
Fixmengen jeweiliger Autobolismen, f"ur Attraktoren und f"ur sensitive Attraktoren sind bei
\cite{male}, \cite{wirr}, \cite{garr}, \cite{ecke} oder bei \cite{walt} zu finden.
\newline Das f"ur uns sehr Bemerkenswerte ist also,
dass Zimmer,
diese elementaren Attraktoren, den zu einer trajektoriellen Partition geh"orenden Zustandsraum partionieren,
wenn derselbe kompakt ist:
\newline\newline
Der Zustandsraum einer trajektoriellen Partition ist in seine elementaren Attraktoren partioniert -
die im allgemeinen sensitive
Attraktoren sind.
\newline\newline
{\footnotesize Der folgende Hinweis irritiert vielleicht. Er soll aber neugierig machen und er will dazu
anregen, weiterzufragen und weiterzulesen:
Wir isolieren aus dem 8. Lehrsatz aus Leibnizens Monadologie\index{Monadologie}\index{Leibniz, Gottfried Willhelm}\newline\newline
\glqq [...]{\em Wenn nun die Monaden\index{Monade} ohne Qualit"aten w"aren, so w"urden sie nicht voneinander zu unterscheiden sein;
denn quantitative Unterschiede
gibt es bei ihnen ja ohnehin nicht. Folglich w"urde -- unter der Voraussetzung,
dass alles voll ist -- jeder Ort bei der Bewegung immer nur das wieder ersetzt erhalten, was er soeben schon gehabt hatte, 
und der Zustand der Dinge w"urde vom anderen ununterscheidbar sein.}\grqq  (siehe \cite{mona})
\newline\newline
mit der Konotation des Heraklitschen 
\begin{center}
\glqq $\pi\alpha\nu\tau\alpha\ \rho\epsilon\iota$\grqq
\end{center} 
die Wendung
\begin{center} 
\glqq {\em dass alles voll ist}\grqq 
\end{center}
und denken dabei nicht nur an das Liouville-Theorem:
F"ur komanent-stetige Kollektive ist der kontinuierliche Determinismus  
so beschaffen, \glqq dass\grqq im Zustandsraum  \glqq alles voll ist\grqq mit - sensitiven Attraktoren -- und zwar im Zustandsraum!\newline
Die Leibnizsche 
Monadologie m"usste so manchem Exegeten
die Schwierigkeit bereitet haben,
dass Leibniz 
seine Beschreibung der Monaden
einerseits
sehr unmissverst"andlich  
in einem Raum durchaus mit geometrischen Begriffen formuliert,
der aber offensichtlich nicht der Konfigurationsraum sein kann.
Beispielsweise ist -- wie wir hier lesen -- in dem Raum, in dem  
Leibniz die Monaden ansiedelt,
alles voller Monaden, die allerdings zugleich vollkommen hermetische Gebilde sind und gleichsam
\glqq [...]{\em keine Fenster} [...]{\em haben}\grqq.
Als den Raum, in dem alles voll ist, denkt sich Leibniz schwerlich den Konfigurationsraum:
Evangelista Torricelli wies bereits 1643 nach, dass Luft Gewicht hat. 
Otto von-Guericke 
hatte im Sommer 1657 sein ber"uhmtes Halbkugelexperiment durchgef"uhrt:
Leibniz muss also wissen, dass Vakuum im "ublichen Sinne hergestellt werden kann.
\newline
Wenn also die Monaden wesentlich auf der Grundlage einer Raumstruktur
beschrieben werden sollen, so muss doch gekl"art werden,
innerhalb welchen Raumes die Monaden liegen.
Die zahlreichen Interpretationen der Monadologie
stellen sich aber dieser Frage nicht.
\newline
Alle Monaden sind "uberdies einmalig. Sie sind auch nicht als
\glqq Dinge\grqq  gedacht, die in der physikalischen Raum-Zeit liegen. Darin,
nicht in der Raum-Zeit zu liegen, sind sie
so wie unsere Zimmer,
die im Phasenraum liegen.\newline
Der Gedanke, die Monaden als eine Vision sensitiver Attraktoren zu identifizieren,
k"onnte sich einerseits als der Schl"ussel zur Monadologie erweisen,
die vielfach als dunkel empfunden wurde. Haben wir hier eine 
einzigartige interdisziplin"are exegetische Chance? Voltaire "ausserte "uber Leibnizens
Monadologie sogar, dass sie den \glqq {\em gesunden Menschenverstand}\grqq  [...] \glqq {\em beleidige}\grqq  und Voltaire
verlieh damit seinem Befremden durch die Monadologie aufrichtig Ausdruck;
anders als viele andere, die von der Beschw"orung des Dunklen leben.
Voltaire unterscheidet sich hier auch
beispielsweise von Kant, der Leibniz in der
Kritik der reinen Vernunft \cite{kant} recht fraglos kritisiert.}

\section{Zweiter Schritt: Maximale Quasiergodik und kinezentrische Felder}\label{secap}
Wir haben nun eine der beiden am 
Schluss des ersten Abschnittes exponierten Aufgaben erf"ullt. Wir haben n"amlich
die "Aquivalenz (\ref{iddito}) 
bewiesen, indem wir den Satz von der Existenz der Zimmer\index{Satz von der Existenz der Zimmer} 2.1.2 zeigten.
Die andere, verbleibende der beiden exponierten Aufgaben, die nach den Ausf"uhrungen des ersten Abschnitts die Essenz der
Quasiergodenfrage bilden,
ist die Aufgabe, die
Gleichheit (\ref{eines}) 
zu beweisen. 
Wir m"ussen 
f"ur jede jeweilige $\mathcal{C}^{1}$-Flussfunktion binnen eines endlichdimensional-reellen Kompaktums  
die st"uckweise glatten minimalen invarianten Mannigfaltigkeiten der Menge $\mathcal{M}^{ 1/2} ([\Psi])$ mit den abgeschlossen
H"ullen der Trajektorien aus der Menge $[\Psi]$ identifizieren.
Wenn wir diese Aufgabe erf"ullen, erreichen wir unser Ziel.
Unser Satz von der Existenz von Zimmern liegt dabei auf dem restlichen Weg zu diesem unserem Ziel.
\newline
Wenn die Zimmer aus der Menge $[[\Psi]]$ minimale invariante Mannigfaltigkeiten 
der Menge $\mathcal{M}^{ 0} ([\Psi])$
schneiden,
so liegen die Zimmer offensichtlich immer ganz in jenen: Dass
die 
Inklusionen
\begin{equation}
\mathcal{M}^{0}  ([\Psi])\ \subset\ \Bigl\{\ \bigcup \Xi:\ \Xi\ \subset\ [[\Psi]]\Bigr\}
\end{equation}
und
\begin{equation}\label{aklaeib}
\mathcal{M}^{ 1} ([\Psi])\ \subset\ \mathcal{M}^{ 1/2} ([\Psi])\ \subset\ \mathcal{M}^{0}([\Psi])
\end{equation}
und, dass die Implikation
\begin{equation}
\chi\in[[\Psi]]\ \Rightarrow\ \downarrow\Bigl\{\mu\in\mathcal{M}^{0}([\Psi]):\mu\cap\chi\not=\emptyset\Bigr\} \supset\ \chi
\end{equation}
und deren zu ihr gegenl"aufiges Pendant, die Implikation
\begin{equation}\label{aklaeia}
\mu\in\mathcal{M}^{0}([\Psi])\ \Rightarrow\ \Bigl\{\chi\in[[\Psi]]:\mu\cap\chi\not=\emptyset\Bigr\}\in\mathbf{part}(\mu)
\end{equation}
f"ur jede trajektorielle $\mathcal{C}^{1}$-Partition eines Kompaktums $\bigcup[\Psi]$ eines endlichdimensionalen reellen Raumes
gelten, ist nun nach dem Satz von der Existenz von Zimmern 2.1.2 evident. Der Operator 
$\downarrow$ bezeichnet dabei wieder das Element jeder jeweiligen einelementigen Menge.
Wir betonen aber sehr, dass,
was wir zeigen wollen und was uns bedeutend zu sein scheint, eben nicht so offensichtlich ist, n"amlich, dass die nat"urliche Partition  
\begin{equation}\label{aklei}
[[\Psi]] = \mathcal{M}^{ 1/2} ([\Psi])
\end{equation}
ist. Dass diese Gleichheit gilt, impliziert hier n"amlich letztlich die Identit"at 
\begin{equation}\label{aklaei}
[[\Psi]] = \mathcal{M}^{0} ([\Psi]) = \mathcal{M}^{ 1/2} ([\Psi])\ ,
\end{equation}
wie wir dies sogleich darlegen werden. In der Identit"at (\ref{aklei}) "offnet sich uns der Blick auf den folgenden Sachverhalt:
\newline \newline 
{\em Die Zimmer sind so gross, wie sie nach unseren Ausf"uhrungen vor der Formulierung des originalen Quasiergodensatzes} 1.2.2
{\em trivialermassen, d.h.,
von vornherein, h"ochstens sein k"onnen.}
\newline\newline 
Deshalb nennen wir
diese Eigenschaft von Zimmern die 
maximale Quasiergodik\index{maximale Quasiergodik} derselben, um so einen Extremal- und Aussch"opfungsaspekt der Identifizierung
(\ref{aklei}), die zugleich das Hinf"alligwerden einer Unterscheidung ist, hervorheben zu k"onnen. Aussagenlogisch gesehen geht hier 
eine Unterscheidung im
Hinf"alligwerden einer Differenzierung unter:\newline
\newline
{\em Maximale Quasiergodik ist die Identit"at von Zimmern mit minimalen invarianten $\mathcal{C}^{0}$-Mannigfaltigkeiten.}
\newline 
\newline
Wie kommt es denn nun zu der Identit"at der st"uckweise glatten minimalen invarianten Mannigfaltigkeiten\footnote{Dieselben nennen wir ja schlicht
minimale invariante Mannigfaltigkeiten, was durch die Indifferenz (\ref{aklaei}) berechtigt ist.}
mit den jeweiligen minimalen invarianten $\mathcal{C}^{0}$-Mannigfaltigkeiten? 
Wenn wir zeigen, dass f"ur jede trajektorielle $\mathcal{C}^{1}$-Partition eines Kompaktums $[\Psi]$ die Identit"at
(\ref{aklei})
gilt, so zeigen wir nicht nur die maximale Quasiergodik, sondern wir zeigen dann auch, dass minimale invariante Mannigfaltigkeiten und minimale invariante $\mathcal{C}^{0}$-Mannigfaltigkeiten f"ur die 
trajektoriellen $\mathcal{C}^{1}$-Partitionen eines Kompaktums eines endlichdimensionalen und reellen Raumes
"ubereinstimmen. Dies lehrt folgendes Reduktionslemma, das innerhalb unserer Argumentation die Position eines Korollares 
der zu beweisenden Identit"at (\ref{aklei})
einnimmt. Wir bezeichnen es dennoch als ein Lemma und heben es als einen eigenen Sachverhalt heraus, obwohl sein Beweis trivial ist.
Das Reduktionslemma hilft n"amlich dabei, die maximale Quasiergodik zu beweisen,
indem es erlaubt, von der G"ultigkeit der Identit"at (\ref{aklei})
auf die G"ultigkeit der Identit"at (\ref{aklaei})
zu schliessen, wobei
wir die Gleichheit (\ref{aklei}) anschliessend erst noch unabh"angig von dem Reduktionslemma zu zeigen haben.
Als ein Lemma sehen wir das Reduktionslemma auch deshalb an, weil es ganz allgemein
eine Hilfe dabei ist, sich innerhalb von trajektoriellen $\mathcal{C}^{1}$-Partitionen eines Kompaktums
zu orientieren.
\newline
\newline
{\bf 2.3.1 Reduktionslemma:}\index{Reduktionslemma}
\newline
{\em Ist $[\Psi]$ eine trajektorielle $\mathcal{C}^{1}$-Partition eines Kompaktums eines endlichdimensionalen reellen Raumes, 
so gilt die "Aquivalenz}
\begin{equation}
[[\Psi]] = \mathcal{M}^{0}([\Psi])\   \Leftrightarrow\ [[\Psi]] = 
\mathcal{M}^{1/2}([\Psi])\ .
\end{equation}
\newline
{\bf Beweis:}\newline
Wir gehen von der Richtigkeit der Identit"at (\ref{aklei})
aus und brauchen daher nur die Implikation
$$[[\Psi]] = \mathcal{M}^{1/2}([\Psi])\Rightarrow[[\Psi]] = \mathcal{M}^{0}([\Psi])$$
zu zeigen.
Gem"ass der Gleichheit (\ref{aklaeia}) gibt es f"ur jede minimale invariante 
$\mathcal{C}^{0}$-Mannigfaltigkeit $\mu\in\mathcal{M}^{0}([\Psi])$
die Partition ${\rm P}(\mu)\subset[[\Psi]]$ in Zimmer 
$\chi\in[[\Psi]]=\mathcal{M}^{1/2}([\Psi])$, die aber 
nach unserer Vorraussetzung, dass die Identit"at (\ref{aklei}) gilt,
jeweils selber minimale invariante $\mathcal{C}^{1/2}$-Mannigfaltigkeiten
sind; und daher sind dieselben gem"ass (\ref{aklaeib})
jeweils alle auch minimale invariante $\mathcal{C}^{0}$-Mannigfaltigkeiten. Dabei ist
die Implikation 
$$a\in \mathcal{M}^{0}([\Psi])\ \land\ b\in \mathcal{M}^{1/2}([\Psi]) \land\
a\cap b\not=\emptyset$$
$$\Rightarrow a\subset b$$
trivial, sodass 
die Partition ${\rm P}(\mu)$ nur ein einziges Element haben kann. 
\newline
{\bf q.e.d.}
\newline
\newline
Wenn wir nun 
f"ur jede jeweilige $\mathcal{C}^{1}$-Flussfunktion $\Psi$ binnen einer 
abgeschlossenen Teilmenge $\zeta$ eines
reellen Raumes $\mathbb{R}^{n}$ der endlichen Dimension $n\in\mathbb{N}$
eine Abbildung $f$ 
des Zustandsraumes $\zeta$ in einen
endlichdimensionalen und reellen Raum
$\mathbb{R}^{\nu}$ f"ur $\nu\in\mathbb{N}$ genau dann identitiv bez"uglich der 
zu $\Psi$
geh"orenden trajektoriellen Partition $[\Psi]$ nennen,\index{identitive Abbildung bez\"uglioh 
einer trajektoriellen Partition}\index{Identitivit\"at einer Abbildung} 
wenn f"ur alle
$x_{ 1}, x_{ 2} \in \zeta$ die Implikation
\begin{equation}
f ( x_{ 1} ) = f ( x_{ 2} ) \Rightarrow\ \mathbf{cl} (\Psi(x_{ 1} ,\mathbb{R}))\ 
=\ \mathbf{cl} (\Psi(x_{ 2},\mathbb{R}))
\end{equation}
wahr ist, so l"asst sich die Konsequenz des Reduktionslemmas so formulieren: Um die maximale Quasiergodik der trajektoriellen
$\mathcal{C}^{1}$-Partitionen eines Kompaktums eines endlichdimensionalen und reellen Raumes
zu zeigen, gen"ugt es,
vorzuf"uhren, dass sich f"ur jede solche trajektorielle Partition $[\Psi]$
eine global
Lipschitz-stetige und identitive Abbildung finden l"asst, die auf deren Zustandsraum $\zeta$ definiert ist.
\newline  
\newline
Es ist kein Exotikum, das die Eigenschaft der Identitivit"at bez"uglich einer jeweiligen trajektoriellen
$\mathcal{C}^{1}$-Partition $[\Psi]$ eines Kompaktums $\zeta$ eines endlichdimensionalen und reellen Raumes
hat:
\newline
\newline
Es sei f"ur jede
$\mathcal{C}^{1,-1/2}$-Flussfunktion $\Phi$ binnen eines 
beschr"ankten Zustandsraumes eines endlichdimensionalen und reellen Raumes
die Abbildung 
\begin{displaymath}
\Omega_{\Phi}:\ \zeta\ \to\ \zeta,\ x \mapsto \Omega_{\Phi} (x)
\end{displaymath} 
\begin{equation}
:=\ \lim_{T\to\infty}\ \frac{1}{2 T}\ \int_{-T}^{T}\ \Phi(x,t)\ dt
\end{equation}
erkl"art.
Wir wollen exakt diese 
Abbildung $\Omega_{\Phi}$
das kinezentrische Feld
bez"uglich der Flussfunktion $\Phi$ nennen.\index{kinezentrisches Feld}
Ferner nennen wir exakt
alle kinezentrischen Felder $\Omega_{\Phi}$
bez"uglich
einer der $\mathcal{C}^{1,-1/2}$-Flussfunktionen $\Phi$,
f"ur die $[\Phi]=[\Psi]$ gilt,
die kinezentrischen Felder der trajektoriellen Partition $[\Psi]$.
Ein kinezentrisches Feld $\Omega_{\Phi}$ einer trajektoriellen Partition ist offensichtlich nicht nur von seiner trajektoriellen 
Partition $[\Psi]$, 
sondern auch von der Wahl der 
das kinezentrische Feld $\Omega_{\Phi}$ indizierenden Flussfunktion abh"angig.
F"ur kinezentrische Felder gilt folgende
\newline
\newline
{\bf 2.3.2 Bemerkung:}\newline {\bf Die Eigenschaften kinezentrischer Felder von trajektoriellen Partitionen:}
\index{Eigenschaften kinezentrischer Felder von trajektoriellen Partitionen}
\newline
{\em Alle kinezentrischen Felder jeder trajektoriellen $\mathcal{C}^{1}$-Partition $[\Psi]$ 
eines Kompaktums eines endlichdimensionalen und reellen Raumes
existieren. Diese kinezentrischen Felder sind dabei sowohl identitiv als auch Lipschitz-stetig.}
\newline
\newline
Der folgende Beweis der Existenz der jeweiligen kinezentrischen Felder ist leicht generalisierbar
und dessen erster Teil, der Beweis, dass kinezentrische Felder existieren, kann auch als L"osung einer "Ubungsaufgabe der Analysis 1 eines Analysiskurses aufgefasst werden.
\newline
\newline {\bf Beweis:}\newline
Sei $\alpha \in \mathcal{C}([0, \infty), [1, 2])$ stetig. F"ur alle $t \in ]0,\infty)$ existiert die Zahl
$$m (t)\ \alpha\ :=\ \frac{1}{t}\ \int_{0}^{t}\ \alpha (s)\ ds\ \in\ [1,\ 2]\ .$$
$T^{ \infty}$ sei die Menge aller streng monotoner und divergenter positiv-reeller Folgen.
F"ur $\{t (j)\}_{ j \in\mathbb{N}} \in\ T^{ \infty}$ existiert eine Teilfolge $\{t (j (k)\}_{ k \in \mathbb{N}}$, f"ur die
$\{m\ (t (j (k)))\}_{ k \in \mathbb{N}}$ gegen $a \in [1, 2] $ konvergiert.
Es sei $b \in [1, 2]$ ein H"aufungspunkt der Folge $\{m\ ( b (j))\ \alpha\}_{ j\in \mathbb{N}}$ f"ur
eine Folge $\{b (j)\}_{ j \in\ \mathbb{N}}\ \in\ T^{ \infty} $,
sodass die Folge $\{b (j (l))\}_{\ l \in \mathbb{N}}$ existiert, f"ur die 
$$\lim_{j\to\infty}\ m\ ( b ( j (l)))\ \alpha = b$$
ist.
Es gibt eine Teilfolge von $\{t (j (k))\}_{ k \in \mathbb{N}}$, n"amlich $\{t (j ( k (l)))\}_{ l \in \mathbb{N}}$, von der Art, dass
f"ur alle $l \in \mathbb{N}$
$$t (j ( k (l)))\ <\ b (j (l))\  <\ t (j ( k (l+1)))$$ ist.
$m\ (\mbox{id})\  \alpha$ ist dabei eine $\mathcal{C}^{ 1}$-Funktion. Die Werte ihrer Ableitung
$$\frac{d}{dt}\ m\ (t)\ \alpha = \frac{1}{t^{ 2}}\ (t\ \alpha\ (t)\ -\ m\ (t)\ \alpha)$$
existieren
f"ur alle $t\in [0,\infty)$.
Wenn $a$ und $b$ verschieden w"aren, so k"ame
$m\ (\mbox{id})\ \alpha$ f"ur beliebig grosse Argumente sowohl der Zahl $a$
als auch der Zahl $b$ beliebig nahe,
sodass eine Folge $\{a (j)\}_{ j \in\ \mathbb{N}}\ \in\ T^{ \infty}$ existierte von der Art, dass
$$a (j)\ \alpha\ (a (j))\ -\ m\ (a (j))\ \alpha = 0$$
w"are; sodass
f"ur alle $j \in \mathbb{N}$ die Ungleichung
$$2\ \geq\ m\ (a (j))\ \alpha\ \geq\ \alpha\ (a (j))\ \geq\ a (j)$$
g"ultig w"are.
Also existiert auch f"ur zwei Zahlen $v, w \in \mathbb{R}$ mit $v < w$ und f"ur eine Funktion $f\in \mathcal{C}([0,\infty), [v, w])$
der Grenzwert
$$\lim_{t\to\infty}\ \frac{1}{t}\ \int_{0}^{t}\ f (s)\ ds\ ,$$
denn es ist
$$(w-v)\ \mathcal{C}([0,\infty), [1, 2])\ +\ (2v-w)\  =\  \mathcal{C}([0,\infty), [v, w])$$
und daher existieren auch alle kinezentrischen Felder $\Omega_{\Phi}$ f"ur alle 
$\mathcal{C}^{1,-1/2}$-Flussfunktionen $\Phi$,
f"ur die $[\Phi]=[\Psi]$ gilt. 
Soweit der Beweis, dass kinezentrische Felder existieren.
\newline\newline
Und nicht nur diese kinezentrischen Felder $\Omega_{\Phi}$ existieren:
F"ur jede auf dem Zustandsraum $\zeta$ st"uckweise stetige, reelle Funktion $f$ und f"ur das Lebesgue-Mass $\lambda$
auf dem Zimmer $\downarrow [[\Psi]]_{\{x\}}$, in dem $x$ ist, existiert also f"ur jeden Punkt $x\in\zeta$
die Flussfunktions-Invariante
\begin{equation}\label{Evinv}
\lim_{T\to\infty}\ \frac{1}{2 T}\ \int_{-T}^{T}\ f(\Psi(x,t))\frac{f(\Psi(x,t))}{||\partial_{2}\Psi(x,t)||}\ dt = 
\int_{\downarrow [[\Psi]]_{\{x\}}}\ f(x)\ d\lambda(x)\ 
\end{equation}
$$=\lim_{T\to\infty}\ \frac{1}{2 T}\ \int_{-T}^{T}\ f(\Phi(x,t))\frac{f(\Phi(x,t))}{||\partial_{2}\Phi(x,t)||}\ dt\ ,$$
falls
die Ableitung $\partial_{2}\Phi(x,t)$
und der Bruch
$$\frac{f(\Phi(x,t))}{||\partial_{2}\Phi(x,t)||}$$
f"ur alle $t\in\mathbb{R}$ existieren. 
Damit die Invarianz (\ref{Evinv}) gilt,
muss $\Phi$ also eine  parziell nach der zweiten Ver"anderlichen differenzierbare Flussfunktion $\Phi$ sein,
f"ur die $[\Phi]=[\Psi]$ gilt.
\newline
Die Identit"at (\ref{Evinv}) ist die ber"uhmte 
Boltzmannsche 
Scharmittel-Zeitmittel-Identifizierung.
\index{Scharmittel-Zeitmittel-Identifizierung}
Wir erkennen schon jetzt, was im Fall allgemeinerer Gegebenheiten als der Gegebenheit
von trajektoriellen $\mathcal{C}^{1}$-Partitionen $[\Psi]$ 
eines Kompaktums eines endlichdimensionalen und reellen Raumes
gilt:
Ein Kriterium daf"ur,
dass es die entsprechend verallgemeinerten
Invarianten (\ref{Evinv}) gibt, ist neben der Lebesgue-Messbarkeit von $f$
das Kriterium,
dass "uberhaupt Lebesgue-messbare Zimmer
als die abgeschlossenen H"ullen in ihnen liegender Trajektorien vorliegen
und, dass die in ihnen liegenden Trajektorien 
st"uckweise glatt sind. 
Dabei wissen wir schon, dass Zimmer beschr"ankter komanent-immanenter Trajektorienkollektive existieren
und Kompakta sind.
\newline
Wir verweisen darauf,
die Aussage (\ref{Evinv})
mit dem sogenannten Bowen-Ruelle-Theorem 
\index{Bowen-Ruelle-Theorem} zu vergleichen,
das z.B. in \cite{ecke} zu finden ist.
\newline\newline
F"ur alle $ x\in \zeta$ und $t \in \mathbb{R}$ ist
\begin{equation}
\Omega_{\Phi} (\Psi(x,0)) = \Omega_{\Phi} (\Psi(x,t))\ ,
\end{equation}
weil f"ur alle $T\in\mathbb{R}^{+}$ die Ungleichung
$$||\ \frac{1}{2\ T}\ \int_{-T}^{T} \Psi(x,s) \ ds\ -\ \frac{1}{2\ T}\ \int_{-T}^{T}\ \Psi(x,s+t)\ ds\ ||$$
$$\leq\ \frac{1}{2\ T}\ \max\{||\Psi(x,s)||:s \in \mathbb{R}\}$$
gilt: Wenn $\Omega_{\Phi}$ stetig ist,
so ist $\Omega_{\Phi}$ auch identitiv.
\newline
\newline
Zeigen wir also die globale Lipschitz-Stetigkeit
des kinezentrischen Feldes 
$\Omega_{\Phi}$ auf dem gesamten Zustandsraum $\zeta$,
so ist damit auch die Stetigkeit und Identitivit"at von $\Omega_{\Phi}$ und damit die maximale 
Quasiergodik gezeigt.
\newline
Um die Lipschitz-Stetigkeit von $\Omega_{\Phi}$ auf $\zeta$ zu beweisen,
betrachten wir den Differenzenquotienten
\begin{equation}
Q\ (x, y)\ :=\ \frac{||\Omega_{\Phi} (x) - \Omega_{\Phi} (y)||}{||x\ -\ y||}
\end{equation}
f"ur $x\neq y,\ x,\ y\in\ \zeta$. F"ur $x \in\ \Psi(y,\mathbb{R}) \setminus \{\ y\}$ sei gesetztermassen $Q\ (x,y):=0$. Sei  
$x \not\in\ \Psi(y,\mathbb{R})$:
Da die Absch"atzung
$$||\Omega_{\Phi} (x)\ -\Omega_{\Phi} (y) ||\ 
\leq\ \max\Bigl\{||\Psi(x,t)- \Psi(y,t)||: t\in \mathbb{R}\Bigr\}\ =:\ A(x, y)$$
gilt und die Komanenz von $[\Psi]$ gegeben ist und ein Flussexponent $\kappa(\Psi)$
gem"ass (\ref{vovo}) existiert, f"ur den f"ur alle $t\in \mathbb{R}$ die Beschr"ankung
\begin{equation}
||\Psi(x,t)-\ \Psi(y,t)||\ \leq\ ||x-y||\ e^{\kappa (\Psi) | t |}
\end{equation}
gilt, ist die Beschr"ankung
$$Q\ (x, y)\ \leq\ \frac{A (x, y)}{||\Psi(x,t) - \Psi(y,t)||}\ e^{ \kappa (\Psi) | t |} $$
f"ur alle $t \in \mathbb{R}$ eingehalten, wobei $A$ die  
von (\ref{vonvon}) her bereits bekannte Metrik und stetige Funktion ist.
Sei
$$\Upsilon:\ (\zeta)^{ 2} \to\ \mathbb{R}^{ +} ,\
(x,y)\ \mapsto\ \Upsilon (x,y)$$
\begin{equation}
:=\ \inf \Bigl\{t\in \mathbb{R}^{ +}: ||\Psi(x,t) - \Psi(y,t)||\ >\ \frac{1}{2} A (x,y)\ \Bigr\}\ .
\end{equation}
Wenn diese Funktion $\Upsilon$ stetig ist,
so nimmt sie auf $\zeta\times\zeta$ ihr Maximum 
\begin{equation}
\Upsilon  (\Gamma):=\max\Bigl\{\Upsilon (x,y): x,y \in \zeta\Bigr\}
\end{equation}
an und es gilt f"ur alle $x,y \in \zeta$
\begin{equation}
Q\ (x, y)\ \leq\ 2\ e^{\ \kappa (\Psi)\ \Upsilon (\Gamma)}\ .
\end{equation}
$\Omega_{\Phi}$ ist dann Lipschitz-stetig und es ist die thematisierte maximale Quasiergodik gegeben.
\newline
\newline
Die Unstetigkeit der Abbildung $\Upsilon$ anzunehmen, f"uhrt aber auf einen Widerspruch:\newline
\newline
Denn dann existierten 
zwei Zust"ande
$v\in \zeta$ und $w\in \zeta$ im Zustandsraum und eine positive Zahl $\alpha\in\mathbb{R}^{+}$ von der Art,
dass f"ur jede positive Zahl $\beta \in \mathbb{R}^{ +}$ die positive Definitheit
\begin{equation}
\max\Bigl\{||w (t) - v(t)||: t\in [0, \Upsilon (w,v) + \alpha ]\Bigr\}\ -\ \frac{1}{2} A (w,v)\ >\ 0
\end{equation}
wahr w"are,
wobei dennoch f"ur alle $\delta \in \mathbb{R}^{ +}$ zwei Zust"ande
$\omega, \nu\in\zeta$ zu finden w"aren, f"ur die
$$||w - \omega||\ <\delta\ > ||\nu - v||$$ 
und zudem
\begin{equation}
\max\Bigl\{||\omega (t)- \nu (t)||: t\in [0, \Upsilon (\omega, \nu)\ +\ \beta ]\Bigr\}\ -\ \frac{1}{2} A (\omega,\nu)\ 
<\ 0
\end{equation}
g"alte,
wo doch die Komanenz der trajektoriellen Partition $[\Psi]$ vorausgesetzt ist und wir schon wissen,
dass die Funktion $A$ auf ihrem Definitionsbereich
$\zeta\times\zeta$ stetig ist.
\newline
{\bf q.e.d.}
\newline
\newline
Wegen dem Satz von der Existenz der Zimmer
2.1.2, dem Reduktionslemma 2.3.1 und schliesslich wegen der nunmehr gezeigten Bemerkung 2.3.2
folgt nun also die Identit"at der Menge der
minimalen invarianten $\mathcal{C}^{1/2}$-Mannigfaltigkeiten mit
der Menge der minimalen invarianten $\mathcal{C}^{1}$-Mannigfaltigkeiten 
einer jeweiligen trajektoriellen
$\mathcal{C}^{1}$-Partition $[\Psi]$ 
eines Kompaktums eines endlichdimensionalen und reellen Raumes und der
\newline
\newline
{\bf 2.3.3 Identit"atssatz:}\index{Identit\"atssatz}\newline
{\em F"ur jede trajektorielle $\mathcal{C}^{1}$-Partition $[\Psi]$ 
eines Kompaktums eines endlichdimensionalen und reellen Raumes gilt}
\begin{equation}\label{iden}
\mathcal{M}^{0} ([\Psi])=
\mathcal{M}^{ 1/2 }([\Psi])=[[\Psi]]\ .
\end{equation}

\chapter{Anhang: kommentierende Erg"anzungen}
\begin{flushright}
Allen Besinnungsvollen: Was war?\\
\end{flushright}
Wir verwiesen bereits hierher auf diesen Anhang. Dies geschah erstmals,
als wir nach der Konstruktion minimaler invarianter Mannigfaltigkeiten einer trajektoriellen Partition behaupteten,
dass
diese minimalen invarianten Mannigfaltigkeiten 
nicht 
einfach die Trajektorien der jeweiligen trajektoriellen Partition sind. Sowohl der Sachverhalt,
dass eine jeweilige minimale invariante Mannigfaltigkeit mit einer Trajektorie "ubereinstimmt
als auch die Negation desselben Sachverhaltes sind topologische Invariante. Daher 
gen"ugt es, 
trajektorielle $\mathcal{C}^{1}$-Partitionen zu betrachten, als deren Elemente es
Trajektorien gibt, die von den glatten minimalen invarianten Mannigfaltigkeiten differieren,
in denen diese Trajektorien
liegen.
Das Kernst"uck einer Erl"auterung der Differenz zwischen 
minimalen invarianten Mannigfaltigkeiten
und Trajektorien bildet demnach die Behauptung des folgenden Satzes:
\newline\newline 
{\bf Satz A.1: Die Nicht-Identit"at von glatten minimalen invarianten Mannigfaltigkeiten
mit glatten Trajektorien }
\newline
{\em Es sei 
$[\Psi]$ eine trajektorielle $\mathcal{C}^{1}$-Partition des Zustandsraumes $\bigcup[\Psi]$, der
ein Kompaktum eines endlichdimensionalen und reellen
Raumes bez"uglich dessen nat"urlicher Topologie sei.
Alle Trajektorien der trajektoriellen Partition $[\Psi]$ sind entweder geschlossen oder quasiperiodisch.
Wenn eine ihrer glatten minimalen invarianten Mannigfaltigkeiten der Menge $\mathcal{M}^{ 1}([\Psi])$ mit einer Trajektorie $\tau$ aus $[\Psi]$ identisch ist,
so ist $\tau$ geschlossen.}
\newline
\newline
{\bf Beweis:}\newline
Es sei $m\in\mathbb{N}$ und
$$\beta=(\beta_{k})_{\ 1\leq k\leq m}\in (\mbox{Inv}^{1}([\Psi]))^{m}$$
ein $m$-Tupel glatter Invarianter der kompakten trajektoriellen Partition 
$[\Psi]$. F"ur $\omega=(\omega_{k})_{\ 1\leq k\leq m}\in\mathbb{R}^{m}$
sei 
$$\beta^{-1}(\{\omega\})\in \mathcal{M}^{ 1}([\Psi])$$ eine 
glatte
minimale invariante Mannigfaltigkeit
dieser kompakten trajektoriellen Partition
$[\Psi]$.
Wenn diese minimale invariante Mannigfaltigkeit
$\beta^{-1}(\{\omega\})$ mit einer Trajektorie $\tau\in [\Psi]$ koinzidiert, d.h., wenn
$\beta^{-1}(\{\omega\})$
so beschaffen ist,
dass es eine Trajektorie $\tau\in [\Psi]$ gibt,
sodass $$\beta^{-1}(\{\omega\})=\tau$$  ist,
so ist 
diese Trajektorie
$\tau=\mathbf{cl}(\tau)$ 
abgeschlossen,
mithin kompakt,
weil die trajektorielle Partition $[\Psi]$ kompakt ist.
\newline
Betrachten wir die 
kompakte trajektorielle Partition $[\Psi]$ mit Trajektorien im $\mathbb{R}^{\nu}$ 
f"ur eine nat"urliche Dimension $\nu$ weiter:
Jede Trajektorie $\tau\in [\Psi]$ hat in jedem ihrer Punkte
beidseitig unendliche L"ange:
F"ur alle $x\in\zeta$ existiert eine Funktion 
$\tilde x\in \mathcal{C}^{ 1}(\mathbb{R},\zeta)$ mit $\{\tilde \Psi(x,\mathbb{R})\}=[\Psi]_{\{x\}}$ und $\tilde \Psi(x,0)=x$.
Jede Folge $\{t_{j}\}_{\ j\in\mathbb{N}}\in \mathbb{R}^\mathbb{N}$ ist so beschaffen,
dass f"ur alle H"aufungspunkte $\tilde \Psi(x,\iota)$
der Folge $\{\tilde \Psi(x,t_{j})\}_{\ j\in\mathbb{N}}$ eine Abbildung $\iota\in \mathbb{N}^\mathbb{N}$
von der Art existiert,
dass dieser H"aufungspunkt $\lim_{j\to\infty}\tilde \Psi(x,t_{\iota(j)})=\tilde \Psi(x,\iota)$ der Grenzwert einer 
durch die Abbildung $\iota$
beschriebenen Teilfolge von $\{\tilde \Psi(x,t_{j})\}_{\ j\in\mathbb{N}}$ ist. Die Trajektorie
$\tilde \Psi(x,\mathbb{R})$ ist also immer entweder geschlossen und periodisch oder quasiperiodisch,
wobei der Grenzwert $\tilde \Psi(x,\iota)$ in einer Teilmenge $\kappa\subset\bigcup[\Psi]$ 
des Zustandsraumes ist; 
letztere ist in dem Sinn $n$-dimensional,
dass es eine lokale Basis von $n\in\mathbb{N}$ Elementen der Menge $B(\tilde \Psi(x,\iota)):=\{b_{1},.. b_{n}\}\subset \mathbb{R}^{\nu}$ gibt,
die linear unabh"angig ist; und diese Basis ist dabei von der Art,
dass f"ur
eine positive reelle Zahl $\varepsilon$ der in $\tilde \Psi(x,\iota)$ zentrierte Spat
$\{\tilde \Psi(x,\iota)\}+\sum_{j=1}^{n}b_{j}[-\varepsilon,\varepsilon]$ in $\kappa$ enthalten ist,
wobei aber noch die folgende Bedingung eingehalten ist:
Es gibt
keine andere linear unabh"angige Menge $\{c_{1},.. c_{n+1}\}$ mit $n+1$ Elementen mit der Eigenschaft,
dass der Spat $\{\tilde \Psi(x,\iota)\}+\sum_{j=1}^{n+1}c_{j}[-\varepsilon,\varepsilon]$ in $\kappa$ enthalten ist.
Dabei ist
$$\Theta_{\iota(j)}^{z}:=\lim_{h\to 0}
\frac{\tilde \Psi(x,t_{\iota(j)})-\delta_{1z}\tilde \Psi(x,\iota)-\delta_{0z}\tilde \Psi(x,t_{\iota(j)}-h)}{||
\tilde \Psi(x,t_{\iota(j)})-\delta_{1z}\tilde \Psi(x,\iota)-\delta_{0z}\tilde \Psi(x,t_{\iota(j)}-h) ||}\ \in\ \mathbb{S}^{n-1}$$
lokalisiert in der zentrierten Einheitssph"are $\mathbb{S}^{n-1}$ des $\mathbb{R}^{n}$,
sodass alle H"aufungspunkte der Folge $\{\Theta_{\iota(j)}^{z}\}_{\ j\in \mathbb{N}}$ in der Einheitssph"are $\mathbb{S}^{n-1}$ liegen.
All dieser H"aufungspunkte Menge bezeichne f"ur beide $z\in\{0,1\}$
der Ausdruck 
$\Theta_{z}^{\star}\subset\mathbb{S}^{n-1}$.
F"ur jeden dieser H"aufungspunkte $s\in\Theta_{z}^{\star}$ gibt es eine Funktion 
$\iota(z,s)\in \iota(\mathbb{N})^{\mathbb{N}}$
von der Art,
dass f"ur beide Indizes $z\in\{0,1\}$ die Teilfolge $\{\Theta_{\iota(z,s)(j)}\}_{j\in \mathbb{N}}$ gegen $s$ konvergiert. 
$\Theta_{0}^{\star}$ hat nur ein einziges Element.
Dieses Element 
ist die Einheitstangente 
in dem Punkt $\tilde \Psi(x,\iota)$
an die Trajektorie, die durch den Punkt $\tilde \Psi(x,\iota)$
verl"auft.  
Sonst w"are das Einheitstangentenfeld in $\tilde \Psi(x,\iota)$ nicht stetig.
Dabei ist die Trajektorie $\tilde \Psi(x,\mathbb{R})$ geschlossen oder quasiperiodisch.
\newline
Nehmen wir an,
die Trajektorie
$\tilde \Psi(x,\mathbb{R})$ sei abgeschlossen, also $\tilde \Psi(x,\mathbb{R})=\mathbf{cl}(\tilde \Psi(x,\mathbb{R}))$:
F"ur jede minimale invariante Mannigfaltigkeit $\beta=(\beta_{k})_{\ 1\leq k\leq m}\in \mathcal{M}^{ 1}([\Psi])$ mit
$\omega=(\omega_{k})_{\ 1 \leq k\leq m}\in\mathbb{R}^{m}$ mit
$\beta^{-1}(\{\omega\})\subset\tilde \Psi(x,\mathbb{R})$ gilt dann nach dem Taylorschen Satz f"ur $s\in\Theta^{\star}_{1}$
f"ur alle $k\in\{1,..m\}$ und $j\in \mathbb{N}$
\begin{displaymath}
\omega_{k}=\beta_{k}\left(\tilde \Psi(x,\iota)+(\tilde x_{\iota(1,s)(j)}-\tilde \Psi(x,\iota))\right)
=\beta_{k}\left(\tilde \Psi(x,t_{\iota(1,s)(j)})\right)
\end{displaymath}
\begin{displaymath}
=\omega_{k}\ +\ \nabla\beta_{k}(\tilde \Psi(x,\iota))^{\top}\left(\tilde \Psi(x,t_{\iota(1,s)(j)})-\tilde \Psi(x,\iota)\right)\ +\
\end{displaymath}
$$o\ \left(||\tilde \Psi(x,t_{\iota(1,s)(j)})-\tilde \Psi(x,\iota)||^{2}\right)$$
\begin{displaymath}
=\omega_{k}\ +\ \nabla\beta_{k}(\tilde \Psi(x,\iota))^{\top}\ s\ ||\tilde \Psi(x,t_{\iota(1,s)(j)})-\tilde \Psi(x,\iota)||\ +\
\end{displaymath}
$$o\ \left(||\tilde \Psi(x,t_{\iota(1,s)(j)})-\tilde \Psi(x,\iota)||^{2}\right)\ +$$
\begin{displaymath}
\nabla\beta_{k}(\tilde \Psi(x,\iota))^{\top}\left(
\frac{\tilde \Psi(x,t_{\iota(1,s)(j)})-\tilde \Psi(x,\iota)}{||\tilde \Psi(x,t_{\iota(1,s)(j)})-\tilde \Psi(x,\iota)||}\ -\ s\right)\ 
||\tilde \Psi(x,t_{\iota(1,s)(j)})-\tilde \Psi(x,\iota)||\ ,
\end{displaymath}
sodass 
\begin{displaymath}
\nabla\beta_{k}(\tilde \Psi(x,\iota))^{\top}\ s\ +\ o\ \left(||\tilde \Psi(x,t_{\iota(1,s)(j)})-\tilde \Psi(x,\iota)||^{ 2}\right)
\end{displaymath}
\begin{displaymath}
=\nabla\beta_{k}(\tilde \Psi(x,\iota))^{\top}\ \left(s -\frac{\tilde \Psi(x,t_{\iota(1,s)(j)})-
\tilde \Psi(x,\iota)}{||\tilde \Psi(x,t_{\iota(1,s)(j)})-\tilde \Psi(x,\iota)||}\right),
\end{displaymath}
also 
\begin{equation}
\nabla\beta_{k}\ (\tilde \Psi(x,\iota))^{\top}\ \Theta^{\star}_{1} = 0
\end{equation}
gilt. Daher ist dann der Rang des Wertes der Jacobi-Matrix $\frac{\partial\beta}{\partial x}$
an der Stelle $\tilde{x}(\iota)$
\begin{equation}
\mbox{Rg}\ \left(\frac{\partial\beta}{\partial x}(\tilde{x}(\iota))\right) = n-\dim \Theta^{\star}_{1}\  
\end{equation}
und 
\begin{equation}
\dim \Theta^{\star}_{1} = \dim \beta^{-1}(\{\omega\})\ .
\end{equation}
Wenn also die Trajektorie $\tilde \Psi(x,\mathbb{R})$ abgeschlossen und quasiperiodisch, aber nicht geschlossen ist,
ist diejenige glatte minimale invariante Mannigfaltigkeit $\beta^{-1}(\{\omega\})$ aus $\mathcal{M}^{ 1}([\Psi])$,
in der die Trajektorie $\tilde \Psi(x,\mathbb{R})$ liegt, mindestens zweidimensional und nicht mit $\tilde \Psi(x,\mathbb{R})$ identisch.
Wenn die Trajektorie $\tilde \Psi(x,\mathbb{R})$ nicht abgeschlossen ist,
kann sie mit keiner minimalen invarianten Mannigfaltigkeit identisch sein.\newline {\bf q.e.d.}
\newline
\newline
Zu der 
von diesem Satz
behaupteten
Nicht-Identit"at von minimalen invarianten Mannigfaltigkeiten
mit Trajektorien 
im Sinne der
konzeptionellen Differenz, die zwischen den 
Trajektorien der Menge $[\Psi]$
und den 
glatten minimalen invarianten Mannigfaltigkeiten der Menge $\mathcal{M}^{ 1}([\Psi])$
besteht,
gibt es insofern ein Gegen"uber:\newline
Es gibt n"amlich die in der folgenden Definition A.2 konzipierte Stetigkeit der allgemeinen Flussfunktionen $\Psi$ binnen 
deren jeweiligem Zustandsraum $\mathbf{P}_{2}\Psi$.
Diese Form der Stetigkeit erweist sich als so beschaffen,
dass genau dann, wenn die Restriktion
$\Psi|(\chi\times \mathbb{R})$
die besagte Form von Stetigkeit hat,
es auch wahr ist,
dass jede Trajektorie 
$$\Psi(x,\mathbb{R})\in [\Psi]:=\{\Psi(x,\mathbb{R}):x\in\mathbf{P}_{2}\Psi\} $$
geschlossen ist, wenn sie in einem 
Unterzustandsraum $\chi\subset\mathbf{P}_{2}\Psi$
mit
\begin{equation}\label{aplace}
\bigcup\Psi(\chi,\mathbb{R})=\bigcup\{\Psi(z,\mathbb{R}):z\in\chi\}=\chi
\end{equation}
enthalten ist.
Als Flussfunktionen im weiteren Sinn binnen einer beliebigen Menge $\mathbf{P}_{2}\Psi$ bezeichnen
wir dabei ja exakt diejenigen Wellenfunktionen gem"ass (\ref{viiergo}), die 
(\ref{viergo}) erf"ullen.\index{Wellenfunktion}
Diese allgemeinen Flussfunktionen modellieren den
Determinismus, der "uber die klassische Punktteilchenmechanik hinausgeht.\index{Flussfunktion im weiteren Sinn binnen einer Menge}
Die besagte
Form der Stetigkeit
von allgemeinen Flussfunktionen im weiteren Sinn binnen deren jeweiligem Zustandsraum $\mathbf{P}_{2}\Psi$, 
die wir nach Pierre Simon de Laplace benennen,
setzt voraus, dass auf dem Zustandsraum $\mathbf{P}_{2}\Psi$ eine Metrik $d$ eingerichtet ist:
\newline
\newline
{\bf Definition A.2: 
Die Laplace-Stetigkeit einer Flussfunktion}\index{Laplace-Stetigkeit einer Flussfunktion}
\newline
{\em Es sei} $\Psi$ {\em eine Flussfunktion, 
deren Zustandsraum $\mathbf{P}_{2}\Psi$ durch die 
Metrik $d$
metrisiert sei.
Jede Restriktion $\Psi|(\chi\times \mathbb{R})$ einer solchen Flussfunktion $\Psi$ nennen wir
genau dann Laplace-stetig bez"uglich der Metrik $d$, wenn erstens
$\chi\subset\mathbf{P}_{2}\Psi$ ein Unterzustandsraum ist, f"ur den} (\ref{aplace}) {\em gilt,
und falls zweitens die
Aussage
\begin{equation}\label{bplace}
\begin{array}{c}
\forall\ (z,\delta)\in \chi\times\mathbb{R}^{+} \ \exists\varepsilon(z,\delta) \in \mathbb{R}^{+}\ \forall\ 
(x,t)\in \chi\times\mathbb{R} \\
d(z,x)<\varepsilon(z,\delta)\Rightarrow d(\Psi(z,t),\Psi(x,t))<\delta
\end{array}
\end{equation}
wahr ist.}
\newline
\newline
Genau dann, wenn $\Psi|(\chi\times \mathbb{R})$ Laplace-stetig bez"uglich der Metrik $d$ ist, sagen wir
auch, dass $\Psi$
auf 
dem Unterzustandsraum 
$\chi\subset\mathbf{P}_{2}\Psi$ bez"uglich der Metrik $d$
Laplace-stetig ist. Als
Laplace-stetige Trajektorien bez"uglich einer Metrik $d$
bezeichnen wir dabei exakt die 
Mengen $\Psi(x,\mathbb{R})$ f"ur Zust"ande $x\in\mathbf{P}_{2}\Psi$, die in einem 
Unterzustandsraum liegen, auf dem $\Psi$
bez"uglich $d$
Laplace-stetig ist.\index{Laplace-stetige Trajektorie} Und wenn wir 
einfach von 
einer Flussfunktion $\Psi$ reden, die 
Laplace-stetig bez"uglich einer Metrik ihres Zustandsraumes $\mathbf{P}_{2}\Psi$ sei,
so meinen wir exakt dies,
dass diese Flussfunktion $\Psi$ auf dem gesamten Zustandsraum $\mathbf{P}_{2}\Psi$ Laplace-stetig bez"uglich seiner Metrik ist.
\newline
{\footnotesize Dass wir dabei die Laplace-Stetigkeit einer Flussfunktion
nach Pierre Simon de Laplace benennen, hat folgenden Grund:
Die Laplace-Stetigkeit einer Flussfunktion gibt die Stetigkeitsauffassung
des Determinismus
wieder, der der 
Kinematik 
mechanischer Punktteilchensysteme unkritisch unterstellt worden war,
bis Poincar\`e und 
Einstein\index{Einstein, Albert}\index{Poincar\`e, Henri}
dar"uber erstaunten, dass innerhalb der klassischen 
Kinematik 
mechanischer Punktteilchensysteme
sehr wohl ein Entwicklungsverhalten 
auftreten kann, das eben nicht 
durch Flussfunktionen beschrieben werden kann, die 
Laplace-stetig sind. Zur Historie siehe beispielsweise \cite{entd}.\newline
Gerade Laplace darf dabei keineswegs als
Propagator der explizit formulierten These angesehen werden, dass
alle klassischen Kinematiken 
mechanischer Punktteilchensysteme
durch Flussfunktionen $\Psi$ beschrieben werden k"onnen,
f"ur die die Laplace-Stetigkeit gem"ass
(\ref{bplace})
auf allen ihren jeweiligen Unterzustandsr"aumen 
$\chi\subset\mathbf{P}_{2}\Psi$
gegeben ist.
Pierre Simon de Laplace hat diese Behauptung nie aufgestellt. 
H"atte er die Laplace-Stetigkeit gem"ass
(\ref{bplace}) explizit formulieren k"onnen, so
w"are sie ihm damit sogleich fragw"urdig geworden. Dieser Irrealis
redet dabei von 
der M"oglichkeit, dass
Laplace die Laplace-Stetigkeit 
in der modernen, post-Cantorschen Form
gem"ass
(\ref{bplace}) explizit formulieren kann.
Diese M"oglichkeit bestand aber f"ur Laplace damals nicht.
\newline
Wie alle mathematischen Physiker des Intervalles von Newton
bis Poincar\`e und Einstein hat sich aber auch Laplace
so verhalten, als ob die
Laplace-Stetigkeit
eine nat"urliche Gegebenheit sei.
Und er ist dabei neben Euler und Hamilton einer der herausragendsten Arbeiter unter
jener unformulierten
Annahme, dass
die Laplace-Stetigkeit innerhalb der klassischen 
Kinematik 
mechanischer Punktteilchensysteme
vorliege.
\newline
Dabei ist es Laplace, der ausserdem
den Determinismus erstmalig durchaus problematisiert, jedoch 
nicht im Hinblick auf seine topologische Beschaffenheit an sich:
Im Zusammenhang mit dem Thema des Determinismus
bringt selbst der Laie den Determinismus
mit dem von Laplace imaginierten 
Laplaceschen D"amon in Verbindung, der als perfekter Historiker
alles
zur"uckberechnen kann und der als Wahrsager alles 
vorausberechnen kann, falls ihm die sowohl ihrem Umfang 
als auch ihrer Qualtit"at nach
totale einmalige Kenntnis
des
Status Quo 
der Welt
gegeben ist.
Die Totalit"at der Qualtit"at seiner 
Kenntnis ist hierbei, dass ihm die Bestandsaufnahme mit idealer Pr"azision aller Messungen 
vorliegt.
\newline
Laplace denkt sich dabei den 
Determinismus noch als den Determinismus gem"ass
jener
d"ammerigen
Unterstellung, dass irgendwie
die Laplace-Stetigkeit innerhalb der klassischen 
Kinematik 
mechanischer Punktteilchensysteme
gilt. 
\newline
Damit, dass Laplace den nach ihm 
benannten D"amon imaginiert und damit, dass Laplace denselben
der Wissenschaftsgeschichte der folgenden Jahrhunderte
als die Flaschenpost, die jede einpr"agsame Denkw"urdigkeit ja ist,
aufgibt,
inspiriert er zur Hinterfragung des Determinismus.
Und damit regt Laplace
auch dazu an, denselben auf seine Beschaffenheit an sich hin
zu bedenken. Dadurch beispielsweise, dass 
Pierre Simon de Laplace den Kenntnisstand 
seines D"amons als einen
Kenntnisstand anlegt,
der nach einer Bestandsaufnahme mit idealer Pr"azision der Messungen 
vorliegt, 
leitet er auf nat"urliche Weise zu der Frage, was denn aber gelten w"urde, wenn selbst diesem D"amon
die Idealit"at der Bestandsaufnahme nicht realisierbar w"are.\newline
Gilt dann, wenn 
nur eine ann"aherungsweise 
Bestandsaufnahme m"oglich ist,
dass eine ann"aherungsweise R"uckw"arts- und Vorw"artsentwicklung bestimmt werden kann?
Und exakt diese Frage ist die Frage, die 
Poincar\`e und Einstein schliesslich negieren, indem sie
die G"ultigkeit der Laplace-Stetigkeit innerhalb der klassischen 
Kinematik 
mechanischer Punktteilchensysteme ausdr"ucklich und durchaus elektrisiert
verneinen.\newline
So gesehen eignet sich wohl niemand besser, der  
Stetigkeit 
gem"ass (\ref{bplace})
den Namen zu geben, als Pierre Simon de Laplace,\index{Laplace, Pierre Simon de} 
die man ansonsten alternativ als die 
Pr"a-Poincar\`e-Einstein-Stetigkeit 
bezeichnen m"usste.}\newline
Der besagte
Korrespondenzpartner  
der 
Behauptung des
Satzes A.1 
ist die nun 
folgende Aussage "uber  
Laplace-stetige Flussfunktionen und deren 
Trajektorien: 
\newline
\newline
{\bf Bemerkung A.3:\newline Die topologische Abgeschlossenheit Laplace-stetiger Trajektorien}
\newline
{\em Es sei} $\Psi$ {\em eine Laplace-stetige Flussfunktion, 
deren Zustandsraum $\mathbf{P}_{2}\Psi$ durch die 
Metrik $d$
metrisiert sei, welche 
die Topologie $\mathbf{T}(d)$
induziere. Dabei sei die topologisierte Flussfunktion $(\Psi,\mathbf{T}(d))$ so beschaffen,
dass die zu} (\ref{aalzb}){\em analoge Form zustandsweiser Sensitivit"at f"ur und nur f"ur deren nicht
triviale Zimmer 
$$\chi\in \{\mathbf{cl}_{\mathbf{T}(d)}(\Psi(z,\mathbb{R})):z\in\mathbf{P}_{2}\Psi \}\setminus
\{\Psi(z,\mathbb{R}):z\in\mathbf{P}_{2}\Psi \}$$
gegeben sei.
Dann gilt die Identit"at
\begin{equation}\label{cplace} 
\Psi(z,\mathbb{R})=\mathbf{cl}_{\mathbf{T}(d)}(\Psi(z,\mathbb{R}))
\end{equation}
f"ur alle Zust"ande $z\in \mathbf{P}_{2}\Psi$ des Zustandsraumes, wobei $\mathbf{cl}_{\mathbf{T}(d)}$ der
auf die induzierte Topologie $\mathbf{T}(d)$ bezogene H"ullenoperator sei.}
\newline
\newline
{\bf Beweis:}\newline
Es sei
$$x\in\mathbf{cl}_{\mathbf{T}(d)}(\Psi(z,\mathbb{R}))\ ,$$
sodass eine monoton wachsende Folge
$$\{t_{j}\}_{j\in\mathbb{N}}\in\ (\mathbb{R}^{+})^{\mathbb{N}}$$
und die Zustandsfolge
$$\{z_{j}\}_{j\in\mathbb{N}}:=\{\Psi(z, t_{j})\}_{j\in\mathbb{N}}\in(\mathbf{P}_{2}\Psi)^{\mathbb{N}}$$
existiert, f"ur die
$$\lim_{j\to\infty}d(z_{j},x)=0$$
ist.
Wegen der Laplace-Stetigkeit der Flussfunktion $\Psi$ gilt
die Aussage
$$\forall\delta\in \mathbb{R}^{+}\ \exists\ j(\delta)\in \mathbb{N}:\ 
\forall\ j\in]j(\delta),\infty[\cap \mathbb{N},t\in \mathbb{R}$$
$$d(\Psi(z_{j},t),\Psi(x,t))<\delta\ .$$
Daher ist es ausgeschlossen, dass die zu (\ref{aalzb}) analoge Form zustandsweiser Sensitivit"at 
auf $\mathbf{cl}_{\mathbf{T}(d)}(\Psi(z,\mathbb{R}))$
vorliegt. Voraussetzungsgem"ass ist
$x$ also ein Zustand des trivialen Zimmers $\mathbf{cl}_{\mathbf{T}(d)}(\Psi(z,\mathbb{R}))=\Psi(z,\mathbb{R})$.
F"ur alle $j\in \mathbb{N}$ ist also
$$\Psi(x,\mathbb{R})=\Psi(z_{j},\mathbb{R})=\Psi(z,\mathbb{R})$$
und
$$\mathbf{cl}_{\mathbf{T}(d)}(\Psi(z,\mathbb{R}))\setminus\Psi(z,\mathbb{R})=\emptyset\ .$$
\newline {\bf q.e.d.}
\newline
\newline
Wie die Zyklizit"at der Trajektorie $\Psi(z,\mathbb{R})=\mathbf{cl}_{\mathbf{T}(d)}(\Psi(z,\mathbb{R}))$
bez"uglich der metrischen Topologie $\mathbf{T}(d)$ die Periodizit"at aller ihrer Zust"ande $y\in\Psi(z,\mathbb{R})$
bez"uglich der stetigen Flussfunktion $\Psi$ impliziert, haben wir bereits im ersten Kapitel 
beispielhaft f"ur 
den Fall er"ortert, dass $\mathbf{T}(d)$
eine nat"urliche Topologie eines endlichdimensionalen reellen Raumes ist. Offenbar ist jene Betrachtung 
analog anstellbar.  
\newline
Es ist bemerkenswert, dass das f"urwahr einfache
Argument des Beweises der Bemerkung A.3 blosslegt, dass im Falle der Laplace-Stetigkeit alle Trajektorien
abgeschlossen sind und damit periodische Zustandsentwicklungen beschreiben.
Freilich, vorformale Vorstellungen sind 
psychologische Ph"anomene und keine logischen Objekte. Sie liegen als
solche jenseits formaler Fassbarkeit und daher kann auch die
Laplace-Stetigkeit nicht ins Dunkel der Kontinuit"atsvorstellungen 
"uber den Determinismus
der Laplaceschen Epoche zur"uck.
Um so bemerkenswerter erscheint 
nichtsdestotrotz
dieser Sachverhalt, wenn wir allem kritischen Streuben zum Trotz einmal so tun, als k"onnte
jene starke Form der Stetigkeit einer Flussfunktion, die in der Definition {\bf A.2} als deren Laplace-Stetigkeit 
formuliert ist, dennoch die dunklen Kontinuit"atsvorstellungen 
"uber den Determinismus
der Laplaceschen Epoche wiedergeben:
Dann h"atte diese Kontinuit"atsvorstellung die Periodizit"at aller Zustandsentwicklungen impliziert. H"atte Laplace diese Implikation
universeller Periodizit"at geglaubt?
Gewiss, Laplacens 
Blick auf das Planetensystem, 
das man damals noch leicht f"ur vollendet periodisch halten konnte,
markiert einen Interessenfokus gerade 
dieses Autors des f"unfb"andigen $Trait\acute{e}$
$de$ $m\acute{e}chanique$ $c\acute{e}leste$,
das als Laplacens Hauptwerk gilt. 
Das Interesse an der Himmelsmechanik d"urfen wir dabei aber durchaus
als das allgemeine paradigmatische Interesse Laplacens Zeit
ansehen.
Wird hier ein wenig ahnbar, wie sehr die mathematische Physik vor den formalen Fortschritten
der das Kontinuum erforschenden Mathematik nach Laplace noch im Dunkel lag?
\newline
Wir warfen in diesem Anhang noch einen Blick zur"uck: 
N"amlich auf die Nostalgika des Determinismus. Erstens auf das
Nostalgikum, sich die abgeschlossenen H"ullen jeweiliger Trajektorien 
nicht anders
als die jeweiligen Trajektorien selbst zu denken.\footnote{Sich die abgeschlossenen H"ullen jeweiliger Trajektorien
anders
als die jeweiligen Trajektorien selbst {\em vorzustellen}, ist zumindest
eine gewisse Strapaze f"ur die Imagination.}
Und zweitens auf das 
Nostalgikum der den Determinismus betreffenden Kontinuit"atsvorstellung,
die der 
Kinematik 
mechanischer Punktteilchensysteme 
unformalisiert und
unkritisch unterstellt worden war und 
die die Laplace-Stetigkeit formuliert.\newline
Dieser Blick sollte
uns die n"otige
Distanzierung von den 
in uns immer noch verfestigten nostalgischen 
Auffassungen vermitteln.
Nach diesem distanzierenden Blick 
auf die Nostalgika des Determinismus
sind wir nun ger"ustet,
die Aussage des elementaren Quasiergodensatzes auch imaginativ zu verinnerlichen, die
die Beantwortung der objektivierten Quasiergodenfrage im Sinne P. und T. Ehrenfests ist -- und insofern die
Best"atigung 
der minimalinvasiven Abwandlung von Boltzmanns 
unhaltbarer
Ergodenhypothese, deren tragische Genialit"at uns aufgehe.\index{Ergodenhypothese}
\newpage
{\Large {\bf Symbolverzeichnis}}\newline\newline
Neben den v"ollig etablierten Bezeichnungen benutzen wir in diesem Traktat 
auch weniger gel"aufige oder spezielle, f"ur dieses Traktat idiomatische Notationen. 
Diese weniger gel"aufigen Schreibweisen f"uhren wir
im fortlaufenden
Text ein und die Stellen, an denen diese Schreibweisen vorgestellt werden, sind im Index unter dem Stichwort \glqq Notationskonvention\grqq  vermerkt.\newline
Es sei $n,m\in\mathbb{N}$, $r\in \mathbb{R}^{+}$, $t\in \mathbb{R}$, $x\in\mathbb{R}^{n}$, ${\rm A}$ eine Menge, ${\rm E}=\{e\}$
eine einelementige Menge, $j\in\{1,2,\dots n\}$,
$\Theta=(\theta_{1},\theta_{2},\dots \theta_{n})$
ein $n$-Tupel, 
$\Lambda=(\lambda_{1},\lambda_{2},\dots \lambda_{m})$ ein $m$-Tupel,
$\Psi$ eine Flussfunktion, $\phi$ eine Funktion, $q\in\{-1,-1/2,0,1/2,1,\dots\}$
ein
Kontinuit"atsindex und es sei
$\alpha$ ein Mengensystem: Es ist dann
\newline
\begin{tabbing}
Wir benutzen zudem folgende \= Etwas \kill
$\Psi^{t}$ \> der durch $\Psi$ und $t$ festgelegte Phasenfluss,\\   
\quad \> d.h., die Bijektion $\Psi(\mbox{id},t)$\\ 
\quad \>  des Zustandsraumes auf sich,\\ 
$[\Psi]$  \> die durch $\Psi$ festgelegte trajektorielle Partition,\\ 
\quad \> d.h., $[\Psi]:=\{\Psi(x,\mathbb{R}):x\in\mathbf{P}_{2}\Psi\}\ ,$\\ 
$[[\Psi]]$ \> die durch $\Psi$ festgelegte nat"urliche Partition, \\ 
\quad \> d.h., $[[\Psi]]:=\{\mathbf{cl}(\tau):\tau\in[\Psi]\}\ ,$\\ 
${\rm Inv}^{q}([\Psi])$ \> die Menge Invarianter der Kontinuit"at $q$\\ 
\quad \> der trajektoriellen Partition $[\Psi]$ ,\\
$\mathcal{M}^{q}([\Psi])$ \> die Menge deren minimaler invarianter\\ 
\quad \> Mannigfaltigkeiten der Kontinuit"at $q$,\\ 
$\mathbf{part}({\rm A})$  \> die Menge aller Partitionen der Menge ${\rm A}$,\\
$\alpha^{\cup}$ \> das Mengensystem aller Vereinigungen\\ \quad \> "uber Teilmengen von $\alpha$,\\
$\alpha_{{\rm A}}$ \> die ${\rm A}$-Auswahl aus $\alpha$, d.h.,\\
\quad \> $\{a\in \alpha:a\cap{\rm A}\not=\emptyset\}\ ,$\\
$\downarrow {\rm E} =e$, \> das Element der einelementige Menge ${\rm E}$, \\
$\mathbf{P}_{j}\Lambda=\lambda_{j}$ , \> \quad\\
$\mathbf{P}_{1}\phi$ \> die Definitionsmenge von $\phi$,\\
$\mathbf{P}_{2}\phi$ \>  die Wertemenge von $\phi$,\\
${\rm A}^{\mathbb{N}}$ \> die Menge aller Folgen $\{a_{j}\}_{j\in\mathbb{N}}$,\\ 
\quad \> deren Glieder Elemente der Menge ${\rm A}$ sind,\\
$\mathbb{B}_{r}(x)$ \> die offene Kugel des $\mathbb{R}^{n}$, in dem $x$ ist, deren\\  
\quad \>  Mittelpunkt $x$ ist und
die den Radius $r$ hat,\\
$\sim$ \> die Hom"oomorphie zweier Teilmengen des $\mathbb{R}^{n}$,\\
$\sim\sim$ \> die relative Hom"oomorphie zweier Teilmengen\\ 
\quad \> des $\mathbb{R}^{n}$ oder des $\mathbb{R}^{m}$,\\
$\oplus$ \> die Kokatenation,\\ \quad \> sodass $\Theta\oplus\Lambda=(\theta_{1},\theta_{2},\dots \theta_{n},
\lambda_{1},\lambda_{2},\dots \lambda_{m})$ ist,\\
$\olessthan$ \quad \>das Kroneckerprodukt,\\
\end{tabbing}
das so erl"autert sei:
Ist
\begin{displaymath}
A=(a_{jk})_{1\leq j\leq n,1\leq k\leq m} = \left( \begin{array}{cccc}
a_{11} & a_{12} &\dots & a_{1m}\\
\dots &\quad &\quad &\quad\\
a_{n1} & a_{n2} &\dots & a_{nm}  \end{array} \right)
\end{displaymath} 
eine $m\times n$-Matrix $A\in {\rm K}^{m\times n}$ "uber dem K"orper $({\rm K},+,\cdot)$
und 
\begin{displaymath}
B=(b_{jk})_{1\leq j\leq u,1\leq k\leq v} = \left( \begin{array}{cccc}
b_{11} & b_{12} &\dots & b_{1v}\\
\dots &\quad &\quad &\quad\\
b_{u1} & b_{u2} &\dots & b_{uv}  \end{array} \right)\ ,
\end{displaymath}
eine
$u\times v$-Matrix $B\in {\rm K}^{u\times v}$ "uber dem selben K"orper $({\rm K},+,\cdot)$,
so ist diejenige $mu\times nv$-Matrix "uber dem selben K"orper $({\rm K},+,\cdot)$, deren
Matrixelemente 
f"ur alle  $$(j,k,\alpha,\beta)\in\{1,2\dots n\}\times\{1,2\dots m\}\times\{1,2\dots u\}\times\{1,2\dots v\}$$
zwischen der $j(\alpha-1)$-ten und der $j\alpha+1$-ten Zeile und der 
$k(\beta-1)$-ten und der $k\beta+1$-ten Spalte auch hinsichtlich ihrer Anordnung mit den 
Matrixelementen der Matrix
\begin{displaymath}
a_{jk}B=a_{jk}(b_{rs})_{1\leq r\leq u,1\leq s\leq v} = \left( \begin{array}{cccc}
a_{jk}b_{11} & a_{jk}b_{12} &\dots & a_{jk}b_{1v}\\
\dots &\quad &\quad &\quad\\
a_{jk}b_{u1} & a_{jk}b_{u2} &\dots & a_{jk}b_{uv}  \end{array} \right)
\end{displaymath}
identisch sind, das Kroneckerprodukt\index{Kroneckerprodukt} $$A \olessthan B$$ der $m\times n$-Matrix $A$ und der $u\times v$-Matrix $B$. Dadurch
ist die Verkn"upfung $\olessthan$, die
offenbar im Allgemeinen eklatant nicht kommutiert, 
f"ur alle $m,n,u,v\in\mathbb{N}$
auf jeder Menge ${\rm K}^{m\times n}\times{\rm K}^{u\times v}$
definiert.

\renewcommand{\indexname}{Begriffsregister}
\printindex
\end{document}